\documentclass[preprint,12pt]{elsarticle}

\journal{Computer Methods in Applied Mechanics and Engineering}

\usepackage[utf8]{inputenc}
\usepackage{natbib}
\usepackage[centerlast]{subfigure}

\usepackage{amssymb}
\usepackage{amsmath,mathrsfs}
\usepackage{math}
\usepackage{color}
\usepackage{bbm}
\usepackage{float}
\usepackage{graphicx}
\usepackage{siunitx}
\usepackage{stfloats}
\usepackage{bm}
\usepackage{caption}
\usepackage{subcaption}
\usepackage{multicol}
\usepackage{array}
\usepackage{booktabs} 
\usepackage{multirow}
\usepackage{listings}

\usepackage{chngcntr}
\usepackage{icomma}
\DeclareMathSymbol{,}{\mathord}{letters}{"3B}

\usepackage{algorithm}
\usepackage{algpseudocode}

\usepackage[toc,page]{appendix}

\usepackage{comment}

\usepackage[usenames,dvipsnames,svgnames,table]{xcolor}
\usepackage[version=4]{mhchem}

\usepackage{tikz}
\usetikzlibrary{patterns}
\usetikzlibrary{external}
\usetikzlibrary{decorations.pathreplacing}
\usetikzlibrary{shapes.misc}
\usetikzlibrary{backgrounds}

\usepackage{pgfplots}
\pgfplotsset{compat=newest}
\usepgfplotslibrary{groupplots}
\usepgfplotslibrary{fillbetween}
\usepgfplotslibrary{colormaps}

\tikzset{
    ultra thin/.style= {line width=0.1pt},
    very thin/.style=  {line width=0.2pt},
    thin/.style=       {line width=0.4pt},
    semithick/.style=  {line width=0.6pt},
    thick/.style=      {line width=1.0pt},
    very thick/.style= {line width=1.2pt},
    ultra thick/.style={line width=1.6pt}
}

\usepackage{wrapfig}

\usepackage{sectsty}  
\usepackage[colorlinks,linkcolor=blue,anchorcolor=green,citecolor=blue]{hyperref}
\usepackage[capitalise]{cleveref}
\renewcommand{\eqref}[1]{Eq.~$($\ref{#1}$)$}  

\biboptions{sort&compress} 

\usepackage{xcolor}
\definecolor{green}{RGB}{87,171,39}
\definecolor{blue}{RGB}{0,84,159}
\definecolor{red}{RGB}{204,7,30}

\begin{document}

\begin{frontmatter}
    \title{A robust finite strain isogeometric solid-beam element}
    \author[els]{Abdullah Shafqat\corref{cor1}}
    \ead{abdullah.shafqat@tu-darmstadt.de}
    \author[rtt]{Oliver Weeger}
    \author[els]{Bai-Xiang Xu}

    \address[els]{Mechanics of Functional Materials Division, Institute of Materials Science, Technical University of Darmstadt, Otto-Berndt-Str.~3, 64287 Darmstadt, Germany}

    \address[rtt]{Cyber-Physical Simulation, Department of Mechanical Engineering, Technical University of Darmstadt, Dolivostr.~15, 64293 Darmstadt, Germany}
    \cortext[cor1]{Corresponding author}

    \begin{abstract}
        
        In this work, an efficient and robust isogeometric three-dimensional solid-beam finite element is developed for large deformations and finite rotations with merely displacements as degrees of freedom. The finite strain theory and hyperelastic constitutive models are considered and B-Spline and NURBS are employed for the finite element discretization. Similar to finite elements based on Lagrange polynomials, also NURBS-based formulations are affected by the non-physical phenomena of locking, which constrains the field variables and negatively impacts the solution accuracy and deteriorates convergence behavior. To avoid this problem within the context of a Solid-Beam formulation, the Assumed Natural Strain (ANS) method is applied to alleviate membrane and transversal shear locking and the Enhanced Assumed Strain (EAS) method against Poisson thickness locking. Furthermore, the Mixed Integration Point (MIP) method is employed to make the formulation more efficient and robust. The proposed novel isogeometric solid-beam element is tested on several single-patch and multi-patch benchmark problems, and it is validated against classical solid finite elements and isoparametric solid-beam elements. The results show that the proposed formulation can alleviate the locking effects and significantly improve the performance of the isogeometric solid-beam element. With the developed element, efficient and accurate predictions of mechanical properties of lattice-based structured materials can be achieved. The proposed solid-beam element inherits both the merits of solid elements e.g. flexible boundary conditions and of the beam elements i.e. higher computational efficiency.
    \end{abstract}
    
\end{frontmatter}

\section{Introduction}
\label{Sec:Introduction}

Beams are the structural components of mechanical systems, with applications ranging from the meter-scale such as in the frame structures of construction buildings,  to the micro- or even nanometer scale in fibrous materials. 
In particular, due to the recent development of advanced and additive manufacturing technologies, lattice structures and periodic meta-materials have gained immense attention \cite{Greer.2019}. 
Beyond mechanical meta-materials, additive manufacturing techniques can also be used to efficiently fabricate complex, micro-structured energy devices, such as batteries or solar cells, due to good geometry controllability, as well as cost and material waste reduction \cite{Pang.2020,Narita.2022}. The building blocks of such lattice structures are often truss- or beam-like struts. Hence, accurate simulation of such strut behavior is important in order to recapture the mechanical and multi-functional behaviour of such lattice-structured materials \cite{Greer.2019,jamshidian2020,zhang2022a,weeger2022}. 

Beams with a high length-to-thickness aspect ratio can be efficiently modeled by  elementary beam theories \cite{antman2005,eugster2015} and analysed using beam finite elements and similar numerical techniques \cite{Bathe.1979,Simo.1986,jung2011,meier2014,WEEGER2017100}.
However, beam theories, which generally abstract the three-dimensional geometry to the beam axis, assume undeformed cross-sections, and employ reduced kinematics and averaged stress measures, have their limitations. 
Constitutive models have to use the reduced stress state based on the assumption that the normal transversal stresses are much smaller than the axial stress and transversal shear stresses, which induces errors when hyperelastic material models are to be used at large strains. 
These assumptions and simplifications further pose challenges when inelastic material behaviors, such as viscoelasticity or plasticity \cite{simo1984,gruttmann2000b,lestringant2020,herrnbock2022,Weeger.2022b}, or multiphysical couplings such as chemo-, electro-, or thermo-elasticity are to be considered \cite{favata2016,ebrahimi2016,smriti2019}.
For instance, Xia \textit{et al.} \cite{Xia.2019} show the sinusoidal behaviour of a three-dimensional tetragonal micro-lattice electrode due to the in-plane buckling of the beams when lithiated. During the lithiation, the beams are strongly driven from the both-sided contact, due to which they also expand radially and longitudinally as explained in \cite{Stein.2014,Zhao.2015}. 
To capture the radial expansion behaviour of beam-like electrodes during lithiation, the precise discretization of the cross-sectional dimensions is also needed, which cannot be included in conventional beam theories. Likewise, in lattice structures, the lithium concentration can be induced along the beam axis or through the cross-section, which is again challenging for beam models.

While beam models possess several challenges regarding the modeling of slender structures with complex, inelastic, and multi-physical behaviours, solid finite elements are unable to capture the behaviour of slender structures such as beams (and shells) due to the presence of spurious numerical effects, i.e., \emph{locking}. Generally, locking leads to an overestimation of the stiffness and thus an underestimation of field variables, which is highly undesirable for engineering analysis. 
To circumvent the limitations of beam element and cure the numerical issues of solid elements, \emph{solid-beam elements} have been developed. Similar to solid elements, they rely on 3D geometries and employ only displacement degrees of freedom, but due to the high aspect ratios, they must apply techniques to alleviate locking. 
Throughout the years, several approaches have been developed to cure locking in Lagrange-based isoparametric elements, such as a mixed displacement-pressure formulation \cite{SUSSMAN1987357}, reduced and selective reduced integration \cite{Zienkiewicz1971ReducedIT,MALKUS197863}, the ${\bar{B}}$ strain projection method \cite{BbarMethod}, the assumed natural strain (ANS) method \cite{Bathe1985,Simo1986}, and the enhanced assumed strain (EAS) method \cite{Simo.1990, Simo1992}.

Isogeometric analysis (IGA) \cite{Hughes.2005} provides the opportunity to unify the computer-aided design (CAD) geometry and numerical discretization by using B-spline or NURBS functions in the finite element framework. In IGA, the same basis functions are used to model the CAD geometry and the approximation of the unknown fields, leading to an exact geometry representation regardless of mesh density. Furthermore, compared to standard FEM, high-order and smooth basis functions can be easily used to improve the accuracy and robustness of the numerical analysis \cite{Cottrell.2007}. Further  benefits of IGA can be also found in shape optimization problems \cite{Wall.2008,Nagy.2010}, since IGA allows direct utilization of control points from the spline geometry as design variables for optimization. 

However, as discussed by Echter \textit{et al.} \cite{Echter.2010}, NURBS-based element formulations are also affected by numerical locking. IGA formulations of slender structures also suffer from transversal shear and membrane locking, which occur due to the inability of the interpolation functions to replicate ``shearless bending'' and ``inextensible bending'', respectively.
In the last two decades, a significant amount of work was done in IGA element formulations to remove the locking phenomenon. Valuable contributions in this regard can be found for shear-deformable Timoshenko beams \cite{Echter.2010,DaBeiraoVeiga.2012,Adam.2014,Kiendl.2015,Weeger.2022b}, curved beams \cite{Bouclier.2012,Vo.2020}, spatial rods \cite{Auricchio.2013,WEEGER2017100}, shells \cite{Elguedj.2008,Echter.2013b,MI2021113693} and solid-shell elements \cite{Bouclier.2013b,Bouclier.2013,Hosseini.2013,Hosseini.2014,caseiro2014assumed,caseiro2015assumed,Bouclier.2015,LEONETTI2018159,Antolin.2020}. 

This work aims to develop a robust solid-beam finite element formulation, which combines the advantages of IGA with techniques to alleviate locking.
Frischkorn and Reese developed an isoparametric solid-beam element \cite{Frischkorn.2013} based on the solid-shell element of Schwarze and Reese \cite{Schwarze.2011b}. The relevant locking in the Lagrange-based solid-beam element is removed by the combination of the assumed natural strain method (ANS) and the enhanced assumed strain method (EAS), which are embedded in a reduced integration concept. 
Recently, Choi \textit{et al.} \cite{Choi.2021c,Choi.2023} developed the nonlinear isogeometric beam formulations with extensible directors with only displacement degrees of freedom. In \cite{Choi.2023}, locking is alleviated by a mixed formulation with the EAS method. 
Furthermore, Magisano \textit{et al.} \cite{Magisano.2017} proposed the mixed integration point (MIP) method to enhance the robustness  and efficiency of displacement-based finite element formulations in geometrically non-linear problems. Later, Pfefferkorn \textit{et al.} \cite{Pfefferkorn.2021} extend this idea to make the EAS elements more robust, as the finite elements based on EAS method require a large number of Newton iterations and small load increments. The MIP method was also recently extended to hyperelastic Kirchhoff-Love shells for nonlinear static and dynamic analysis in \cite{Leonetti.2023}.

Employing these aforementioned state-of-the-art techniques, a robust isogeometric solid-beam element is proposed in this work. It uses the ANS method to alleviate transversal shear and membrane locking, the EAS method for removing Poisson thickness locking, and the MIP method to make the formulation robust and efficient. 
The choice of the ANS method is due to its simplicity and efficiency in implementation, even for higher-order NURBS. For the solid-beam element, the idea of using the ANS method is inspired from IGA solid-shell elements \cite{caseiro2014assumed,caseiro2015assumed}, the EAS from the isoparametric solid-beam element \cite{Frischkorn.2013}, and the MIP method for the EAS elements from \cite{Pfefferkorn.2021}. 

The paper is structured in the following way: In \cref{sec:Finite strain elasticity}, the finite strain theory of the solid-beam element, which is based on the 3D continuum theory, is briefly described. \Cref{sec: Isogeometric finite element discretizaton} gives the overview of the basics of IGA and its numerical implementation on the field variables. The particular methods to enable a locking-free and robust isogeometric solid-beam formulation (ANS, EAS, MIP)  are explained in \cref{Sec:Isogeometric solid-beam formulation}, along with details on the numerical discretization of the variational formulation. Then, in \cref{Sec:Numerical Examples}, the performance of the proposed formulation is evaluated on several benchmark problems and comparisons with existing methodologies. In the end, \cref{Sec:Conclusion} concludes the paper with a summary and outlook.   


\section{Finite strain elasticity}
\label{sec:Finite strain elasticity}

\subsection{Kinematics}
\label{sec:Kinematics}

\begin{figure}[t!]
	\centering        \includegraphics[width=0.7\textwidth]{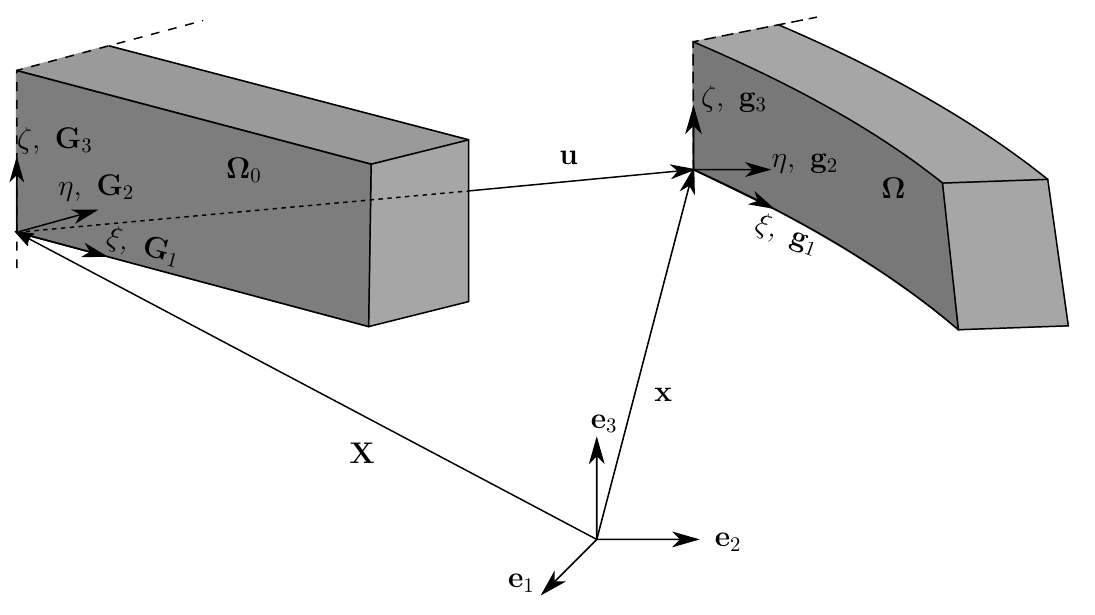}
	\caption{ Reference and current configuration of solid-beam element in Euclidean space}
\label{fig:Kinematics_1}
\end{figure}

In solid-beam formulations, beams are geometrically and kinematically modelled as conventional 3D continua.  
\Cref{fig:Kinematics_1} shows the solid-beam in its reference configuration $\Omega_{0}\subset\mathbb{R}^3$ at time $t={0}$ on the left and in its current configuration $\Omega\subset\mathbb{R}^3$ at time ${t}>0$ on the right.
The position vectors in the reference and in the current configurations of the solid-beam element are represented by $\mathbf{X}\subset\Omega_{0}$ and $\mathbf{x}\subset\Omega$, respectively. The transformation from reference to the current configuration is described by a displacement field $\mathbf{u}:\Omega_{0}\to\mathbb{R}^3$. The position vector in the deformed configuration is calculated as $\mathbf{x}=\mathbf{X}+\mathbf{u}(\mathbf{X})$. 

The mapping of an infinitesimal solid-beam element from the reference to the current configuration can be described by the deformation gradient. The deformation gradient $\mathbf{F}$, its variation $\delta\mathbf{F}$, and its linearization $\Delta\mathbf{F}$ can be written as:
\begin{equation}
\label{eqn:Kinematics.1}
\mathbf{F}(\mathbf{u})=\frac{\partial \mathbf{x}}{\partial \mathbf{X}}=\mathbf{I}+\frac{\partial \mathbf{u}}{\partial \mathbf{X}}, \qquad
\delta\mathbf{F}=\frac{\partial \delta\mathbf{u}}{\partial \mathbf{X}}, \qquad
\Delta\mathbf{F}=\frac{\partial \Delta\mathbf{u}}{\partial \mathbf{X}}.
\end{equation}

In the finite strain theory, objective and symmetric strain measures such as the right Cauchy-Green tensor $\mathbf{C} = \mathbf{F}^{\mathrm{T}}\mathbf{F}$ and the Green-Lagrange strain tensor $\mathbf{E}=\tfrac{1}{2}(\mathbf{C}-\mathbf{I})$ are often employed. The Green-Lagrange strain and its variation $\delta \mathbf{E}$, as well as their linearizations $\Delta \mathbf{E}$ and $ \Delta \delta\mathbf{E}$, are expressed as:
\begin{equation}
\label{eqn:Kinematics.2}
\begin{aligned}
\mathbf{E}(\mathbf{u}) &= \frac{1}{2}(\mathbf{F}^{\mathrm{T}}\mathbf{F}-\mathbf{I}), \qquad & 
\Delta \mathbf{E}(\mathbf{u}) &= \frac{1}{2}(\Delta \mathbf{F}^{\mathrm{T}} \mathbf{F}+\mathbf{F}^{\mathrm{T}} \Delta \mathbf{F}),\\
\delta \mathbf{E}(\mathbf{u}) &= \frac{1}{2}(\delta \mathbf{ F}^{\mathrm{T}} \mathbf{F}+\mathbf{F}^{\mathrm{T}}\delta \mathbf{ F}),
\qquad & 
 \Delta \delta\mathbf{E} &= \frac{1}{2}(\Delta \mathbf{ F}^{\mathrm{T}} \delta \mathbf{ F}+\delta \mathbf{ F}^{\mathrm{T}} \Delta\mathbf{ F}).
\end{aligned}
\end{equation}
As a symmetric second-order tensor, the Green-Lagrange strain can be represented in Voigt notation as:
\begin{equation}
\label{eqn:Kinematics.Ev}
\mathbf{\hat{E}}=\begin{bmatrix}
E_{\xi\xi} & E_{\eta \eta} & E_{\zeta \zeta} & 2E_{\xi \eta} & 2E_{\eta \zeta} & 2E_{\xi \zeta } \\
\end{bmatrix}^\mathrm{T}.
\end{equation}
Using \eqref{eqn:Kinematics.1}, the Voigt form of the variational and linearized Green-Lagrange strain tensors can be written as functions of the position vector and variational and linearized displacement fields, respectively, as:
\begin{equation}
\label{eqn:Kinematics.3}
\delta \hat{\mathbf{E}}=\begin{bmatrix}
\delta E_{\xi\xi} \\[6pt]
\delta E_{\eta \eta} \\[6pt]
\delta E_{\zeta \zeta} \\[6pt]
2\delta E_{\xi \eta} \\[6pt]
2\delta E_{\eta \zeta} \\[6pt]
2\delta E_{\xi \zeta}
\end{bmatrix}=\begin{bmatrix}
\mathbf{x}_{,1}^{\mathrm{T}}\cdot \delta \mathbf{u}_{,1} \\[6pt] 
\mathbf{x}_{,2}^{\mathrm{T}}\cdot \delta \mathbf{u}_{,2} \\[6pt] 
\mathbf{x}_{,3}^{\mathrm{T}}\cdot \delta \mathbf{u}_{,3} \\[6pt] 
\mathbf{x}_{,1}^{\mathrm{T}}\cdot \delta \mathbf{u}_{,2} + \mathbf{x}_{,2}^{\mathrm{T}}\cdot \delta \mathbf{u}_{,1} \\[6pt]
\mathbf{x}_{,2}^{\mathrm{T}}\cdot \delta \mathbf{u}_{,3} + \mathbf{x}_{,3}^{\mathrm{T}}\cdot \delta \mathbf{u}_{,2} \\[6pt]
\mathbf{x}_{,1}^{\mathrm{T}}\cdot \delta \mathbf{u}_{,3}+\mathbf{x}_{,3}^{\mathrm{T}}\cdot \delta \mathbf{u}_{,1}
\end{bmatrix}, \quad
\Delta \hat{\mathbf{E}}=\begin{bmatrix}
\Delta E_{\xi\xi} \\[6pt]
\Delta E_{\eta \eta} \\[6pt]
\Delta E_{\zeta \zeta} \\[6pt]
2\Delta E_{\xi \eta} \\[6pt]
2\Delta E_{\eta \zeta} \\[6pt]
2\Delta E_{\xi \zeta}
\end{bmatrix}=\begin{bmatrix}
\mathbf{x}_{,1}^{\mathrm{T}}\cdot \Delta \mathbf{u}_{,1} \\[6pt] 
\mathbf{x}_{,2}^{\mathrm{T}}\cdot \Delta \mathbf{u}_{,2} \\[6pt] 
\mathbf{x}_{,3}^{\mathrm{T}}\cdot \Delta \mathbf{u}_{,3} \\[6pt] 
\mathbf{x}_{,1}^{\mathrm{T}}\cdot \Delta \mathbf{u}_{,2} + \mathbf{x}_{,2}^{\mathrm{T}}\cdot \Delta \mathbf{u}_{,1} \\[6pt]
\mathbf{x}_{,2}^{\mathrm{T}}\cdot \Delta \mathbf{u}_{,3} + \mathbf{x}_{,3}^{\mathrm{T}}\cdot \Delta \mathbf{u}_{,2} \\[6pt]
\mathbf{x}_{,1}^{\mathrm{T}}\cdot \Delta \mathbf{u}_{,3}+\mathbf{x}_{,3}^{\mathrm{T}}\cdot \Delta \mathbf{u}_{,1}
\end{bmatrix}.
\end{equation}
As the variational Green Lagrange strain tensor $\delta \hat{\mathbf{E}}$ is a nonlinear function of the displacement field $\mathbf{u}$, its linearization reads in the Voigt notation as:
\begin{equation}
\label{eqn:Kinematics.5}
\Delta \delta \hat{\mathbf{E}}=\begin{bmatrix}
\Delta\delta E_{\xi\xi} \\[6pt] 
\Delta\delta E_{\eta \eta} \\[6pt] 
\Delta\delta E_{\zeta \zeta} \\[6pt] 
2\Delta\delta E_{\xi \eta} \\[6pt] 
2\Delta\delta E_{\eta \zeta} \\[6pt] 
2\Delta\delta E_{\xi \zeta}
\end{bmatrix}=\begin{bmatrix}
\delta \mathbf{u}_{,1}^{\mathrm{T}}\cdot \Delta \mathbf{u}_{,1} \\[6pt] 
\delta \mathbf{u}_{,2}^{\mathrm{T}}\cdot \Delta \mathbf{u}_{,2} \\[6pt] 
\delta \mathbf{u}_{,3}^{\mathrm{T}}\cdot \Delta \mathbf{u}_{,3} \\[6pt] 
\delta \mathbf{u}_{,1}^{\mathrm{T}}\cdot \Delta \mathbf{u}_{,2} + \delta \mathbf{u}_{,2}^{\mathrm{T}}\cdot \Delta \mathbf{u}_{,1} \\[6pt]
\delta \mathbf{u}_{,2}^{\mathrm{T}}\cdot \Delta \mathbf{u}_{,3} + \delta \mathbf{u}_{,3}^{\mathrm{T}}\cdot \Delta \mathbf{u}_{,2} \\[6pt]
\delta \mathbf{u}_{,1}^{\mathrm{T}}\cdot \Delta \mathbf{u}_{,3}+\delta \mathbf{u}_{,3}^{\mathrm{T}}\cdot \Delta \mathbf{u}_{,1}
\end{bmatrix}.
\end{equation}

\subsection{Hyperelastic constitutive models}

In hyperelastic material modeling, the relation between stresses and strains is expressed through a strain energy function $W$:
\begin{equation} \label{eqn:Constitutive models.W}
    \mathbf{P} = \frac{dW}{d\mathbf{F}} 
    \qquad\Leftrightarrow\qquad
    \mathbf{S} =  \frac{dW}{d\mathbf{E}} = 2 \frac{dW}{d\mathbf{C}}. 
\end{equation}
Here, $\mathbf{P}$ denotes the first and $\mathbf{S}$ the second Piola-Kirchhoff stress tensor, respectively.

In analogy to the linear elastic Hooke's law, the simplest hyperelastic material model is the St. Venant-Kirchhoff (SVK) model. The strain energy function and constitutive relation for the SVK material are given by:
\begin{equation}
\label{eqn:Constitutive models.1}
\begin{aligned}
W(\mathbf{E})&=\frac{1}{2}\lambda [\operatorname{tr}(\mathbf{E})]^{2}+\mu \mathbf{E}:\mathbf{E}, \\ 
\mathbf{S}(\mathbf{E})&=\frac{\partial W}{\partial \mathbf{E}}= \lambda \operatorname{tr}(\mathbf{E})\mathbf{I}+2\mu \mathbf{E}.
\end{aligned}
\end{equation}
The Lamé constants $\lambda=\frac{E \nu}{(1+\nu)(1-2 \nu)}$ and $\mu=\frac{E}{2(1+\nu)}$ can be expressed in terms of the Young's modulus $E>0$ and Poisson's ratio $-1<\nu<0.5$.

Furthermore, in this work, we also employ the compressible Neo-Hookean (NH) material model, which can be expressed as:
\begin{equation}
\label{eqn:Constitutive models.2}
\begin{aligned}
W(\mathbf{C}) & =\frac{\lambda}{2}(\ln J)^2+\frac{\mu}{2}(\operatorname{tr} \mathbf{C}-3)+\mu \ln J, \\
\mathbf{S}(\mathbf{C}) & =2 \frac{d W}{d \mathbf{C}}=\lambda \ln J \mathbf{C}^{-1}+\mu\left(\mathbf{I}-\mathbf{C}^{-1}\right),
\end{aligned}
\end{equation}
where $J=\det\mathbf{F}$ is the determinant of the deformation gradient.

\subsection{Governing equations}
\label{sec:Governing equations}

The static equilibrium equation for the hyperelastic solid element is given by the balance of linear momentum, which can be expressed in the reference configuration as:
\begin{equation}
\label{eqn:Governing equations.1}
    \text{Div}_{\mathbf{X}}\mathbf{P}+\rho_{0}\mathbf{b} = \mathbf{0} \quad \text{on} \quad   \Omega_0,  
\end{equation}
where $\rho_{0}$ is the mass density and $\rho_{0}\mathbf{b}$ the body force. To complete the boundary value problem, Dirichlet and Neumann boundary conditions are provided:
\begin{equation}
\label{eqn:Governing equations.2}
\mathbf{u}=\mathbf{u}^{*} \quad \text{on} \quad  \Gamma_{u} \qquad \text{and} \qquad \mathbf{P}\cdot\mathbf{n}=\mathbf{t}^{*} \quad \text{on} \quad  \Gamma_{t},
\end{equation}
respectively. Here, $\mathbf{u}^{*}$ is the given displacement on the boundary $\Gamma_{u}\subset\partial\Omega_0$, $\mathbf{n}^{*}$ is the prescribed traction on the boundary $\Gamma_{t}\subset\partial\Omega_0$, and $\mathbf{n}$ is the outward normal vector on $\Gamma_t$. 

For the finite element implementation, the strong form of the governing equations as given in \cref{eqn:Governing equations.1,eqn:Governing equations.2} is typically transformed into its weak form.
Considering $\delta \mathbf{u}$ as a virtual displacement, variation, or test function with $\delta \mathbf{u} \in \{H^{1}(\Omega_0)\,|\,\delta \mathbf{u}=\mathbf{0}~ \text{on}~ \Gamma_{u}\}$, multiplying the static equilibrium equation with this test function, doing integration by parts, and applying divergence theorem along with Cauchy's theorem leads to the weak or variational form of the equilibrium equation as:
\begin{equation}
\label{eqn:Governing equations.3}
g_u(\mathbf{u},\delta \mathbf{u}) = g_{u}^{int}(\mathbf{u},\delta \mathbf{u}) - g_{u}^{ext}(\mathbf{u},\delta \mathbf{u}) = 0 \quad\forall\delta\mathbf{u}.
\end{equation}
The internal energy variation contribution to the weak form $g_{u}^{int}$ can likewise be expressed in terms of the first or second Piola-Kirchhoff tensors, since $\delta W=\delta\mathbf{F} : \mathbf{P} = \delta\mathbf{E} : \mathbf{S}$, as: 
\begin{equation}
\label{eqn:Governing equations.4}
g_{u}^{int}(\mathbf{u},\delta \mathbf{u}) =
\int_{\Omega _{0}}\delta \mathbf{F} : \mathbf{P}(\mathbf{F})~dV =\int_{\Omega _{0}}\delta \mathbf{E} : \mathbf{S}(\mathbf{E})~dV 
\end{equation}
and the external energy variation contribution as:
\begin{equation}
\label{eqn:Governing equations.5}
g_{u}^{ext}(\mathbf{u},\delta \mathbf{u}) = \int_{\Gamma_{t}}\delta \mathbf{u}\cdot\mathbf{t}~dA + \int_{ \Omega _{0}}\delta \mathbf{u}\cdot \rho_{0}\mathbf{b}~dV.
\end{equation}

As the weak form \eqref{eqn:Governing equations.3} is a nonlinear function in terms of the sought displacement field $\mathbf{u}\in\{ H^1(\Omega)\,|\,\mathbf{u}=\mathbf{u}^*~\text{on}~\Gamma_u\}$, its discretized counterpart, c.f.~\cref{sec: Isogeometric finite element discretizaton}, must typically be solved iteratively, e.g., using a Newton-Raphson scheme, which requires its linearization. Thus, the directional derivative of the variational form is given as:
\begin{equation}
\label{eqn:Governing equations.6}
D[g_u]\cdot \Delta \mathbf{u}=\Delta g_u = \underset{Geometric~part}{\underbrace{\int_{\Omega _{0}}\Delta \delta \mathbf{E} : \mathbf{S}(\mathbf{E})~dV }} + \underset{Material~part}{\underbrace{\int_{\Omega _{0}}\delta \mathbf{E} : \mathbb{C}(\mathbf{E}) : \Delta \mathbf{E}~dV}} - \Delta{g}_{u}^{ext},
\end{equation}
where $\mathbb{C}=\frac{d^2W}{d\mathbf{E}^2}$.
In the case of conservative loading, 
where the external loads are independent of the deformations, i.e., ${g}_{u}^{ext}\equiv{g}_{u}^{ext}(\delta\mathbf{u})$, it is $\Delta{g}_{u}^{ext} = 0$.
However, when non-conservative loading is applied, i.e., ${g}_{u}^{ext}\equiv{g}_{u}^{ext}(\mathbf{u},\delta\mathbf{u})$, it is  $\Delta g_{u}^{ext} \neq  0$ and the tangent matrix obtained from that term needs to be considered. A detailed description can be found in~\ref{subsec: Deformation dependent}. 


\section{Isogeometric finite element discretizaton}
\label{sec: Isogeometric finite element discretizaton}

The core idea of isogeometric analysis is to model the geometry and discretize the kinematic fields both by using spline basis functions such as B-splines or NURBS, which are commonly used in computer-aided design \cite{Hughes.2005}.
Here, IGA is applied for the discretization and approximation of the solution of the weak form  \cref{eqn:Governing equations.3}.

\subsection{B-splines and NURBS}

Detailed information on splines and IGA can be found in the books of Piegl and Tiller \cite{PiegTill96} and Cottrell \textit{et al.} \cite{Cottrell.2009}, respectively.
Here, we briefly summarize the main definitions relevant to solid-beam modeling.

B-spline basis functions can be defined and computed using the Cox-de Boor recursion formulas.  
For a polynomial degree of zero ($p=0$), the B-spline basis functions are given as:
\begin{equation}
\label{eqn:NURBS.1}
B_{i,0}(\xi )=\left\{\begin{matrix}
1 & \text{if}~ \xi _{i}\leq \xi < \xi _{i+1} \\ 
0 & \text{otherwise}
\end{matrix}\right.,
\end{equation}
and for higher degrees ($p\geq 1$) as:
\begin{equation}
\label{eqn:NURBS.2}
B_{i,p}(\xi )=\frac{\xi - \xi_{i}}{\xi_{i+p} - \xi_{i}}B_{i,p-1}(\xi )+\frac{\xi_{i+p+1} - \xi}{\xi_{i+p+1} - \xi_{i+1}}B_{i+1,p-1}(\xi ),
\end{equation}
where the knot vector $\mathbb{I}_{\xi }=[\xi_{1},\xi_{2},...,\xi_{n+p+1}]$ is a non-decreasing sequence of real numbers. 
3D solid geometries can be modelled as B-spline volumes:
\begin{equation}
\label{eqn:NURBS.3}
\textbf{V}(\xi,\eta,\zeta ) = \sum_{i=1}^{n}\sum_{j=1}^{m}\sum_{k=1}^{l} B_{i,p}(\xi)\, C_{j,q}(\eta)\, D_{k,r}(\zeta)\, \textbf{b}_{ijk},
\end{equation}
where $\textbf{b}_{ijk}\in\mathbb{R}^3$ are the control points and $B_{i,p}(\xi)$, $C_{j,q}(\eta)$ and $D_{k,r}(\zeta)$ are the B-spline functions of order $p$, $q$ and $r$, corresponding to knot vectors $\mathbb{I}_{\xi }=[\xi_{1},\xi_{2},...,\xi_{n+p+1}]$, $\mathbb{I}_{\eta }=[\eta_{1},\eta_{2},...,\eta_{m+q+1}]$ and $\mathbb{I}_{\zeta }=[\zeta_{1},\zeta_{2},...,\zeta_{l+r+1}]$, respectively.

NURBS are rational B-Spline functions that allow to describe conic geometries like circles and ellipsoids, which B-spline functions are not able to model. NURBS of order $p$ are defined as:
\begin{equation}
\label{eqn:NURBS one dim}
N_{i,p}(\xi )=\frac{B_{i,p}(\xi )\cdot w_{i}}{\sum\limits_{i=1}^{n}B_{i,p}(\xi )\cdot w_{i}},
\end{equation}
here $w_{i}\in\mathbb{R}$ are the weights of the control points, which allow the proper representation of curves.
Introducing weights $w_{ijk}\in\mathbb{R}$, the trivariate NURBS basis functions are defined as:
\begin{equation}
\label{eqn:NURBS.4}
N_{ijk,pqr}(\xi,\eta ,\zeta ) = \frac{B_{i,p}(\xi )\, C_{j,q}(\eta )\, D_{k,r}(\zeta )\, w_{ijk}}{\sum\limits_{i=1}^{n} \sum\limits_{j=1}^{m} \sum\limits_{k=1}^{l} B_{i,p}(\xi)\,  C_{j,q}(\eta)\, D_{k,r}(\zeta )\, w_{ijk}},
\end{equation}
and then NURBS solids analogous to \cref{eqn:NURBS.3}.

The implementation of isogeometric finite elements is quite similar to the standard isoparametric FEM, except that the computation of the basis functions and their derivatives. 

\subsection{Discretaization of field variables}
\label{sec: Discretaization of field variables}
In IGA, the approximated position vector in the reference configuration and the displacement field are discretized using trivariate NURBS functions, and consequently also the position vectors in the current configuration and the variational displacement field:
\begin{equation}
\label{eqn:Isogeometric finite element discretization:1}
\begin{aligned}
\mathbf{X}^{h} &= \sum_{I = 1}^{nel}N_{I}(\xi,\eta,\zeta )\,\mathbf{X}_{I}, \qquad & 
\mathbf{u}^{h} &= \sum_{I = 1}^{nel} N_{I}(\xi,\eta,\zeta )\, \mathbf{u}_{I},\\
\mathbf{x}^{h} &= \sum_{I = 1}^{nel} N_{I}(\xi,\eta,\zeta )\, (\mathbf{X}_{I} + \mathbf{u}_{I}), \qquad &
\delta \mathbf{u}^{h} &= \sum_{I = 1}^{nel} N_{I}(\xi,\eta,\zeta )\, \delta \mathbf{u}_{I}
.
\end{aligned}
\end{equation}
where $nel=n\cdot m\cdot l$ is the total number of trivariate basis functions and control points of the discretized solid element. In comparison to the definition of NURBS basis functions in \cref{eqn:NURBS.4}, for the simplification of the notation, $I\sim ijk$ is to be understood as a multi-index and the dependency on the degrees is dropped in the notation.

Then, the approximated form of the variational Green-Lagrange strains \eqref{eqn:Kinematics.3} is given as:
\begin{equation}
	\label{eqn:Isogeometric finite element discretization:2}
\delta 	\hat{\mathbf{E}}^{h}=\sum_{I = 1}^{nel}\begin{bmatrix}
F_{11}N_{I,1} & F_{21}N_{I,1} & F_{31}N_{I,1} \\[6pt]
F_{12}N_{I,2} & F_{22}N_{I,2} & F_{32}N_{I,2} \\[6pt]
F_{13}N_{I,3} & F_{23}N_{I,3} & F_{33}N_{I,3}\\[6pt]
F_{11}N_{I,2}+F_{12}N_{I,1} & F_{21}N_{I,2}+F_{22}N_{I,1} & F_{31}N_{I,2}+F_{32}N_{I,1}\\[6pt]
F_{12}N_{I,3}+F_{13}N_{I,2} & F_{22}N_{I,3}+F_{23}N_{I,2} & F_{32}N_{I,3}+F_{33}N_{I,2}\\[6pt]
F_{13}N_{I,1}+F_{11}N_{I,3} & F_{23}N_{I,1}+F_{21}N_{I,3} & F_{33}N_{I,1}+F_{31}N_{I,3} 
\end{bmatrix}\! \delta \mathbf{u}_{I} =\sum_{I = 1}^{nel} \mathbf{B}_{I}\, \delta \mathbf{u}_{I},
\end{equation}
where $\mathbf{B}_{I}$ is the corresponding strain-displacement operator  of control point $I$. The B-matrix of the discretized element can be assembled as $\mathbf{B} = [\mathbf{B}_{1},...\mathbf{B}_{I},...\mathbf{B}_{nel}]$. In addition, the B-operator can be used to get the incremental Green-Lagrange strains as:
\begin{equation}
	\label{eqn:Isogeometric finite element discretization:3}
 \Delta \hat{\mathbf{E}}^{h}=\sum_{I = 1}^{nel} \mathbf{B}_{I}\, \Delta \mathbf{u}_{I}.
\end{equation}
Similarly, the interpolated form of the linearized variational Green-Lagrange strains can be obtained as
\begin{equation}
	\label{eqn:Isogeometric finite element discretization:4}
\Delta \delta \hat{\mathbf{E}}^{h}=\sum_{I=1}^{nel}\sum_{J=1}^{nel} \delta\mathbf{u}_{I}^{\mathrm{T}}\, \mathfrak{B}_{IJ}\, \delta\mathbf{u}_{J},
\end{equation}
where $\mathfrak{B}_{IJ}$ is the third-order operator of dimension $nel\times 6\times  nel$ with:
\begin{equation}
\label{eqn:Isogeometric finite element discretization:5}
\mathfrak{B}_{IJ}=\begin{bmatrix}
N_{I,1}N_{J,1}\\[6pt] 
N_{I,2}N_{J,2}\\[6pt] 
N_{I,3}N_{J,3}\\[6pt] 
N_{I,1}N_{J,2}+N_{I,2}N_{J,1}\\[6pt]
N_{I,2}N_{J,3}+N_{I,3}N_{J,2}\\[6pt]
N_{I,1}N_{J,3}+N_{I,3}N_{J,1} 
\end{bmatrix}.
\end{equation}

In the standard isogeometric analysis for nonlinear elasticity, the discretizations and discretized expressions of \crefrange{eqn:Isogeometric finite element discretization:2}{eqn:Isogeometric finite element discretization:5} are substituted into the weak form \cref{eqn:Governing equations.3} and the resulting nonlinear system of algebraic equations is solved for the displacement control points $\mathbf{u}_I$ \cite{Cottrell.2009}.


\section{Isogeometric solid-beam formulation}
\label{Sec:Isogeometric solid-beam formulation}
Using solid elements in the analysis of slender beam structures leads to  locking effects such as shear, membrane, and Poisson thickness locking, which deteriorate the accuracy and convergence behavior \cite{Echter.2010}. To alleviate these locking effects and ensure the robustness of the numerical solution procedure, an isogeometric solid-beam element is developed using ANS, EAS, and MIP methods.  

\subsection{Assumed natural strain method}
\label{Sec:ANS method}
The assumed natural strain (ANS) method alleviates transversal shear and membrane locking, which typically appear in thin structures such as shells and beams and are caused by using low-order shape functions \cite{Frischkorn.2013,caseiro2014assumed}.

\begin{itemize}
  \item[\textbf{--}]Transversal shear locking occurs due to the development of shear stresses when the element is subjected to out-of-plane bending, see \cref{fig:ANS_1}. In \cref{fig:ANS_1a}, the exact bending behaviour of the element is represented and the approximated behaviour with locking due to low order shape function is shown in  \cref{fig:ANS_1b}. In the case of a pure out-of-plane bending moment,  shear stresses are developed, because the edges of the element are not perpendicular to the mid-plane as shown in \cref{fig:ANS_1b}. These undesirable shear stresses can lead the element to a nonphysical equilibrium with smaller displacements.

  \begin{figure}[t]
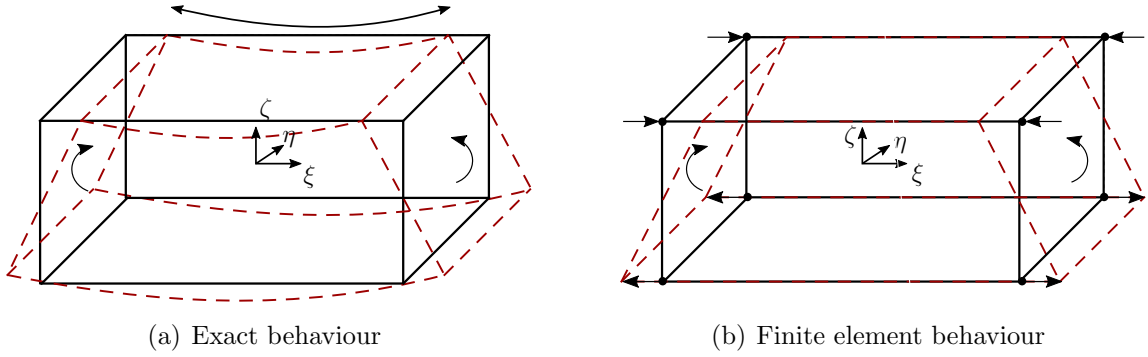

	\centering
    \subfigure[ Exact behaviour \label{fig:ANS_1a}]{\resizebox{0.4\textwidth}{!}{
        \input{exact_shear_locking}}}
    ~~~~~~\subfigure[Finite element behaviour \label{fig:ANS_1b}]{\resizebox{0.4\textwidth}{!}{
        \input{fem_shear_locking}}}	
    \caption{ Transversal shear locking due to low-order shape functions (solid lines - undeformed, dashed lines - deformed) }
    \label{fig:ANS_1}
\end{figure}

\item[\textbf{--}] Membrane locking occurs in slender curved beams due to the inability of curved solid finite elements to represent pure bending (inextensional) deformation. The curved finite elements overestimate the bending stiffness and exhibit parasitic membrane stresses. Membrane locking is also named inextensional locking, which is not a problem in straight beams, but already very slight curvature in thinner beams is sufficient to produce this artefact.
\end{itemize}

To cure these locking effects, the main idea of the ANS method is to calculate the strains (which are causing locking) at the so-called tying points and then a projection matrix is used to associate these strains at tying points back to the standard integration points. To alleviate the locking and avoid numerical instabilities, an appropriate selection of tying point locations is very important. Here, the tying points are defined with the help of reduced integration, i.e., they are located at the Gauss points of a reduced quadrature rule with $p$ instead of $p+1$ integration points. 

\begin{figure}[t]
	\centering
	\resizebox{0.7\textwidth}{!}{
        \input{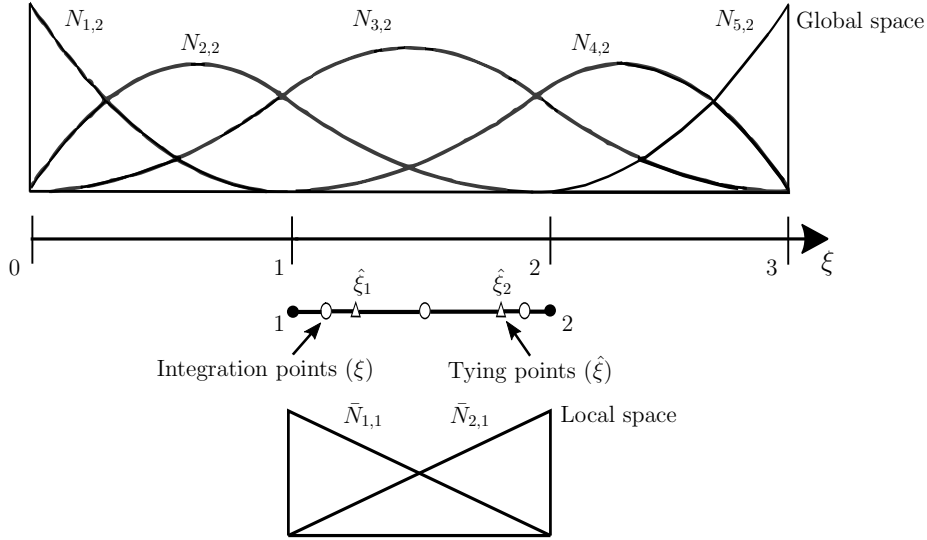}
    }
    \caption{Global and local spaces for the quadratic NURBS element}
    \label{fig:ANS_2}
\end{figure}

Following the work of Caseiro \textit{et al.} \cite{caseiro2014assumed,caseiro2015assumed} for the isogeometric analysis of solid-shell finite elements, the ANS method for isogeometric analysis of solid-beam elements is defined accordingly. For the second-order NURBS functions  along the beam  direction $\xi$, the tying points for the assumed natural strain fields (membrane strain $E_{\xi \xi }$ and transversal shear strains $E_{\xi \eta }$, $E_{\xi \zeta }$)  are defined by a two-point Gaussian quadrature rule (i.e., a reduced integration rule) in the $\xi$-direction of the element, as can be seen in \cref{fig:ANS_2}.
The univariate NURBS basis functions defined along the beam direction $\xi$ by the knot vector $\mathbb{I}_{\xi }$ are depicted as the ``global space''. On each element, i.e., each knot span, a local knot vector $\mathbb{\bar{I}}_{\xi }$ defines the ``local space'' for the NURBS basis functions of lower order $p-1$. 

Consider as an example the second-order ($p=2$) NURBS basis functions with knot vector $\mathbb{I}_{\xi }=[0,0,0,1,2,3,3,3]$ with three elements (knot spans) representing the global space as shown in \cref{fig:ANS_2}. The local space of order $p-1=1$ of each element can be determined accordingly, for example, the local knot vector of the middle element is $\mathbb{\bar{I}}_{\xi }=[1,1,2,2]$. Accordingly, the integration points (circles) and the tying points (triangles) of an element are obtained using Gaussian quadrature rules with $p+1=3$ and $p=2$ points, respectively.

Within each element, the assumed natural strain fields for $E^{ANS}_{\xi \xi }$, $E^{ANS}_{\xi \eta }$ and $E^{ANS}_{\xi \zeta }$ are evaluated at the current Gauss points as:
\begin{equation}
\label{eqn:IGA.1}
E_{\square }^{ANS}(\xi,\eta,\zeta)
=\sum_{j=1}^{p} \bar{N}_{j,p-1}(\xi)\, \bar{E}_{\square,j }(\eta,\zeta) 
=\sum_{j=1}^{p} P_{j}(\xi) \, E_{\square }(\hat{\xi}_{j},\eta,\zeta)
\quad  \text{for}\quad  \square \in\{\xi\xi,\,\xi\eta,\,\xi\zeta\}.
\end{equation}
Here, $\bar{N}_{j,p-1}(\xi)$ are the shape functions of the local space, evaluated at the quadrature point (or arbitrary coordinate) $\xi$, and   $\bar{E}_{\square,j }(\eta,\zeta)$ are the ANS  strain coefficients.
These are calculated such that the ANS strains interpolate the global strains at the tying points with coordinates $(\hat{\xi}_{j},\eta,\zeta)$:
\begin{equation}
\begin{aligned}
\label{eqn:IGA.2}
E_{\square }(\hat\xi_j,\eta,\zeta) &= \sum_{k=1}^{p} \bar{N}_{k,p-1}(\hat\xi_j)\, \bar{E}_{\square,k }(\eta,\zeta) \quad  \text{for}\quad j=1,\ldots,p~,\\
\Leftrightarrow \qquad
\mathbf{E}_{\square }(\eta,\zeta) &= \mathbf{M} \cdot \bar{\mathbf{E}}_{\square }(\eta,\zeta) , 
\end{aligned}
\end{equation}
where $\mathbf{E}_{\square },\, \bar{\mathbf{E}}_{\square }\in\mathbb{R}^{p},\, \mathbf{M}\in\mathbb{R}^{p\times p}$ with $\mathbf{M}_{kj}=\bar{N}_{k,p-1}(\hat{\xi}_j)$.
Thus, the ANS strain coefficients can directly be expressed in terms of the tying point evaluations of the global strains, $\bar{\mathbf{E}}_{\square}(\eta,\zeta) = \mathbf{M}^{-1} \cdot \mathbf{E}_{\square}(\eta,\zeta)$, which allows to express  \cref{eqn:IGA.1} in using the projected shape functions $P_{j}(\xi)$, which are the entries of the vector:
\begin{equation}
\label{eqn:IGA.3}
\mathbf{P}(\xi) = \left( \bar{N}_{1,p-1}(\xi), \ldots, \bar{N}_{p,p-1}(\xi) \right) \cdot \mathbf{M}^{-1}
= \bar{\mathbf{N}}(\xi) \cdot \mathbf{M}^{-1}.
\end{equation}

The strain-displacement operator $\mathbf{B}$ and the third-order operator $\bar{\mathfrak{B}}$ for the assumed natural strains can then be defined as: 
\begin{equation}
\label{eqn:IGA.4}
\begin{aligned}
\textbf{B}_{\square}^{ANS}(\xi,\eta,\zeta)&=\sum_{j=1}^{p}P_{j}(\xi)\,\textbf{B}_{\square}(\hat{\xi}_{j},\eta,\zeta),\\
\mathfrak{B}_{\square}^{ANS}(\xi,\eta,\zeta)&=\sum_{j=1}^{p}P_{j}(\xi)\,\mathfrak{B}_{\square}(\hat{\xi}_{j},\eta,\zeta).
\end{aligned}
\end{equation}
In the end, the 1st, 4th, and 6th row of the $\textbf{B}$- and $\mathfrak{B}$-operators corresponding to $\square\in\{\xi\xi,\,\xi\eta,\,\xi\zeta\}$ are replaced by the modified fields according to \cref{eqn:IGA.4}. The variational formulation of a solid-beam element based only on the ANS method is outlined in \cref{sec:Governing equations}. 
In this NURBS-based displacement formulation, the interpolation based on the tying points is independent of the parametric coordinates in the thickness direction $\eta$ and $\zeta$, therefore leading to the so-called solid-beam formulation. 


\subsection{Enhanced assumed strain method }
\label{Sec:EAS method}

The enhanced assumed strain (EAS) method is employed to alleviate the Poisson thickness locking, which occurs for non-zero Poisson's ratio $\nu\neq 0$.
When linear approximation functions are used for the displacement field in the transversal direction of the beam, the transversal normal strains $E_{\eta \eta}$ and $E_{\zeta \zeta}$ become constant along the $\eta$- and $\zeta$-directions, respectively. This causes locking due to the transversal contraction. To incorporate realistic 3D constitutive models in solid-beam formulation, the transversal normal strains need to be at least linear nature over their respective beam thickness directions. 
Additionally, in bending situations for $\nu \neq 0$, a trapezoidal deformation mode of the cross-section of the beam occurs. To avoid the artificial shear strain $E_{\eta \zeta}$ and the associated in-plane shear locking, $E_{\eta \zeta}$ also needs to be enhanced.

To cure this material locking, in the EAS method, the strains are enhanced and the total strain field reads as $\mathbf{E}_{a} = \mathbf{E} + \mathbf{E}_{e}$. The enhancement of the in-plane strain components of the cross-section is defined as:
\begin{equation}
\label{Eq.EAS_1}
 \mathbf{E}_{e} ({\xi},\eta,\zeta)= \begin{bmatrix}
0 \\
E_{e,\eta \eta } \\
E_{e,\zeta \zeta } \\
0 \\
E_{e,\eta \zeta } \\
0\end{bmatrix}= \underset{\mathbf{B}_e(\eta,\zeta)}{\underbrace{ \begin{bmatrix}
0 & 0 & 0 & 0 & 0  \\
\eta & 0 & 0 & 0 & 0  \\
0 & \zeta & 0 & 0 & 0  \\
0 & 0 & 0 & 0 & 0  \\
0 & 0 & \eta & \zeta & \eta\zeta  \\
0 & 0 & 0 & 0 & 0 \\
\end{bmatrix}}} \boldsymbol{\alpha}^{h}(\tilde{\xi}),
\end{equation}
where $\tilde{\xi} \in [-1,~ 1]$ is the reparameterized coordinate of the beam direction from the global coordinate in an element (knot span) $\xi \in [\xi^{e}_{1},~ \xi^{e}_{2}]$ by the linear mapping:
\begin{equation}
\label{Eq.repar_eq}
\tilde{\xi} = 1 -2  \frac{\xi^{e}_{2}-\xi}{\xi^{e}_{2}-\xi^{e}_{1}}.
\end{equation}

The additional enhancement parameters $\boldsymbol{\alpha}^{h}$ are interpolated by Lagrange interpolation functions at the element level and thus they can be  statically condensed out due to inter-element discontinuity. The enhancement parameters and their variations are discretized as:
\begin{equation}
\label{Eq.EAS_2}
\boldsymbol{\alpha}^{h}(\tilde{\xi}) = \sum_{i=1}^{p_{\alpha}+1} L_{i,p_{\alpha}}(\tilde{\xi}) \, \boldsymbol{\alpha}_{i} = \mathbf{L}(\tilde{\xi}) \, \boldsymbol{\alpha}, \ 
\qquad \quad 
\delta\boldsymbol{\alpha}^{h}(\tilde{\xi}) = \sum_{i=1}^{p_{\alpha}+1} L_{i,p_{\alpha}}(\tilde{\xi}) \, \delta\boldsymbol{\alpha}_{i} = \mathbf{L}(\tilde{\xi}) \, \delta\boldsymbol{\alpha},
\end{equation}
where  
\begin{equation}
\label{Eq.lag_func}
\mathbf{L}(\tilde{\xi}) = \begin{cases}
\mathbf{1}_{5\times5}, &  p=1; \\
[L_{1,1}(\tilde{\xi})\mathbf{1}_{5\times5} , \ldots, L_{p,p_{\alpha}}(\tilde{\xi})\mathbf{1}_{5\times5}], &  p\geq 2, 
\end{cases}
\end{equation}
$\boldsymbol{\alpha}_{i},\ \delta \boldsymbol{\alpha}_{i} \in\mathbb{R}^5$ and $L_{i,p_{\alpha}}(\tilde{\xi})$ are the Lagrange polynomials of degree $p_{\alpha}= p -1$, which are described in \ref{subsec: Lagrange interpolation functions}.

Simo and Amero \cite{Simo1992} developed the EAS method based on the Hu-Washizu three-field mixed formulation, in which the independent stresses get eliminated by the $L^{2}$-orthogonality condition between the enhanced strain field and the independent stresses due to which the enhanced strain field does not contribute to the virtual work. Thus, the weak form for EAS elements is:
\begin{equation}
\label{Eq.EAS_3}
\begin{aligned}
g_{u}^{EAS}(\mathbf{u},\boldsymbol{\alpha},\delta\mathbf{u}) &= \int_{\Omega _{0}}\delta \mathbf{E} : \mathbf{S}(\mathbf{E}_{a})~{dV} - g_{u}^{ext} &=0 \quad & \forall\delta\mathbf{u},\\
g_{\alpha}^{EAS}(\mathbf{u},\boldsymbol{\alpha},\delta\boldsymbol{\alpha}) &= \int_{\Omega _{0}}\delta \mathbf{E}_{e} : \mathbf{S}(\mathbf{E}_{a})~{dV} &=0 \quad & \forall\delta\boldsymbol{\alpha}.
\end{aligned}
\end{equation}
Here, $\mathbf{S}(\mathbf{E}_{a})$ are the stresses computed from the constitutive model based on the total strains $\mathbf{E}_{a}$. The discretized variational formulation and the linearized weak form of the solid-beam based only on the EAS method is given in \cref{eqn:Discretizated Variational formulation.4}, where the geometric stiffness $\mathbf{K}_{uu}^{e}$ is based on $\mathbf{S}(\mathbf{E}_{a})$.

\subsection{Mixed integration method}
\label{Sec:Mixed integration method}

One possible drawback of the EAS method is the weak robustness during the Newton-Raphson scheme, which often requires small load increments and a large number of  iterations to converge \cite{Pfefferkorn.2021}. To make the formulation more robust and speed up convergence, the mixed integration point (MIP) method is used here.

The key idea of the MIP is to relax the constitutive model at each integration point by introducing the independent stresses at the Gauss quadrature points. Unlike the classical mixed method, these independent stresses are not continuous and are not discretized but are defined locally at each integration point. Due to the local definition, they can be analytically condensed out by the static condensation, such that only displacement and the additional variables need to be discretized, like in a standard EAS-based element. The MIP method finally leads to a modification of the geometric stiffness matrix, which is defined by the independent stresses, while the residual vector remains unchanged. Thus, the converged results remain the same, but the robustness of geometrically nonlinear problems is greatly improved \cite{Magisano.2017}.

The extension of the presented isogeometric solid-beam formulation to the MIP method leads to the final weak form  based on ANS, EAS, and MIP methods, c.f.~\cite{Pfefferkorn.2021}:
\begin{equation}
\label{eqn:MIP.1}
 \begin{aligned}
g_{u}^{ISB}(\mathbf{S}_\beta,\delta\mathbf{u}) &= \int_{\Omega _{0}}\delta \mathbf{E} : \mathbf{S}_{\beta }~{dV} - g_{u}^{ext} &=0 \quad & \forall\delta\mathbf{u},\\
g_{\alpha}^{ISB}(\mathbf{S}_\beta,\delta\boldsymbol{\alpha}) &= \int_{\Omega _{0}}\delta \mathbf{E}_{e} : \mathbf{S}_{\beta }~{dV} &=0 \quad &  \forall\delta\boldsymbol{\alpha}, \\
g_{S}^{ISB}(\mathbf{u},\boldsymbol{\alpha},\mathbf{S}_\beta,\delta\mathbf{S}_\beta) &= \int_{\Omega _{0}}\delta \mathbf{S}_{\beta} : (\mathbf{E}_{a}-\mathbf{E}_{\beta })~{dV} &=0 \quad & \forall\delta\mathbf{S}_\beta,
\end{aligned}
\end{equation}
where $\mathbf{S}_{\beta}$ are the independent stresses and $\mathbf{E}_{\beta}$ are the constitutive strains, which are calculated by the inverse of the material model such that $\mathbf{S}_{\beta}=\frac{dW}{d\mathbf{E}_\beta}(\mathbf{E}_{\beta})$ holds.
Note that, if the ANS method is used, the strains $\mathbf{E}$ (and thus also the total, enhanced strains $\mathbf{E}_a$ in the EAS method), as well as $\delta\mathbf{E}$, include the modification of the $\xi\xi$-, $\xi\eta$-, and $\xi\zeta$-components according to \cref{Sec:ANS method}.


\subsection{Discretizated variational form and static condensation}
\label{Sec:Discretizated Variational formulation}

For the approximation of the solution of the weak form of the solid-beam formulation with ANS, EAS, and MIP methods, see \cref{eqn:MIP.1}, we apply an isogeometric finite element discretization to the displacement field $\mathbf{u}^h$ and its variation $\delta\mathbf{u}^h$, see \cref{sec: Discretaization of field variables}, an element-wise discretization of the EAS parameters $\boldsymbol{\alpha}^h$ and their variations $\delta\boldsymbol{\alpha}^h$, see \cref{Eq.EAS_2}, and a quadrature-point wise approximation of the MIP stresses $\mathbf{S}_\beta$ and their variations $\delta\mathbf{S}_\beta$, see \cref{Sec:Mixed integration method}.

Using Gaussian quadrature with weights $w_{g}$, the approximation of the weak form \eqref{eqn:MIP.1} on element level read as:
\begin{equation}
\label{eqn:Discretizated Variational formulation.1}
\begin{aligned}
\tilde{g}_{u}&=\sum_{g=1}^{n_{gp}} w_{g}\, \delta \mathbf{u}^{e}\cdot \left(\mathbf{B}_{g}^{\mathrm{T}}\,\mathbf{S}_{\beta,g}\right) - g_{u}^{ext} &= 0,\\
\tilde{g}_{\alpha}&=\sum_{g=1}^{n_{gp}} w_{g}\, \delta \boldsymbol{\alpha}^{e}\cdot \left((\mathbf{B}_{e,g}\mathbf{L}_{g})^{\mathrm{T}}\,\mathbf{S}_{\beta,g}\right)
&=0, \\
\tilde{g}_{S}&=\sum_{g=1}^{n_{gp}} w_{g}\, \delta \mathbf{S}_{\beta,g}\cdot (\mathbf{E}_{a,g}-\mathbf{E}_{\beta,g}) &=0,
\end{aligned}
\end{equation}
where $n_{gp}$ is the number of the Gauss quadrature points. By considering the variation of the independent variables $\delta \mathbf{u}^{e},~\delta \boldsymbol{\alpha}^{e},~ \delta \mathbf{S}_{\beta,g}$ as arbitrary, the linearization of the formulation \cref{eqn:Discretizated Variational formulation.1} at the element level can be written as:
\begin{equation}
\label{eqn:Discretizated Variational formulation.2}
\left[\begin{array}{ccccc}
\sum\limits_{g=1}^{n_{gp}}w_{g} \mathfrak{B}_{g}\mathbf{S}_{\beta,g} & \mathbf{0} & w_{1}\mathbf{B}_{1}^{\mathrm{T}} & \ldots & w_{n_{gp}}\mathbf{B}_{n_{gp}}^{\mathrm{T}} \\
\mathbf{0} & \mathbf{0} & w_{1}\mathbf{B}_{e,1}^{\mathrm{T}} & \ldots & w_{n_{gp}}\mathbf{B}_{e,n_{gp}}^{\mathrm{T}} \\
\mathbf{B}_{1} & \mathbf{B}_{e,1} & -\hat{\mathbb{D}}_1 & & \mathbf{0} \\
\vdots & \vdots & & \ddots & \\
\mathbf{B}_{n_{gp}} & \mathbf{B}_{e,n_{gp}} & \mathbf{0} & & -\hat{\mathbb{D}}_{n_{gp}}
\end{array}\right]\left[\begin{array}{c}
\Delta \mathbf{u}^{e} \\[6pt]
\Delta \boldsymbol{\alpha}^e \\[6pt]
\Delta \mathbf{S}_{\beta,1} \\
\vdots \\[6pt]
\Delta \mathbf{S}_{\beta,n_{gp}}
\end{array}\right]=-\left[\begin{array}{c}
\sum\limits_{g=1}^{n_{gp}} w_g\mathbf{B}_{g}^{\mathrm{T}} \mathbf{S}_{\beta,g} \\
\sum\limits_{g=1}^{n_{gp}} w_g\left(\mathbf{B}_{e,g}\mathbf{L}_{g}\right)^{\mathrm{T}} \mathbf{S}_{\beta,g} \\
\mathbf{E}_{a, 1}-\mathbf{E}_{\beta,1 } \\
\vdots \\
\mathbf{E}_{a, n_{gp}}-\mathbf{E}_{\beta,n_{gp}}
\end{array}\right].
\end{equation}

After the static condensation, the independent stresses at the Gauss points $g$ in each step $j$ of a Newton-Raphson scheme are obtained as:
\begin{equation}
\label{eqn:Discretizated Variational formulation.3}
\mathbf{S}_{\beta,g}^{j+1}= \left( \mathbf{S}_{\beta,g}+\hat{\mathbf{C}}_{g}\mathbf{B}_{g}\Delta \mathbf{u}^{e}+\hat{\mathbf{C}}_{g}\mathbf{B}_{e,g} \mathbf{L}_{g} \Delta \boldsymbol{\alpha}^{e} \right)^{j}.
\end{equation}
After substituting back the independent stresses into \cref{eqn:Discretizated Variational formulation.2}, the linearized weak form of the MIP method for the  solid-beam element reads as:
\begin{equation}
\label{eqn:Discretizated Variational formulation.4}
\begin{bmatrix}
\textbf{K}_{uu}^{e} & \textbf{K}_{u \alpha}^{e} \\
 \textbf{K}_{\alpha u}^{e}& \textbf{K}_{\alpha \alpha}^{e} \\
\end{bmatrix}\begin{bmatrix}
\Delta \textbf{u}^e \\
\Delta \boldsymbol{\alpha}^e
\end{bmatrix}=-\begin{bmatrix}
 \textbf{R}_{u}^{e} \\
 \textbf{R}_{\alpha}^{e}
\end{bmatrix},
\end{equation}
where 
\begin{equation}
\label{eqn:Discretizated Variational formulation.5}
\mathbf{R}_{u}^{e}= \sum_{g=1}^{n_{gp}}w_{g}\,\mathbf{B}_{g}^{\mathrm{T}}\,\mathbf{S}_{g}, 
\qquad \mathbf{R}_{\alpha}^{e}= \sum_{g=1}^{n_{gp}}w_{g} \, (\mathbf{B}_{e,g}\mathbf{L}_{g})^{\mathrm{T}}\, \mathbf{S}_{g},
\end{equation}
and
\begin{equation}
\label{eqn:Discretizated Variational formulation.6}
\begin{aligned}
\mathbf{K}_{uu}^{e} &= \sum_{g=1}^{n_{gp}}w_{g}\mathbf{B}_{g}^{\mathrm{T}}\hat{\mathbf{C}_{g}}\mathbf{B}_{g}+\sum\limits_{g=1}^{n_{gp}}w_{g} \mathbf{S}_{\beta,g}^{\mathrm{T}}\mathfrak{B}_{g}, \qquad &
\mathbf{K}_{u \alpha}^{e} &= \sum_{g=1}^{n_{gp}}w_{g}\mathbf{B}_{g}^{\mathrm{T}}\hat{\mathbf{C}_{g}}\mathbf{B}_{e,g}\mathbf{L}_{g}, \\
\mathbf{K}_{\alpha u}^{e} &= \sum_{g=1}^{n_{gp}}w_{g}(\mathbf{B}_{e,g}\mathbf{L}_{g})^{\mathrm{T}}\hat{\mathbf{C}_{g}}\mathbf{B}_{g},\qquad &
\mathbf{K}_{\alpha \alpha}^{e} &= \sum_{g=1}^{n_{gp}}w_{g}(\mathbf{B}_{e,g}\mathbf{L}_{g})^{\mathrm{T}}\hat{\mathbf{C}_{g}}\mathbf{B}_{e,g}\mathbf{L}_{g}.
\end{aligned}
\end{equation}
Here, the independent stresses $\mathbf{S}_{\beta,g}$ at each integration point contribute only to the geometric stiffness tensor, while the residuals are calculated by the constitutive material model $\mathbf{S}(\mathbf{E}_{a})$. In each Newton step $j$, the enhancement parameters $\boldsymbol{\alpha}^{e}$ are evaluated after the static condensation of \cref{eqn:Discretizated Variational formulation.4} at the element level as:
\begin{equation}
\label{eqn:Discretizated Variational formulation.7}
\boldsymbol{\alpha}^{e,j+1} = \boldsymbol{\alpha}^{e,j} + \Delta \boldsymbol{\alpha}^{e}, \qquad \text{where}
\qquad \Delta \boldsymbol{\alpha}^{e} = -\mathbf{K}_{\alpha \alpha}^{e^{-1}} \left(\mathbf{R}_{\alpha}^{e} + \mathbf{K}_{\alpha u}^{e}\Delta \mathbf{u}^{e} \right).
\end{equation}
After plugging back $\Delta \boldsymbol{\alpha}^{e}$ into \cref{eqn:Discretizated Variational formulation.4}$_{1}$, the element residual vector $\mathbf{R}^{e}$ and the element stiffness matrix $\mathbf{K}^{e}$ become
\begin{equation}
\label{eqn:Discretizated Variational formulation.8}
\mathbf{R}^{e} = \mathbf{R}_{u}^{e} - \mathbf{K}_{u \alpha}^{e}\mathbf{K}_{\alpha \alpha}^{e^{-1}} \mathbf{R}_{\alpha}^{e} , \qquad  
\mathbf{K}^{e} = \mathbf{K}_{uu}^{e} - \mathbf{K}_{u\alpha}^{e}\mathbf{K}_{\alpha \alpha}^{e^{-1}}\mathbf{K}_{\alpha u}^{e}.
\end{equation}
Algorithm \ref{alg:IGA_1} explains the detailed procedure to apply the robust IGA-based solid-beam element formulation in an element routine.

\begin{algorithm}[ht!]
\caption{Element-level evaluations of the robust finite strain isogeometric solid-beam element}
\label{alg:IGA_1}

\begin{algorithmic}[1]
\vspace{3mm}
\State \textbf{Begin} For loop over all elements $e$
\State \hspace*{5mm} Calculate $\Delta \boldsymbol{\alpha}^{e}$ from the stored $\mathbf{K}_{\alpha \alpha}^{e^{-1}}\mathbf{R}_{\alpha}^{e},~ \mathbf{K}_{\alpha \alpha}^{e^{-1}}\mathbf{K}_{\alpha u}^{e}$ matrices and $\Delta \mathbf{u}^{e}$ vector, and update \hspace*{5mm}~~ of $\boldsymbol{\alpha}^{e}$ for the current step, see \cref{eqn:Discretizated Variational formulation.7}. For first Newton step $\boldsymbol{\alpha}^{e}=\mathbf{0}$ and $\Delta \boldsymbol{\alpha}^{e}=\mathbf{0}$
\State \hspace*{5mm} \textbf{Begin} For loop over all Gauss integration points $g$
\State \hspace*{10mm} Calculate strain-displacement operator $\textbf{B}_{g}(\xi,\eta,\zeta)$ and third-order operator $\mathfrak{B}_{g}(\xi,\eta,\zeta)$ at \hspace*{10mm}~~ current integration point
\State \hspace*{10mm} Calculate $\bar{\mathbf{N}}(\xi)$ vector at
current integration point, see \cref{eqn:IGA.3} 
\State \hspace*{10mm} \textbf{Begin} For loop over tying points $\hat\xi$
\State \hspace*{15mm} Calculate $\mathbf{M}$-matrix of local basis functions at the tying points, see \cref{eqn:IGA.2}
\State \hspace*{15mm} Calculate strain-displacement operator $\textbf{B}_{\square}(\hat{\xi},\eta,\zeta)$ and third-order operator $\mathfrak{B}_{\square}(\hat{\xi},\eta,\zeta)$ \hspace*{15mm}~~ by using the global basis functions at tying point coordinates
\State \hspace*{15mm} Calculate modified ANS strain-displacement operator $\textbf{B}_{\square}^{ANS}(\xi,\eta,\zeta)$ and third-order \hspace*{15mm}~~ operator $\mathfrak{B}_{\square}^{ANS}(\xi,\eta,\zeta)$ by \eqref{eqn:IGA.4}
\State \hspace*{10mm} \textbf{End} of tying points routine
\State \hspace*{10mm} Replace the relevant lines of operators $\textbf{B}_g(\xi,\eta,\zeta)$ and $\mathfrak{B}_g(\xi,\eta,\zeta)$ with the modified \hspace*{10mm} ~~ANS operators $\textbf{B}_{\square}^{ANS}(\xi,\eta,\zeta)$ and $\mathfrak{B}_{\square}^{ANS}(\xi,\eta,\zeta)$, respectively 
\State \hspace*{10mm} Calculate the $\textbf{B}_{e,g}$ operator and accordingly the enhanced strains $\textbf{E}_{e,g}$ by \eqref{Eq.EAS_1} 
\State \hspace*{10mm} Calculate the constitutive stresses at each integration point $\mathbf{S}_{g}$ from the total strains \hspace*{10mm} ($\mathbf{E}_{a,g} = \mathbf{E}_{g} + \mathbf{E}_{e,g}$)
\State \hspace*{10mm} \textbf{if} first Newton iteration of a load step, then
\State \hspace*{15mm} $\mathbf{S}_{\beta,g} = \mathbf{S}_{g}$
\State \hspace*{10mm} \textbf{else}
\State \hspace*{15mm} $\mathbf{S}_{\beta,g} = \mathbf{S}_{\beta,g},$ according to \eqref{eqn:Discretizated Variational formulation.3}  
\State \hspace*{10mm} \textbf{end}
\State \hspace*{10mm} Calculate the residuals and stiffness matrices as in \eqref{eqn:Discretizated Variational formulation.5} and \eqref{eqn:Discretizated Variational formulation.6} from the \hspace*{10mm}~~ displacements and additional degrees of freedom
\State \hspace*{5mm} \textbf{End} of Integration points routine
\State \hspace*{5mm} Calculate the element residual vector and stiffness matrix by \eqref{eqn:Discretizated Variational formulation.8}
\State \textbf{End} of element routine
\end{algorithmic}
\end{algorithm}


\section{Numerical examples}
\label{Sec:Numerical Examples}

The developed isogeometric nonlinear solid-beam formulation is now validated on several benchmark examples. The proposed formulation is used to solve a set of numerical problems and the results are compared with reference or analytical solutions. 

A Newton-Raphson iterative procedure is used for all the reported calculations. Residual and energy norm convergence criteria are used in the Newton-Raphson iterations to check the convergence of the solution. The presented formulation is implemented in MATLAB using an open-source nonlinear isogeometric analysis framework \cite{Du.2020}.
The numerical problems show the desired quadratic rate of convergence of the Newton scheme that is even obtained for the large load steps. 

For comparison with other element formulations, the following nomenclature is used:
\begin{description}
  \item[ISB$p$] -- Proposed IGA solid-beam element with degree $p$ NURBS functions along the beam direction and degree 1 in the cross-section directions.
  \item[IGA$p$] -- Standard IGA solid element with degree $p$ NURBS functions along the beam direction and degree 1 in the cross-section directions. IGA$pqq$ indicates that instead of degree 1, degree $q$ is used in the cross-section directions.
  \item[Q1STb] -- Isoparametric solid-beam element with linear shape functions, based on ANS, EAS and reduced integration method, as proposed by \cite{Frischkorn.2013}.
\end{description}
In the following benchmark examples, convergence studies are done by the help of the relative difference of the numerical solution compared to the reference or analytical solution. The relative difference of displacements is calculated by: 
\begin{equation}
\label{Eq.Numerical example_1}
e_{u}^{rel} =\left|\frac{u_{ref}-u}{u_{ref}} \right|.
\end{equation}

\subsection{Cantilever beam under shear loading}
\label{sec:num:1}

\begin{figure}[t]
    \centering
    \subfigure[Cantilever beam under shear loading (\cref{sec:num:1}) \label{fig:Numerical example_1.1a}]{
    ~~~~\resizebox{0.40\textwidth}{!}{
        \begin{tikzpicture}

\coordinate (A1) at (0,0,0);
\coordinate (A2) at (2,0,0);
\coordinate (A3) at (2,2,0);
\coordinate (A4) at (0,2,0);

\coordinate (B1) at (0,0,20);
\coordinate (B2) at (2,0,20);
\coordinate (B3) at (2,2,20);
\coordinate (B4) at (0,2,20);

\coordinate (C1) at (0,0,20/3);
\coordinate (C2) at (2,0,20/3);
\coordinate (C3) at (2,2,20/3);
\coordinate (C4) at (0,2,20/3);

\coordinate (D1) at (0,0,20/1.5);
\coordinate (D2) at (2,0,20/1.5);
\coordinate (D3) at (2,2,20/1.5);
\coordinate (D4) at (0,2,20/1.5);

\coordinate (E1) at (0,0,21);
\coordinate (E2) at (2,0,21);
\coordinate (E3) at (2,2,21);
\coordinate (E4) at (0,2,21);

\coordinate (F1) at (-0.6,0,20);
\coordinate (F2) at (2.6,0,20);
\coordinate (F3) at (2.6,2,20);
\coordinate (F4) at (-0.6,2,20);

\coordinate (G1) at (0,-0.6,20);
\coordinate (G2) at (2,-0.6,20);
\coordinate (G3) at (2,2.6,20);
\coordinate (G4) at (0,2.6,20);

\coordinate (H1) at (2.5,0,20/3);
\coordinate (H2) at (2.5,0,20/1.5);

\coordinate (I1) at (3.8,0,0);
\coordinate (I2) at (3.8,0,20);

\coordinate (J1) at (0,-1,20);
\coordinate (J2) at (2,-1,20);

\coordinate (K1) at (-1,0,20);
\coordinate (K2) at (-1,2,20);

\coordinate (L1) at (0,0.75,0);
\coordinate (L2) at (2,0.75,0);
\coordinate (L3) at (2,2.75,0);
\coordinate (L4) at (0,2.75,0);

\coordinate (M1) at (5.5,0,20);
\coordinate (M2) at (5.5,1,20);
\coordinate (M3) at (5.5,0,18);
\coordinate (M4) at (6.5,0,20);

\coordinate (N1) at (1,1,0);
\coordinate (N2) at (3,1,0);
\coordinate (N3) at (3,2,0);

\draw[black] (A1) -- (A2);
\draw[black] (A2) -- (A3);
\draw[black] (A3) -- (A4);
\draw[black] (A4) -- (A1);

\draw[black] (B1) -- (B2);
\draw[black] (B2) -- (B3);
\draw[black] (B3) -- (B4);
\draw[black] (B4) -- (B1);

\draw[black, dashed] (A1) -- (B1);
\draw[black, rotate =20] (A2) -- (B2);
\draw[black, rotate =20] (A3) -- (B3);
\draw[black, rotate =20] (A4) -- (B4);

\filldraw [black] (A1) circle (2pt);
\filldraw [black] (A2) circle (2pt);
\filldraw [black] (A3) circle (2pt);
\filldraw [black] (A4) circle (2pt);

\filldraw [black] (B1) circle (2pt);
\filldraw [black] (B2) circle (2pt);
\filldraw [black] (B3) circle (2pt);
\filldraw [black] (B4) circle (2pt);

\filldraw [black] (C1) circle (2pt);
\filldraw [black] (C2) circle (2pt);
\filldraw [black] (C3) circle (2pt);
\filldraw [black] (C4) circle (2pt);

\draw[black, dashed] (C1) -- (C2);
\draw[black, dashed] (C2) -- (C3);
\draw[black, dashed] (C3) -- (C4);
\draw[black, dashed] (C4) -- (C1);

\filldraw [black] (D1) circle (2pt);
\filldraw [black] (D2) circle (2pt);
\filldraw [black] (D3) circle (2pt);
\filldraw [black] (D4) circle (2pt);

\draw[black, dashed] (D1) -- (D2);
\draw[black, dashed] (D2) -- (D3);
\draw[black, dashed] (D3) -- (D4);
\draw[black, dashed] (D4) -- (D1);

\draw[red] (B1) -- (E1);
\draw[red] (B2) -- (E2);
\draw[red] (B3) -- (E3);
\draw[red] (B4) -- (E4);

\draw[red] (B1) -- (F1);
\draw[red] (B2) -- (F2);
\draw[red] (B3) -- (F3);
\draw[red] (B4) -- (F4);

\draw[red] (B1) -- (G1);
\draw[red] (B2) -- (G2);
\draw[red] (B3) -- (G3);
\draw[red] (B4) -- (G4);

\draw[black, dashed] (C2) -- (H1);
\draw[black, dashed] (D2) -- (H2);

\draw[black,stealth-stealth] (H1) -- (H2);

\draw[black, dashed] (A2) -- (I1);
\draw[black, dashed] (B2) -- (I2);

\draw[black,stealth-stealth] (I1) -- (I2);

\draw[black, dashed] (B1) -- (J1);
\draw[black, dashed] (B2) -- (J2);
\draw[black,stealth-stealth] (J1) -- (J2);

\draw[black, dashed] (B1) -- (K1);
\draw[black, dashed] (B4) -- (K2);
\draw[black,stealth-stealth] (K1) -- (K2);

\draw[red,-stealth,line width=1.5pt] (A1) -- (L1);
\draw[red,-stealth,line width=1.5pt] (A2) -- (L2);
\draw[red,-stealth,line width=1.5pt] (A3) -- (L3);
\draw[red,-stealth,line width=1.5pt] (A4) -- (L4);

\draw[black,-stealth,line width=1pt] (M1) -- (M2);
\draw[black,-stealth,line width=1pt] (M1) -- (M3);
\draw[black,-stealth,line width=1pt] (M1) -- (M4);

\draw[blue,dashed, line width=1pt] (N1) -- (N2);
\draw[blue,-stealth,line width=1pt] (N2) -- (N3);

\node at (0.5,-4) {\Large $L$};
\node at (-1,-4) {\Large $a$};
\node at (-6.8,-9.05) {\Large $t$};
\node at (-9,-6.75) {\Large $t$};

\node at (-5,1) {\Large $t = 1, ~~a= 2,~100$};
\node at (-5.4,0.3) {\Large $E = 12,~~ \nu = 0,~0.3$};

\node[red] at (2.4,2.5) {\Large $\frac{P}{4}$};

\node[blue] at (3.3,1.5) {\Large $w$};

\node[] at (7,0,20) {\Large $X,u$};
\node[] at (5.9,-0.2,17.5) {\Large $Y,v$};
\node[] at (5.8,1.2,20) {\Large $Z,w$};

\end{tikzpicture}
    }~~~~
    }~
    \subfigure[Cantilever beam under bending moment (\cref{sec:num:2}) \label{fig:Numerical example_1.1b}]{\resizebox{0.5\textwidth}{!}{
        \includegraphics[]{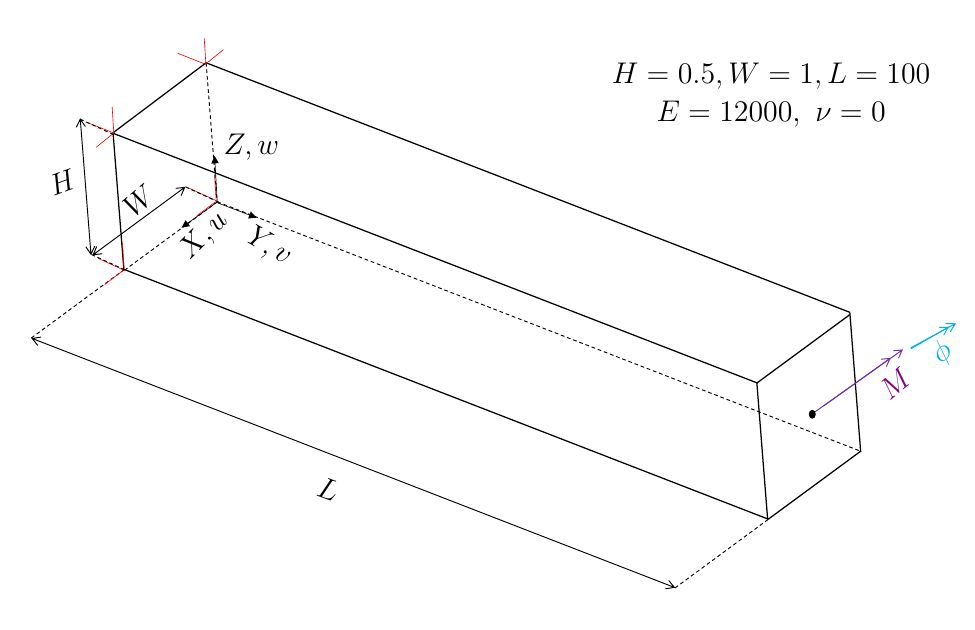}}}
    \caption{Geometries and material properties of cantilever beam examples}
    \label{fig:Numerical example_1.1}
\end{figure}
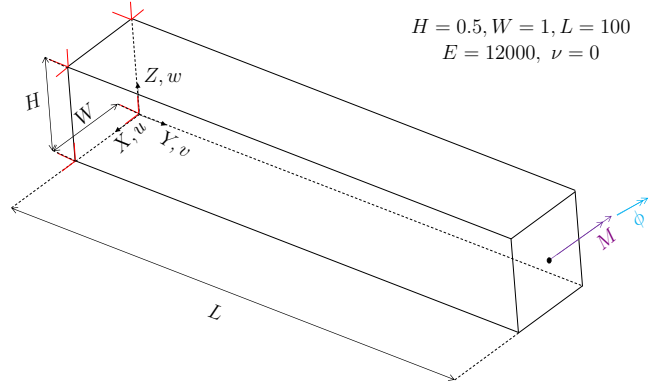

In this example, an elastic cantilever beam with a quadratic cross-section is loaded by a conservative shear force, which results in large deformations and small rotations. This example was also studied by the isoparametric solid-beam elements in \cite{Frischkorn.2013}. 
The St.~Venant-Kirchhoff hyperelastic material model is used. \Cref{fig:Numerical example_1.1a} shows the relevant geometric and material parameters. The beam is always discretized with three elements (knot spans), but two different aspect ratios of an element $a/t=a= 2,~100$ are considered. In accordance with the aspect ratio, the length of the beam is defined as $L=3a=6,~300$, i.e., the beam is either fairly thick or very slender.

\begin{figure}[t!]
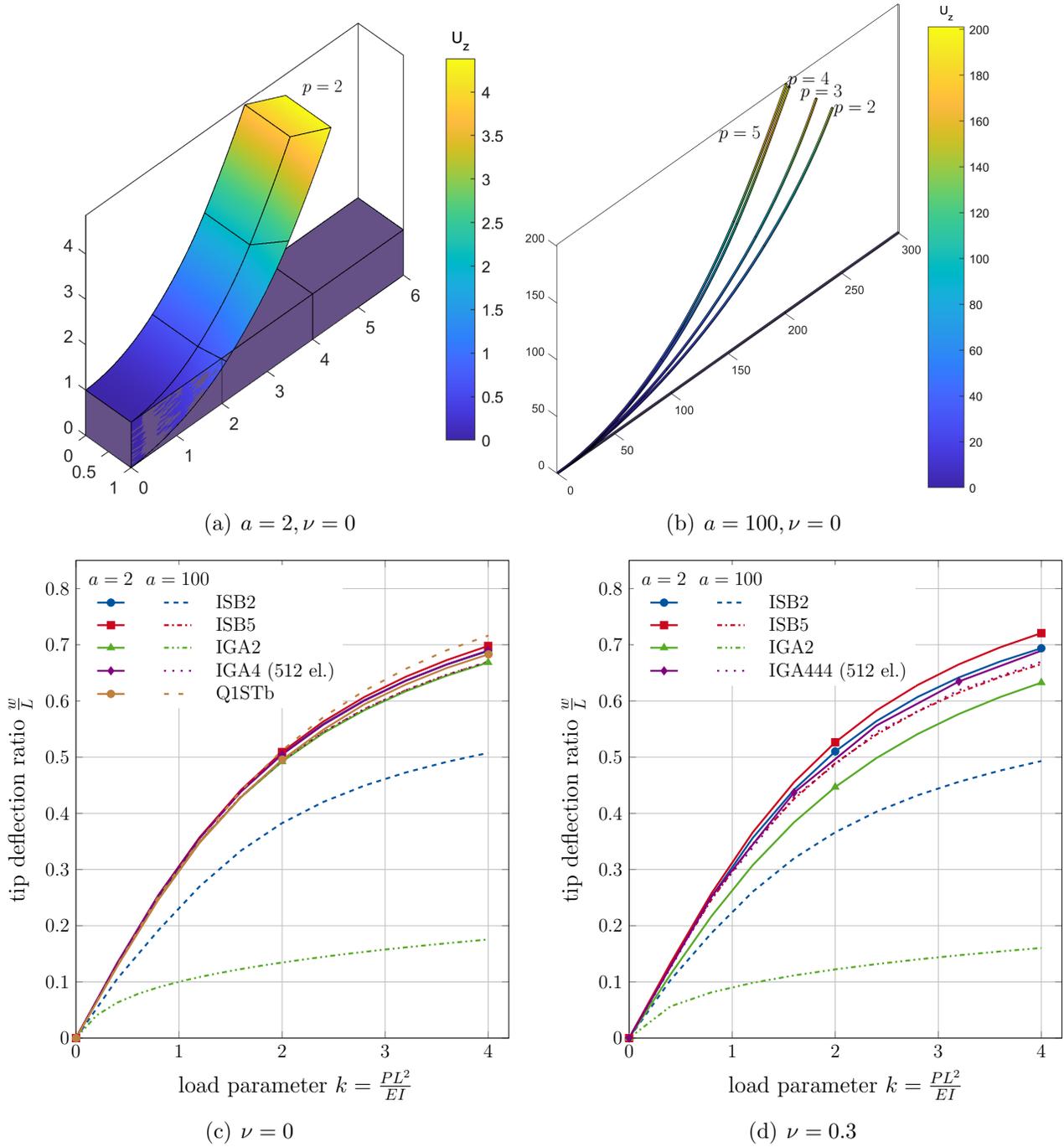

    \centering
    \subfigure[$a = 2,\, \nu = 0$  \label{fig:Numerical example_1.2a}]{\resizebox{0.40\textwidth}{!}{
        \input{defor_undef_at_2_diff_deg}}}
    \subfigure[$a = 100, \,\nu = 0$ \label{fig:Numerical example_1.2b}]{\resizebox{0.42\textwidth}{!}{
        \input{defor_undef_at_100_diff_deg}}}
    \subfigure[ $\nu = 0$\label{fig:Numerical example_1.2c}]{\resizebox{0.48\textwidth}{!}{
        \def\mystrut{\vphantom{hg}}

\pgfplotsset{
    legend image with text/.style={
        legend image code/.code={%
            \node[anchor=center] at (0.3cm,0cm) {#1};
        }
    },
    every axis plot/.append style={very thick},every mark/.append style={mark size=2.5pt}
}
\pgfplotstableread[col sep=comma]{shear_different_element_type.csv}{\loadL}

\begin{tikzpicture}[]
    \begin{axis}[height=25em, 
						width=23em, 
        legend style={at={(0.01,0.99)},anchor=north west,{draw=none},                font=\mystrut,
                legend cell align=left},
        legend columns = 2,
        scale=1.2,
        grid,
        xmin=0,
        xmax=4.2,
        ymin=0,
        ymax=0.85,
        xtick = {0,1,2,3,4,5},
        ytick = {0,0.1,...,0.8},
        xlabel=\text{ \large load parameter $k=\frac{PL^{2}}{EI}$},
        ylabel=\text{\large tip deflection ratio $\frac{w}{L}$},
        ]
        
         \addlegendimage{legend image with text= {$a=2$} };
        \addlegendentry{};
        
        \addlegendimage{legend image with text={$a=100$}};
        \addlegendentry{};

        \addplot[mark=*, color=blue, mark repeat = 5] table[x=load_solid_10, y=sb2_a_t_2_nu_0, col sep=comma]{\loadL};
        \addlegendentry{}
        
        \addplot[mark=none, dashed, color=blue] table[x=load_solid_10, y=sb2_a_t_100_nu_0, col sep=comma]{\loadL};
        \addlegendentry{ISB2}
        
        \addplot[mark=square*, color=red, mark repeat = 5] table[x=load_solid_10, y=sb5_a_t_2_nu_0, col sep=comma]{\loadL};
        \addlegendentry{}
        
        \addplot[mark=none, dashdotted, color=red] table[x=load_solid_10, y=sb5_a_t_100_nu_0, col sep=comma]{\loadL};
        \addlegendentry{ISB5}
        
        \addplot[mark=triangle*, color=green,  mark repeat=5] table[x=load_solid_10, y=solid_a_t_2, col sep=comma]{\loadL};
        \addlegendentry{}
        
        \addplot [mark=none, dashdotdotted, color=green] table[x=ls_a_t_100, y=s_a_t_100, col sep=comma]{\loadL};
        \addlegendentry{IGA2}
        
        \addplot[mark=diamond*, color=violet,  mark repeat=5] table[x=load_solid_10, y=solid4_512_a_t_2_nu_0, col sep=comma]{\loadL};
        \addlegendentry{}
        
        \addplot [mark=none, loosely dotted, color=violet] table[x=load_solid_10, y=solid4_512_a_t_100_nu_0, col sep=comma]{\loadL};
        \addlegendentry{IGA4 ($512$ el.)}
        
        \addplot[mark=otimes*, color=brown, mark repeat = 5] table[x=l_q1stb_2, y=q1stb_2, col sep=comma]{\loadL};
        \addlegendentry{}
        
        \addplot[mark=none, loosely dashed, color=brown] table[x=l_q1stb_100, y=q1stb_100, col sep=comma]{\loadL};
        \addlegendentry{Q1STb}
        
        
        
    \end{axis}
\end{tikzpicture}
    }}
    \subfigure[$\nu = 0.3$\label{fig:Numerical example_1.2d}]{\resizebox{0.48\textwidth}{!}{
\def\mystrut{\vphantom{hg}}

\pgfplotsset{
    legend image with text/.style={
        legend image code/.code={%
            \node[anchor=center] at (0.3cm,0cm) {#1};
        }
    },
    every axis plot/.append style={very thick},every mark/.append style={mark size=1.5pt}
}

\pgfplotstableread[col sep=comma]{shear_different_element_type.csv}{\loadL}

\begin{tikzpicture}[]
    \begin{axis}[height=25em, 
						width=23em, 
        legend style={at={(0.01,0.99)},anchor=north west,{draw=none},font=\mystrut,
                legend cell align=left},
        legend columns = 2,
        scale=1.2,
        grid,
        xmin=0,
        xmax=4.2,
        ymin=0,
        ymax=0.85,
        xtick = {0,1,2,3,4,5},
        ytick = {0,0.1,...,0.8},
        xlabel=\text{\large load parameter $k=\frac{PL^{2}}{EI}$},
        ylabel=\text{\large tip deflection ratio $\frac{w}{L}$},
        ]

    \addlegendimage{legend image with text= {$a=2$} }
        \addlegendentry{}
        
        \addlegendimage{legend image with text={$a=100$}}
        \addlegendentry{}

        \addplot[mark=*, color=blue,  mark repeat=5] table[x=load_solid_10, y=sb2_a_t_2_nu_0.3, col sep=comma]{\loadL};
        \addlegendentry{}
        
        \addplot[mark=none, dashed, color=blue] table[x=load_solid_10, y=sb2_a_t_100_nu_0.3, col sep=comma]{\loadL};
        \addlegendentry{ISB2}
        
        \addplot[mark=square*, color=red,  mark repeat=5] table[x=load_solid_10, y=sb5_a_t_2_nu_0.3, col sep=comma]{\loadL};
        \addlegendentry{}
        
        \addplot[mark=none, dashdotted, color=red] table[x=load_solid_10, y=sb5_a_t_100_nu_0.3, col sep=comma]{\loadL};
        \addlegendentry{ISB5}
        
        \addplot[mark=triangle*, color=green,  mark repeat=5] table[x=load_solid_10, y=solid_a_t_2_nu_0.3, col sep=comma]{\loadL};
        \addlegendentry{}
        
        \addplot [mark=none, dashdotdotted, color=green] table[x=load_solid_10, y=solid_a_t_100_nu_0.3, col sep=comma]{\loadL};
        \addlegendentry{IGA2}
        
        \addplot[mark=diamond*, color=violet,  mark repeat=2] table[x=load_step_5, y=solid4_512_a_t_2_nu_0.3, col sep=comma]{\loadL};
        \addlegendentry{}
        
        \addplot [mark=none, loosely dotted, color=violet] table[x=load_step_5, y=solid4_512_a_t_100_nu_0.3, col sep=comma]{\loadL};
        \addlegendentry{IGA444 ($512$ el.)}
        
        
        
    \end{axis}
\end{tikzpicture}
    }}
    \caption{Deformed configuration of beam under shear and comparison of normalized tip displacements obtained with different formulations with reference solutions}
    \label{fig:Numerical example_1.2}
\end{figure}

Deformed and undeformed configurations of the beam under shear for different aspect ratios $a=2$ and $a=100$ can be seen in \cref{fig:Numerical example_1.2a} and \cref{fig:Numerical example_1.2b}, respectively. For the slender beam  with $a=100$ in \cref{fig:Numerical example_1.2b}, the $p$-convergence of the ISB$p$ formulation can be observed. 
The results of the IGA solid-beam element (ISB$p$) are compared with the isoparametric solid-beam (Q1STb) from \cite{Frischkorn.2013} and the reference IGA solid elements (IGA$p$) as shown in \cref{fig:Numerical example_1.2c,fig:Numerical example_1.2d}. Additionally, the results of the IGA solid-beam (ISB5) are also plotted which converged to the highly refined reference solution of IGA solid element IGA4 with $512\times 1\times 1$ elements for aspect ratio of $a/t~=~100$ and different Poisson's ratios $\nu = 0$ and $\nu = 0.3$, while for aspect ratio $a/t~=~2$ the results are slightly better than the reference solution as the solid element still locks even in higher degree continuous function. The comparison shows that for the aspect ratio of $a/t~=~2$, ISB2 results are in accordance with the isoparametric solid-beam (Q1STb) element and the reference solution in \cref{fig:Numerical example_1.2c}.

\subsection{Cantilever beam under bending moment}
\label{sec:num:2}

To verify the IGA solid-beam formulation for large (finite) rotations and to alleviate transversal shear and membrane locking, a cantilever beam is considered which deforms into a circular arc by applying a pure bending moment. As the solid  and solid-beam formulations are based only on displacement degrees of freedom, the bending moment is applied by means of an equivalent follower load. \ref{subsec: Deformation dependent} shows a detailed description of applying this non-conservative loading.

\begin{figure}[tp!]
	\centering
        \subfigure[Deformed configuration of ISB1\label{fig:Numerical example_2.2a}]{\resizebox{0.48\textwidth}{!}{
        \input{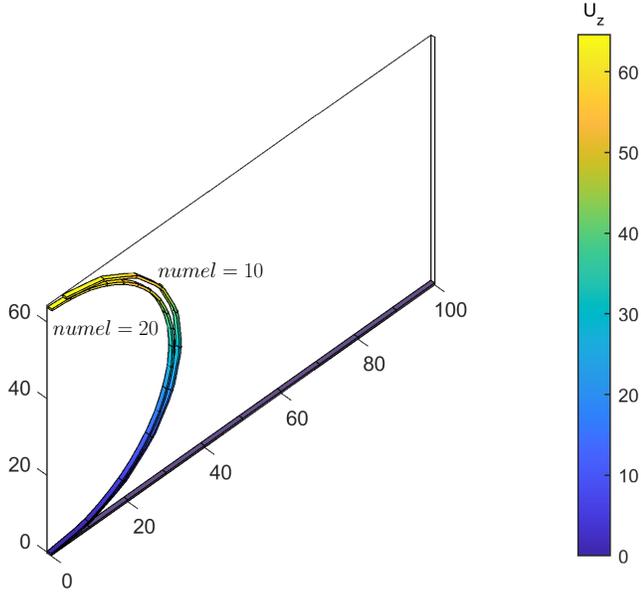}
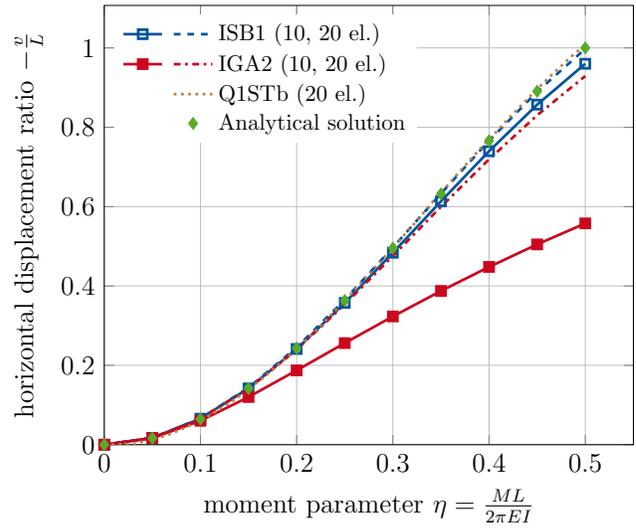
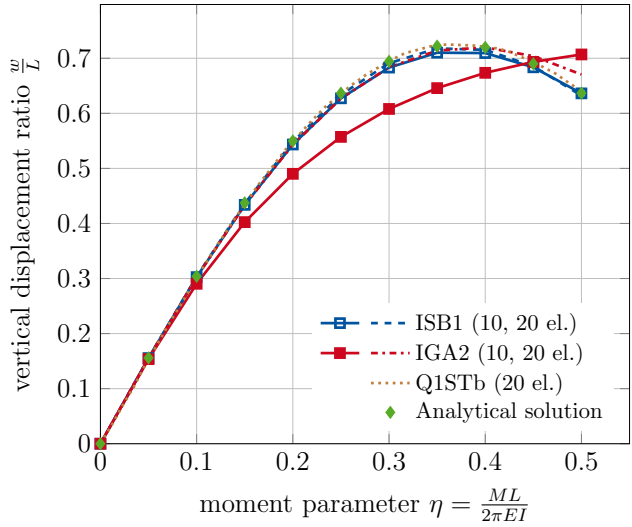
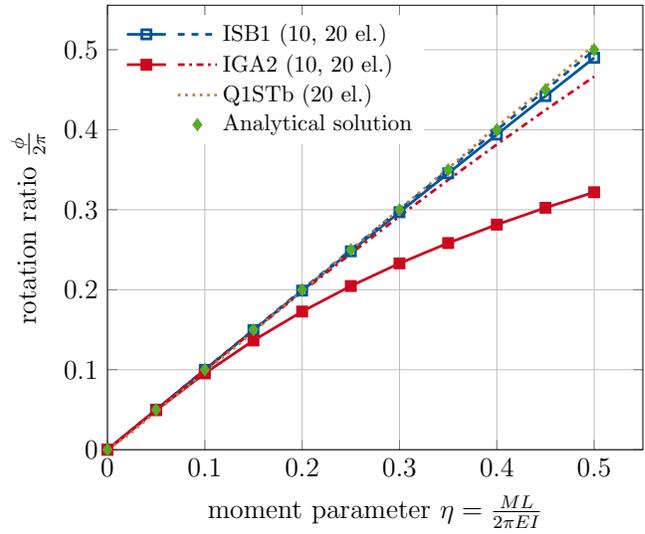}}
        \subfigure[Horizontal displacement response curve \label{fig:Numerical example_2.2b}]{\resizebox{0.49\textwidth}{!}{
\pgfplotstableread[col sep=comma]{horizontal_displacement_ratio.csv}{\loadL}

\pgfplotsset{
    every axis plot/.append style={very thick},every mark/.append style={mark size=2.5pt}
}

\begin{tikzpicture}[]
    \begin{axis}[legend columns=2,
            legend style={
                legend cell align=left,
            },
        scale=1.2,
        grid,
        xmin=0,
        ymin=0,
        xtick = {0,0.1,0.2,0.3,0.4,0.5},
        ytick = {0,0.2,0.4,0.6,0.8,1},
        legend style={at={(0.02,0.98)},anchor=north west,,{draw=none},font=\footnotesize},
        xlabel=\text{moment parameter $\eta=\frac{ML}{2\pi EI} $},
        ylabel=\text{horizontal displacement ratio $-\frac{v }{L }$}
                ]

       \addplot[mark=none, color=blue,mark=square] table[x=x_parameter, y=sb_q_1_numel_10, col sep=comma] {\loadL};
        \addlegendentry{}
        
        \addplot[mark=none, color=blue,dashed] table[x=x_parameter, y=sb_q_1_numel_20, col sep=comma] {\loadL};
        \addlegendentry{ISB1~(10,~20~el.)}
        
        \addplot[mark=none, color=red,mark=square*] table[x=x_parameter, y=s_q_2_numel_10, col sep=comma] {\loadL};
        \addlegendentry{}
        
        \addplot[mark=none, color=red,dashdotted] table[x=x_parameter, y=s_q_2_numel_20, col sep=comma] {\loadL};
        \addlegendentry{IGA2~(10,~20~el.)}
        
        \addlegendimage{empty legend}
        \addlegendentry{}
        
        \addplot[mark=none, color=brown,dotted] table[x=q1stb_numel_20_x, y=q1stb_numel_20_y, col sep=comma] {\loadL};
        \addlegendentry{Q1STb~(20~el.)}
        
        \addlegendimage{empty legend}
        \addlegendentry{}
        
        \addplot[only marks, mark=diamond*, green] table[x=x_parameter, y=ana_y, col sep=comma] {\loadL};
        \addlegendentry{Analytical solution}
        
    \end{axis}
\end{tikzpicture}}}
        \subfigure[Vertical displacement response curve \label{fig:Numerical example_2.2c}]{\resizebox{0.49\textwidth}{!}{
\pgfplotstableread[col sep=comma]{vertical_displacement_ratio.csv}{\loadL}

\pgfplotsset{
    every axis plot/.append style={very thick},every mark/.append style={mark size=2.5pt}
}

\begin{tikzpicture}[]
    \begin{axis}[legend columns=2,
            legend style={
                legend cell align=left,
            },
        scale=1.2,
        grid,
        xmin=0,
        ymin=0,
        xtick = {0,0.1,0.2,0.3,0.4,0.5},
        ytick = {0,0.1,0.2,0.3,0.4,0.5,0.6,0.7},
        legend style={at={(0.4,0.32)},anchor=north west,{draw=none},font=\footnotesize},
        xlabel=\text{moment parameter $\eta=\frac{ML}{2\pi EI} $},
        ylabel=\text{vertical displacement ratio $\frac{w }{L }$}
                ]

        \addplot[mark=none, color=blue,mark=square] table[x=x_parameter, y=sb_q_1_numel_10, col sep=comma] {\loadL};
        \addlegendentry{}
        
        \addplot[mark=none, color=blue,dashed] table[x=x_parameter, y=sb_q_1_numel_20, col sep=comma] {\loadL};
        \addlegendentry{ISB1~(10,~20~el.)}
        
        \addplot[mark=none, color=red,mark=square*] table[x=x_parameter, y=s_q_2_numel_10, col sep=comma] {\loadL};
        \addlegendentry{}
        
        \addplot[mark=none, color=red,dashdotted] table[x=x_parameter, y=s_q_2_numel_20, col sep=comma] {\loadL};
        \addlegendentry{IGA2~(10,~20~el.)}
        
        \addlegendimage{empty legend}
        \addlegendentry{}
        
        \addplot[mark=none, color=brown,dotted] table[x=q1stb_numel_20_x, y=q1stb_numel_20_y, col sep=comma] {\loadL};
        \addlegendentry{Q1STb~(20~el.)}
        
        \addlegendimage{empty legend}
        \addlegendentry{}
        
        \addplot[only marks, mark=diamond*, green] table[x=x_parameter, y=ana_y, col sep=comma] {\loadL};
        \addlegendentry{Analytical solution}
        
    \end{axis}
\end{tikzpicture}}}
	\subfigure[Rotation response curve \label{fig:Numerical example_2.2d}]{\resizebox{0.5\textwidth}{!}{

\pgfplotstableread[col sep=comma]{rotation_ratio.csv}{\loadL}

\pgfplotsset{
    every axis plot/.append style={very thick},every mark/.append style={mark size=2.5pt}
}

\begin{tikzpicture}[]
    \begin{axis}[legend columns=2,
            legend style={
                legend cell align=left,
            },
        scale=1.2,
        grid,
        xmin=0,
        ymin=0,
        xtick = {0,0.1,0.2,0.3,0.4,0.5},
        ytick = {0,0.1,0.2,0.3,0.4,0.5},
        legend style={at={(0.02,0.98)},anchor=north west,,{draw=none},font=\footnotesize},
        xlabel=\text{moment parameter $\eta=\frac{ML}{2\pi EI} $},
        ylabel=\text{rotation ratio $\frac{\phi }{2\pi }$}
                ]

        \addplot[mark=none, color=blue,mark=square] table[x=x_parameter, y=sb_q_1_numel_10, col sep=comma] {\loadL};
        \addlegendentry{}
        
        \addplot[mark=none, color=blue,dashed] table[x=x_parameter, y=sb_q_1_numel_20, col sep=comma] {\loadL};
        \addlegendentry{ISB1~(10,~20~el.)}
        
        \addplot[mark=none, color=red,mark=square*] table[x=x_parameter, y=s_q_2_numel_10, col sep=comma] {\loadL};
        \addlegendentry{}
        
        \addplot[mark=none, color=red,dashdotted] table[x=x_parameter, y=s_q_2_numel_20, col sep=comma] {\loadL};
        \addlegendentry{IGA2~(10,~20~el.)}
        
        \addlegendimage{empty legend}
        \addlegendentry{}
        
        \addplot[mark=none, color=brown,dotted] table[x=q1stb_numel_20_x, y=q1stb_numel_20_y, col sep=comma] {\loadL};
        \addlegendentry{Q1STb~(20~el.)}
        
        \addlegendimage{empty legend}
        \addlegendentry{}
        
        \addplot[only marks, mark=diamond*, green] table[x=x_parameter, y=ana_y, col sep=comma] {\loadL};
        \addlegendentry{Analytical solution}

    \end{axis}
\end{tikzpicture}}}
	\caption{Comparison of displacement ratios of beam under bending with different formulations}
	\label{fig:Numerical example_2.2}
\end{figure}
In \cref{fig:Numerical example_1.1b}, the geometry and the material parameters of the cantilever beam are shown. Again, the St.~Venant-Kirchhoff material model is used. The beam, which has an aspect ratio $L/H=200$, is meshed with 10 or 20 elements (knot spans) in the axial direction. Horizontal and vertical displacements, as well as rotation angle, are calculated at the center of the free-end cross-section of the beam. 
The deformed configurations of the beam and response curves for the normalized displacements and rotation angle with respect to the normalized moment can be seen in \cref{fig:Numerical example_2.2}. The results obtained with the linear ISB1 are compared with the  analytical 2D Timoshenko beam solution from \cite{E.Ramm1976}, quadratic IGA solid elements (IGA2), and  the linear, isoparametric solid-beam element (Q1STb). 
The results obtained from the ISB1 element are in very close agreement with the Q1STB element, as the linear NURBS basis function coincides  with the linear Lagrange shape functions, which serves as another verification of the ISB formulation.
For this slender beam, both are also in good agreement with the analytical solution, while locking can be observed for the IGA2 element, which underestimates the correct solution.

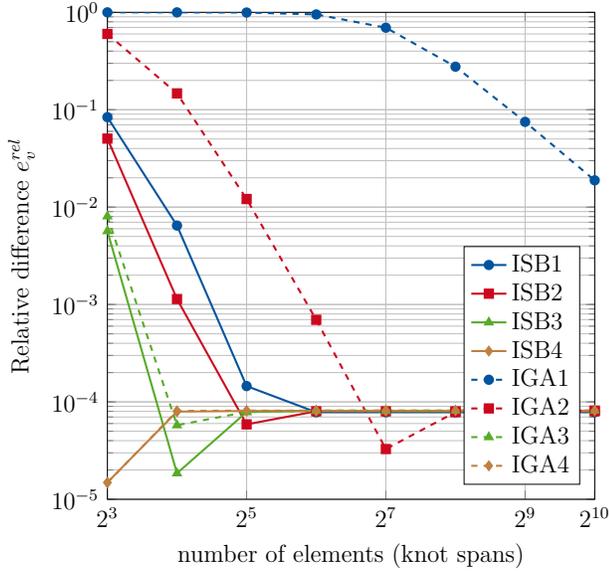
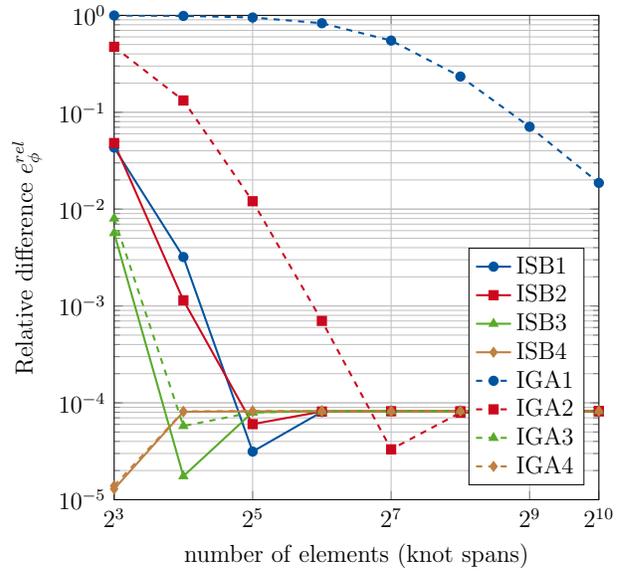
\begin{figure}[t]
    \centering\noindent
    \subfigure[Relative difference of horizontal displacement ratio \newline vs.\ number of elements \label{fig:Numerical example_2.3a}]{\resizebox{0.48\linewidth}{!}{
        \pgfplotstableread[col sep=comma]{Convergence_study_L_100_relative_diff_vs_numel.csv}{\loadLOT}

\pgfplotsset{
    every axis plot/.append style={very thick},every mark/.append style={mark size=2.5pt}
}

\begin{tikzpicture} 
	\begin{loglogaxis}[ height=24em, 
						width=24em, 
						grid=both,
						log basis x=2,
						xmin=8,
						xmax=1024, 
						xlabel={number of elements (knot spans)},
						xtick={8,32,128,512,1024},
						ymin=1e-5,
						ymax=1e0, 
						ylabel={Relative difference $e^{rel}_{v}$},
						legend style={
							cells={anchor=west},
							legend pos=south east, 
						}
					  ] 
	
		
		\addplot[color=blue, solid,thick,mark=*,mark options={solid}] table[x=numel,y=type_9_q_1_hor] {\loadLOT};
 		\addlegendentry{ISB1};
 		\addplot[color=red, solid,thick,mark=square*,mark options={solid}] table[x=numel,y=type_9_q_2_hor] {\loadLOT};
 		\addlegendentry{ISB2};
 		\addplot[color=green, solid,thick,mark=triangle*,mark options={solid}] table[x=numel,y=type_9_q_3_hor] {\loadLOT};
 		\addlegendentry{ISB3};
 		\addplot[color=brown, solid,thick,mark=diamond*,mark options={solid}] table[x=numel,y=type_9_q_4_hor] {\loadLOT};
 		\addlegendentry{ISB4};

   		\addplot[color=blue, dashed,thick,mark=*,mark options={solid}] table[x=numel,y=type_1_q_1_hor_pk_1] {\loadLOT};
 		\addlegendentry{IGA1};
 		\addplot[color=red, dashed,thick,mark=square*,mark options={solid}] table[x=numel,y=type_1_q_2_hor_pk_1] {\loadLOT};
 		\addlegendentry{IGA2};
 		\addplot[color=green, dashed,thick,mark=triangle*,mark options={solid}] table[x=numel,y=type_1_q_3_hor_pk_1] {\loadLOT};
 		\addlegendentry{IGA3};
 		\addplot[color=brown, dashed,thick,mark=diamond*,mark options={solid}] table[x=numel,y=type_1_q_4_hor_pk_1] {\loadLOT};
 		\addlegendentry{IGA4};
	\end{loglogaxis} 
\end{tikzpicture}}
    }~
    \subfigure[Relative difference of rotation ratio vs.\  number of \newline elements \label{fig:Numerical example_2.3b}]{\resizebox{0.48\linewidth}{!}{
        \pgfplotstableread[col sep=comma]{Convergence_study_L_100_relative_diff_vs_numel.csv}{\loadLOT}

\pgfplotsset{
    every axis plot/.append style={very thick},every mark/.append style={mark size=2.5pt}
}

\begin{tikzpicture} 
	\begin{loglogaxis}[ height=24em, 
						width=24em, 
						grid=both,
						log basis x=2,
						xmin=8,
						xmax=1024, 
						xlabel={number of elements (knot spans)},
						xtick={8,32,128,512,1024},
						ymin=1e-5,
						ymax=1e0, 
						ylabel={Relative difference $e^{rel}_{\phi}$},
						legend style={
							cells={anchor=west},
							legend pos=south east, 
						}
					  ] 
	
		
		\addplot[color=blue, solid,thick,mark=*,mark options={solid}] table[x=numel,y=type_9_q_1_phi] {\loadLOT};
 		\addlegendentry{ISB1};
 		\addplot[color=red, solid,thick,mark=square*,mark options={solid}] table[x=numel,y=type_9_q_2_phi] {\loadLOT};
 		\addlegendentry{ISB2};
 		\addplot[color=green, solid,thick,mark=triangle*,mark options={solid}] table[x=numel,y=type_9_q_3_phi] {\loadLOT};
 		\addlegendentry{ISB3};
 		\addplot[color=brown, solid,thick,mark=diamond*,mark options={solid}] table[x=numel,y=type_9_q_4_phi] {\loadLOT};
 		\addlegendentry{ISB4};

   		\addplot[color=blue, dashed,thick,mark=*,mark options={solid}] table[x=numel,y=type_1_q_1_phi_pk_1] {\loadLOT};
 		\addlegendentry{IGA1};
 		\addplot[color=red, dashed,thick,mark=square*,mark options={solid}] table[x=numel,y=type_1_q_2_phi_pk_1] {\loadLOT};
 		\addlegendentry{IGA2};
 		\addplot[color=green, dashed,thick,mark=triangle*,mark options={solid}] table[x=numel,y=type_1_q_3_phi_pk_1] {\loadLOT};
 		\addlegendentry{IGA3};
 		\addplot[color=brown, dashed,thick,mark=diamond*,mark options={solid}] table[x=numel,y=type_1_q_4_phi_pk_1] {\loadLOT};
 		\addlegendentry{IGA4};
	\end{loglogaxis} 
\end{tikzpicture}}
    }
    \caption{Convergence study of tip displacement and rotation of beam under bending}
    \label{fig:Numerical example_2.3}
\end{figure}

A convergence study comparing ISB$p$ and IGA$p$ is performed in \cref{fig:Numerical example_2.3} by increasing the degree of the NURBS basis functions and also the number of elements. 
The convergence study shows that the lower degree isogeometric solid-beam element requires fewer elements (degrees of freedom) to converge to the analytical solution than the IGA solid element, which suffers from artificial stiffening behaviour (transversal shear and membrane locking).
Note that the relative errors are plateauing due to the conceptual difference between the analytical beam model and the 3D solid models.

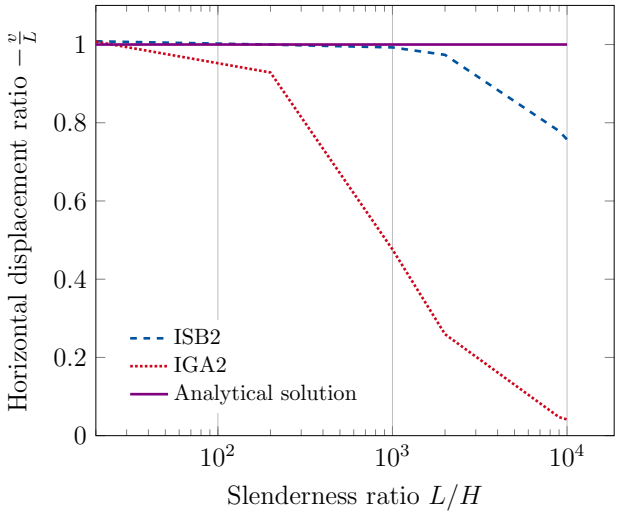
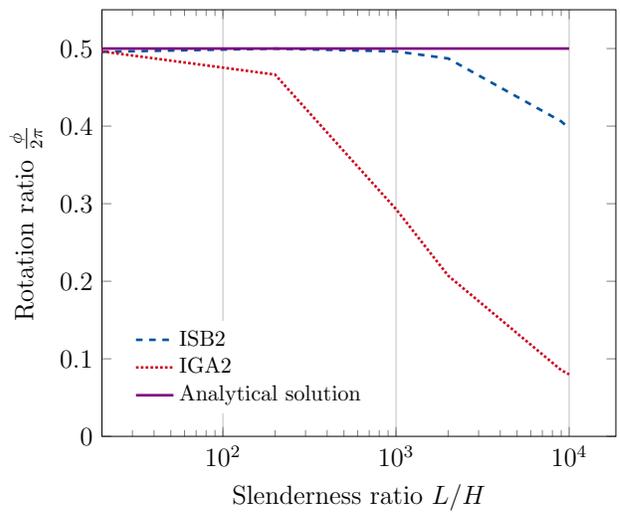
\begin{figure}[t!]
    \centering
    \subfigure[Horizontal displacement ratio vs.\ slenderness ratio \label{fig:Numerical example_2.4a}]{\resizebox{0.48\linewidth}{!}{
        \pgfplotstableread[col sep=comma]{comparison_of_length_of_elements.csv}{\loadL}

\pgfplotsset{
    every axis plot/.append style={very thick},every mark/.append style={mark size=2.5pt}
}

\begin{tikzpicture}[]
    \begin{axis}[
        scale = 1.2,
        xmajorgrids,
        xmin=20,
        ymin=0,
        legend style={at={(0.05,0.27)},anchor=north west,{draw=none}},
        legend style={font=\footnotesize,legend cell align=left},
        xmode=log,
        ytick = {0,0.2,0.4,0.6,0.8,1},
        ymax = 1.1,
        xlabel=\text{Slenderness ratio~$L/H$},
        ylabel=\text{Horizontal displacement ratio $-\frac{v}{L}$}
                ]

        \addplot[blue,mark=none,dashed] table[x=sl_ra, y=hor_ratio_sb, col sep=comma] {\loadL};
        \addlegendentry{ISB2}
        
        \addplot[red,densely dotted] table[x=sl_ra, y=hor_ratio_s, col sep=comma] {\loadL};   
        \addlegendentry{IGA2}

        \addplot[violet,mark=none] coordinates {(20, 1) (10000, 1)};
        \addlegendentry{Analytical solution}
        
    \end{axis}
\end{tikzpicture}}
    }
    \subfigure[Rotation ratio vs.\ slenderness ratio \label{fig:Numerical example_2.4b}]{\resizebox{0.48\linewidth}{!}{
        \pgfplotstableread[col sep=comma]{comparison_of_length_of_elements.csv}{\loadL}

\pgfplotsset{
    every axis plot/.append style={very thick},every mark/.append style={mark size=2.5pt}
}

\begin{tikzpicture}[]
    \begin{axis}[
        scale = 1.2,
        xmajorgrids,
        xmin=20,
        ymin=0,
        ymax=0.55,
        legend style={at={(0.05,0.27)},anchor=north west,{draw=none}},
        legend style={font=\footnotesize, legend cell align=left},
        xmode=log,
        ytick = {0,0.1,0.2,0.3,0.4,0.5},
        xlabel=\text{Slenderness ratio~$L/H$},
        ylabel=\text{Rotation ratio $\frac{\phi }{2\pi }$}
                ]

        \addplot[blue,mark=none,dashed] table[x=sl_ra, y=phi_ratio_sb, col sep=comma] {\loadL};
        \addlegendentry{ISB2}
        
        \addplot[red,densely dotted] table[x=sl_ra, y=phi_ratio_s, col sep=comma] {\loadL};   
        \addlegendentry{IGA2}

        \addplot[violet,mark=none] coordinates {(20, 0.5) (10000, 0.5)};
        \addlegendentry{Analytical solution}
        
    \end{axis}
\end{tikzpicture}}
    }
    \caption{Study of a cantilever beam under bending for varying slenderness ratio}
    \label{fig:Numerical example_2.4}
\end{figure}

Finally, for this beam bending benchmark example, the IGA solid-beam element ISB2 is tested with different slenderness ratios and is compared with the IGA solid element IGA2. In \cref{fig:Numerical example_2.4}, the results for different slenderness ratios $L/H$ ranging from 10 to 10,000 with the discretization of 20 elements (knot span) along the beam direction are presented. It can be seen that the ISB2 element is able to alleviate the locking effects even for larger slenderness ratios when the solid element IGA2 completely fails to obtain the correct solution. 

\subsection{Bathe's circular arc beam}
\label{sec:num:3}

This classical benchmark problem is tested to validate the bending-torsion behaviour of the proposed IGA solid-beam formulation. It was first introduced by Bathe and Belourchi \cite{Bathe.1979}. The geometry of the 45$^{\circ}$ circular arc cantilever beam with unit square cross-section subjected to an out-of-plane conservative load $F$ can be seen in \cref{fig:Numerical example_3.2a}. 
Now, two different cases of isotropic material behaviours are considered: the linearly elastic St.~Venant-Kirchhoff (SVK) material without Poisson effect ($E = 10^{7},~\nu~=~0$) and the compressible Neo-Hooke (NH) material with Poisson effect ($E = 10^{7},~\nu~=~0.3$). The load is also applied in two different levels, $F = 300$ and 600 units. The beam is meshed with 8 elements (knot spans) along the beam direction with NURBS of degree $p = 2$ and the load is applied in 10 equidistant load steps in both material cases. 

\subsubsection{Case 1: Saint-Venant Kirchoff material without Poisson effect}

\begin{figure}[H]
    \centering
    \subfigure[Undeformed and deformed configurations of 45$^{\circ}$-arc cantilever beam \label{fig:Numerical example_3.2a}]{\resizebox{0.47\linewidth}{!}{
        \input{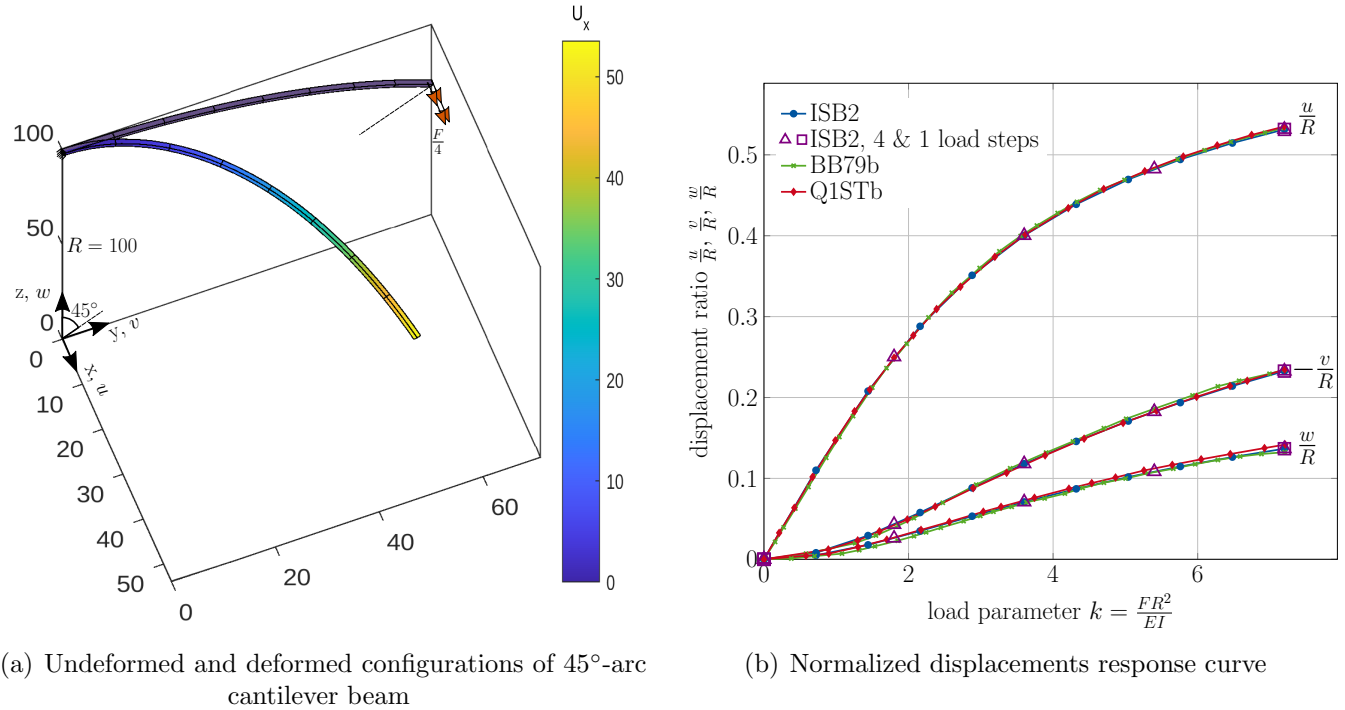}}
    }~
    \subfigure[Normalized displacements response curve \label{fig:Numerical example_3.2b}]{\resizebox{0.50\linewidth}{!}{
        \pgfplotstableread[col sep=comma]{600_units_comparison.csv}{\loadFyL}

\pgfplotsset{
    every axis plot/.append style={very thick},every mark/.append style={mark size=2.5pt}
}

\begin{tikzpicture}[]
\tikzstyle{every node}=[font=\Large]
    \begin{axis}[ legend columns=1,
            legend style={
                legend cell align=left
            },
        num1/.style={only marks,mark options={mark size=4pt},color=violet, mark=triangle},
        num2/.style={only marks,mark options={mark size=3pt},color=violet, mark=square},
        combo legend/.style={
          legend image code/.code={
            \draw [/pgfplots/num1] plot coordinates {(1mm,0cm)};
            \draw [/pgfplots/num2] plot coordinates {(5mm,0cm)};
          }
        },
        scale=2,
        grid,
        xmin=0,
        ymin=0,
        xtick = {0,2,4,6,8},
        ytick = {0,0.1,0.2,0.3,0.4,0.5,0.6},
        legend style={at={(0.02,0.97)},anchor=north west,{draw=none}},
        xlabel=\text{load parameter $k =\frac{FR^{2}}{EI}$},
        ylabel=\text{displacement ratio $\frac{u }{R },~\frac{v }{R },~\frac{w}{R }$}
                ]
         \addlegendimage{color=blue,very thick, mark = *}\addlegendentry{ISB2}
        \addlegendimage{combo legend}\addlegendentry{ISB2, 4 $\&$ 1 load steps}
        \addlegendimage{color=green,mark = x,very thick}\addlegendentry{BB79b}
        \addlegendimage{color=red,mark = diamond*,very thick}\addlegendentry{Q1STb}

        \addplot[color=blue, very thick, mark = *] table[x=load_par, y=I-SB2_u, col sep=comma] {\loadFyL};
        
        \addplot[color=blue,very thick, mark = *] table[x=load_par, y=I-SB2_v, col sep=comma] {\loadFyL};
        
        \addplot[color=blue,very thick, mark = *] table[x=load_par, y=I-SB2_w, col sep=comma] {\loadFyL};
        
        \addplot[color=green,mark = x,very thick] table[x=BB79b_u_x, y=BB79b_u_y, col sep=comma] {\loadFyL};
        
        \addplot[color=green,mark = x,very thick] table[x=BB79b_v_x, y=BB79b_v_y, col sep=comma] {\loadFyL};
        
        \addplot[color=green,mark = x,very thick] table[x=BB79b_w_x, y=BB79b_w_y, col sep=comma] {\loadFyL};

        \addplot[color=red,mark = diamond*,very thick] table[x=Q1STb_u_x, y=Q1STb_u_y, col sep=comma] {\loadFyL};
        
        \addplot[color=red,mark = diamond*,very thick] table[x=Q1STb_v_x, y=Q1STb_v_y, col sep=comma] {\loadFyL};
        
        \addplot[color=red,mark = diamond*,very thick] table[x=Q1STb_w_x, y=Q1STb_w_y, col sep=comma] {\loadFyL};

        \addplot[only marks,mark options={mark size=5pt},color=violet, mark=triangle] table[x=ls_4_x, y=ls_4_u_y, col sep=comma] {\loadFyL};

        \addplot[only marks,mark options={mark size=5pt},color=violet, mark=triangle] table[x=ls_4_x, y=ls_4_v_y, col sep=comma] {\loadFyL};

        \addplot[only marks,mark options={mark size=5pt},color=violet, mark=triangle] table[x=ls_4_x, y=ls_4_w_y, col sep=comma] {\loadFyL};

        \addplot[only marks,mark options={mark size=4pt},color=violet, mark=square] table[x=ls_1_x, y=ls_1_u_y, col sep=comma] {\loadFyL};

        \addplot[only marks,mark options={mark size=4pt},color=violet, mark=square] table[x=ls_1_x, y=ls_1_v_y, col sep=comma] {\loadFyL};

        \addplot[only marks,mark options={mark size=4pt},color=violet, mark=square] table[x=ls_1_x, y=ls_1_w_y, col sep=comma] {\loadFyL};
        
        \node [below right] at (rel axis cs:0.92,0.28) {\LARGE $\frac{w }{R }$};
        \node [below right] at (rel axis cs:0.91,0.44) {\LARGE $-\frac{v }{R }$};
        \node [above right] at (rel axis cs:0.92,0.88) {\LARGE $\frac{u }{R }$};
        
    \end{axis}
\end{tikzpicture}}
    }
    \caption{Bending-torsion of arc beam for $F$ = 600}
    \label{fig:Numerical example_3.2}
\end{figure}

The deformed configuration for an out-of-plane load $F = 600$ can be seen in \cref{fig:Numerical example_3.2a}. In \Cref{fig:Numerical example_3.2b}, the response curve showing the normalized displacements evaluated at the center of free end is plotted over the normalized load parameter. As can be seen, the results of the ISB2 elements are in accordance with the beam theory \cite{Bathe.1979} and the isoparametric solid-beam element (Q1STb) \cite{Frischkorn.2013}. Additionally, to demonstrate the robustness of the formulation, the load is applied in only four load steps and even only one load step. As can also be seen in \cref{fig:Numerical example_3.2b}, the results obtained even with only one load step are the same as with 10 equidistant load steps.

\Cref{table:Numerical_table_3.1} provides a comparison of the tip deflections obtained with the ISB formulation for $F=300,~600$ with other beam element formulations from literature. The results of the other structural beam formulations are taken from \cite{Frischkorn.2013}, except for the beam element proposed by Choi \textit{et al} \cite{Choi.2023}.
Generally, a good agreement can be observed. 
Furthermore, \cref{table:Numerical_table_3.2} details the number of load steps and iterations per load step required to obtain the results presented in \cref{table:Numerical_table_3.1}. The robust behaviour of the ISB formulation can be observed, as it can converge with large load steps and requires only a few Newton iterations. The number of iterations required to solve the problem is provided with and without employing the MIP method (i.e., ISB2 with only ANS and EAS method), which shows that the MIP method is crucial to improve robustness and helps to significantly decrease the computational time. 
   
\begin{table}[t]
\centering
\renewcommand{\arraystretch}{1.1}
\begin{tabular}{lccccccc}
\hline
\multicolumn{1}{l}{Element} & \multicolumn{3}{c}{$F = 300$}      &  & \multicolumn{3}{c}{$F = 600$} \\ \cline{2-4} \cline{6-8} 
                            & $v$        & $w$      & $u$      &  & ${v}$     & ${w}$    & ${u}$    \\ \hline
IGA solid-beam (ISB2)                & -11.81     & 7.08     & 40.02    &  & -23.26   & 13.67  & 53.14  \\
Frischkorn and Reese (Q1STb) \cite{Frischkorn.2013}        & -11.71     & 7.27     & 40.06    &  & -23.38  & 14.11  & 53.50   \\
Choi \textit{et al.} (\text{p=4},~numel =80) \cite{Choi.2023}        & \multicolumn{3}{c}{Not reported} &  & -23.56  & 13.60  & 53.48   \\
Bathe and Bolourchi \cite{Bathe.1979}          & -11.50      & 6.80      & 39.50     &  & -23.50   & 13.40   & 53.40   \\
Rhim and Lee \cite{Rhim.1998}                & -12.03     & 7.13     & 40.31    &  & -23.64  & 13.70   & 53.46  \\
Crisfield \cite{M.A.Crisfield.1990}                 & -12.18     & 7.13     & 40.53    &  & -23.87  & 13.69  & 53.71  \\
Simo and Vu-Quoc \cite{Simo.1986}           & -11.87     & 6.96     & 40.08    &  & -23.48  & 13.50   & 53.37  \\
Cardona and Geradin \cite{Cardona.1988}        & -12.07     & 7.15     & 40.35    &  & -23.67  & 13.74  & 53.50   \\
Dvorkin \textit{et al.} \cite{Dvorkin.1988}               & \multicolumn{3}{c}{Not reported} &  & -23.50   & 13.60   & 53.30   \\ \hline
\end{tabular}

\caption{Comparison with structural beam elements--tip deflections for two different load levels}
\label{table:Numerical_table_3.1}
\end{table}

\begin{table}[t]
\centering
\renewcommand{\arraystretch}{1.1}
\begin{tabular}{lccccc}
\hline
Element                             & \multicolumn{2}{c}{${F} = 300$}      & & \multicolumn{2}{c}{${F} = 600$}   \\ \cline{2-3} \cline{5-6} 
                                    & \# load steps   & \# iter.  & & \# load steps & \# iter. \\ \hline
IGA solid-beam (ISB2)               & 1               & 7             & & 1             & 8             \\
IGA solid-beam (ISB2, without MIP) & 1               & 15           &  & 1             & 21            \\
Frischkorn and Reese (Q1STb) \cite{Frischkorn.2013}         & 1               & 14            & & 1             & 19            \\
Choi \textit{et al.} \cite{Choi.2023}                         & \multicolumn{2}{c}{Not reported} & & 1             & 8             \\
Bathe and Bolourchi \cite{Bathe.1979}                & \multicolumn{5}{c}{60 equal increments}                          \\
Rhim and Lee \cite{Rhim.1998}                       & 1               & 8             & & 1             & 13            \\
Crisfield \cite{M.A.Crisfield.1990}                         & 1               & 9            &  & 3             & 17            \\
Simo and Vu-Quoc \cite{Simo.1986}                   & 1               & 13           &  & 3             & 27            \\
Cardona and Geradin \cite{Cardona.1988}                & \multicolumn{5}{c}{6 equal increments, total 47 iterations}      \\
Dvorkin \textit{et al.}  \cite{Dvorkin.1988}                     & \multicolumn{5}{c}{10 equal increments, total 34 iterations}     \\ \hline
\end{tabular}

\caption{Comparison with structural beam elements--number of load steps and Newton iterations required for two different load levels}
\label{table:Numerical_table_3.2}
\end{table}

\subsubsection{Case 2: Compressible Neo-Hookean  material with Poisson effect}

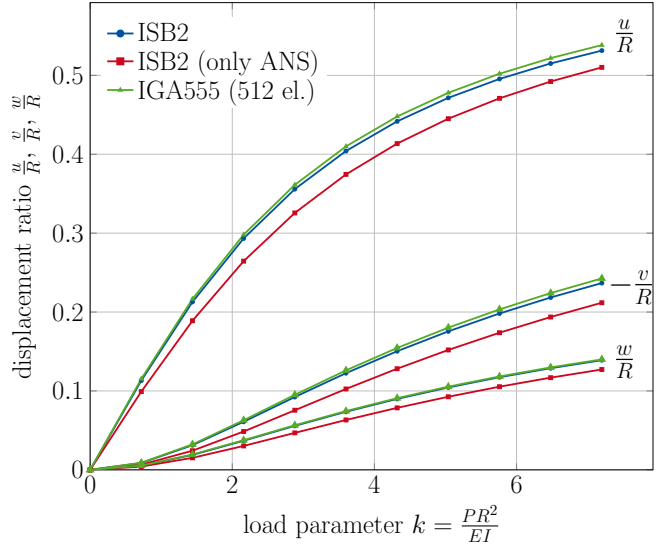
\begin{figure}[ht!]
    \centering
    \resizebox{0.5\textwidth}{!}{
        \pgfplotstableread[col sep=comma]{comparison_ANS_ANSEAS_element.csv}{\loadFyL}

\pgfplotsset{
    every axis plot/.append style={very thick},every mark/.append style={mark size=2.5pt}
}

\begin{tikzpicture}[]
\tikzstyle{every node}=[font=\Large]
    \begin{axis}[ legend columns=1,
            legend style={
                legend cell align=left
            },
        num1/.style={only marks,mark options={mark size=2pt},color=magenta, mark=triangle},
        num2/.style={only marks,mark options={mark size=3pt},color=magenta, mark=square},
        combo legend/.style={
          legend image code/.code={
            \draw [/pgfplots/num1] plot coordinates {(1mm,0cm)};
            \draw [/pgfplots/num2] plot coordinates {(5mm,0cm)};
          }
        },
        scale=1.88,
        grid,
        xmin=0,
        ymin=0,
        xtick = {0,2,4,6,8},
        ytick = {0,0.1,0.2,0.3,0.4,0.5,0.6},
        legend style={at={(0.02,0.97)},anchor=north west,{draw=none}},
        xlabel=\text{load parameter $k =\frac{PR^{2}}{EI}$},
        ylabel=\text{displacement ratio $\frac{u }{R },~\frac{v }{R },~\frac{w}{R }$}
                ]
        \addlegendimage{color=blue, very thick, mark = *}\addlegendentry{ISB2}
        \addlegendimage{color=red,mark = square*,very thick}\addlegendentry{ISB2 (only ANS)}

         \addlegendimage{color=green,mark = triangle*,very thick, mark size=1}\addlegendentry{IGA555 ($512$ el.)}
        
        \addplot[color=blue,mark = *,very thick, mark size=1] table[x=moment_steps, y=type_4_p1_q_2_numel_8_x, col sep=comma] {\loadFyL};
        
        \addplot[color=blue,mark = *,very thick, mark size=1] table[x=moment_steps, y=type_4_p1_q_2_numel_8_y, col sep=comma] {\loadFyL};
        
        \addplot[color=blue,mark = *,very thick, mark size=1] table[x=moment_steps, y=type_4_p1_q_2_numel_8_z, col sep=comma] {\loadFyL};

        \addplot[color=red, very thick, mark = square*, mark size=1] table[x=moment_steps, y=type_2_p1_q_2_numel_8_x, col sep=comma] {\loadFyL};
        
        \addplot[color=red,very thick, mark = square*, mark size=1] table[x=moment_steps, y=type_2_p1_q_2_numel_8_y, col sep=comma] {\loadFyL};
        
        \addplot[color=red,very thick, mark = square*, mark size=1] table[x=moment_steps, y=type_2_p1_q_2_numel_8_z, col sep=comma] {\loadFyL};

        
        

        \addplot[color=green, very thick, mark = triangle*, mark size=1] table[x=moment_steps, y=type_1_p5_q_5_numel_512_x, col sep=comma] {\loadFyL};
        
        \addplot[color=green,very thick, mark = triangle*] table[x=moment_steps, y=type_1_p5_q_5_numel_512_y, col sep=comma] {\loadFyL};
        
        \addplot[color=green,very thick, mark = triangle*] table[x=moment_steps, y=type_1_p5_q_5_numel_512_z, col sep=comma] {\loadFyL};

        
        \node [below right] at (rel axis cs:0.92,0.28) {\LARGE $\frac{w }{R }$};
        \node [below right] at (rel axis cs:0.91,0.44) {\LARGE $-\frac{v }{R }$};
        \node [above right] at (rel axis cs:0.92,0.88) {\LARGE $\frac{u }{R }$};
        
    \end{axis}
\end{tikzpicture}
    }
    \caption{Comparison of results with different formulations}
    \label{fig:Numerical example_3.1.1}
\end{figure}

To see the Poisson thickness locking effect more prominently, the NH material model with Poisson's ratio $\nu~=~0.3$ is considered. The boundary conditions and geometry are the same as in the previous case and the beam is again discretized with 8 elements along the axial direction. 

\Cref{fig:Numerical example_3.1.1} shows the tip displacement of the free end of the beam at the center. Tip displacements of the highly refined reference solid element IGA555 with $512\times 1\times 1$ elements are compared with ISB2 as well as ISB2  without the EAS method (only ANS). The results show that the Poisson thickness locking gets alleviated and that the results are significantly improved with the EAS method. 
This can also be seen from \cref{table:Numerical_table_3.2.1}. 
The table also includes results for the IGA5 element with degree $q=1$ in the cross-section and  $1024\times 1\times 1$ elements, which exhibit artificial stiffening and underestimate the true solution, similar to ISB2 without EAS. This shows that Poisson thickness (material) locking cannot be improved by refinement along the beam direction and underpins the necessity of the EAS method when using low degrees in the cross-section.

\begin{table}[H]
\centering
\renewcommand{\arraystretch}{1.1}
\begin{tabular}{lccccccc}
\hline
\multicolumn{1}{l}{Element} & \multicolumn{3}{c}{$F = 300$}      &  & \multicolumn{3}{c}{$F = 600$} \\ \cline{2-4} \cline{6-8} 
                            & ${w}$        & ${v}$      & ${u}$      &  & ${w}$     & ${v}$    & ${u}$    \\ \hline
ISB2                                & 7.3628   & -12.2516   & 40.4085 & & 13.8976 & -23.6708 & 53.1366 \\
ISB2 (only ANS, without EAS)        & 6.3197   & -10.2504   & 37.4453  & & 12.7167 & -21.1791 & 51.0090  \\
IS5 ($1024\times 1\times 1$ el.) &  6.2371        &    -10.5217        &    37.8288      & & 12.5818 & -21.6890  & 51.5366 \\
IS555 ($512\times 1\times 1$ el.)  &     7.4503     &    -12.6047        &  40.9935        & & 13.9976 & -24.2655 & 53.8399 \\
Choi \textit{et al}  \cite{Choi.2023}                         & \multicolumn{3}{c}{Not reported} & & 13.8148 & -23.9790 & 53.6901 \\ \hline
\end{tabular}%
\caption{Comparison of tip deflections for two different load levels with reference elements}
\label{table:Numerical_table_3.2.1}
\end{table}

\subsection{Lateral buckling of right-angle frame}
\label{sec:num:4}

\begin{figure}[H]
    \centering
    \subfigure[Geometry and boundary conditions  \label{fig:Numerical example_4.1a}]{\resizebox{0.37\linewidth}{!}{
        \input{right_angle_frame_image}}
    }
    \subfigure[Undeformed and deformed configurations (ISB3) \label{fig:Numerical example_4.1b}]{\resizebox{0.55\linewidth}{!}{
        \includegraphics[]{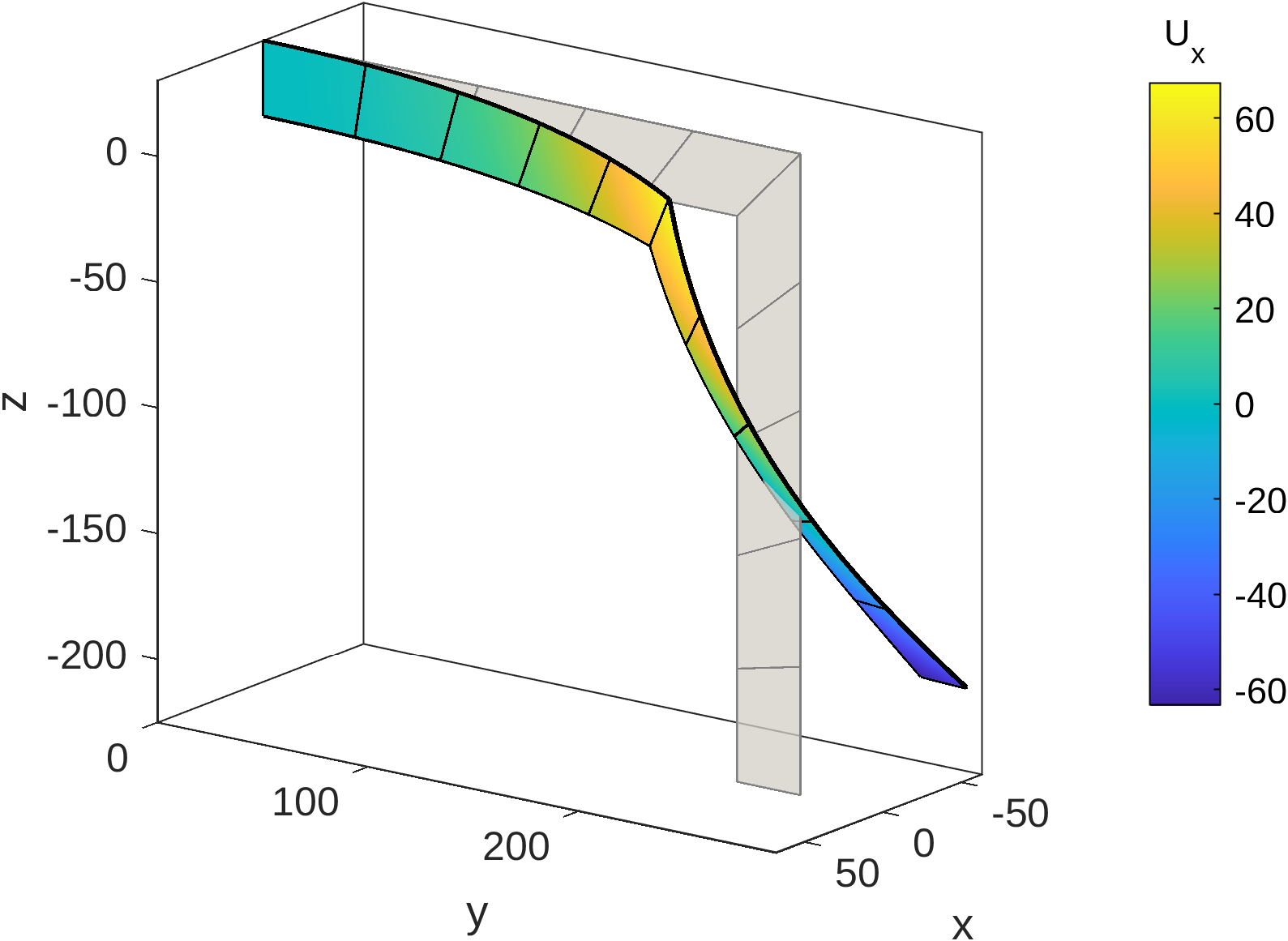}}
    }
    \subfigure[Load-displacement response curves \label{fig:Numerical example_4.1c}]{\resizebox{0.5\linewidth}{!}{
        \def\mystrut{\vphantom{hg}}

\pgfplotsset{
    legend image with text/.style={
        legend image code/.code={%
            \node[anchor=center] at (0.3cm,0cm) {#1};
        }
    },
    every axis plot/.append style={very thick},every mark/.append style={mark size=2.5pt}

}

\pgfplotstableread[col sep=comma]{right_angle_lateral_buckling_formulation_comp.csv}{\loadFyL}

\begin{tikzpicture}[]
\pgfplotsset{every axis plot/.style={}}=[font=\Large]
    \begin{axis}[ legend columns=2,
            legend style={
                font=\mystrut,
                legend cell align=left,
            },
        width=26em,
        height=24em,
        grid,
        xmin=0,
        ymin=0.4,
        xmax = 60,
        ymax = 1.6,
        xtick = {0,10,20,30,40,50,60},
        ytick = {0.4,0.6,0.8,1,1.2,1.4,1.6},
        legend style={cells={anchor=west},legend pos=south east,{draw=none}},
        xlabel=\text{\text{out-of-plane displacement at the center of free end}},
        ylabel=\text{load $P$}
                ]
        
        \addlegendimage{legend image with text=single~~}
        \addlegendentry{}
        
        \addlegendimage{legend image with text=multi-}
        \addlegendentry{}
        
        \addlegendimage{legend image with text=patch}
        \addlegendentry{}

        \addlegendimage{legend image with text=patch}
        \addlegendentry{}
        
        \addplot[color=blue,very thick, mark=*, mark repeat=10] table[x=sb_single_2_numel_10_x_ls_100, y=sb_single_3_numel_10_y_ls_100, col sep=comma] {\loadFyL};
         \addlegendentry{};
         
        \addplot[color=blue,dashed,very thick] table[x=sb_multi_2_numel_10_x_ls_100, y=sb_multi_3_numel_10_y_ls_100, col sep=comma] {\loadFyL};
        \addlegendentry{ISB2};
        
         \addplot[color=red,very thick, mark=square*, mark repeat=10] table[x=sb_single_3_numel_10_x_ls_100, y=sb_single_3_numel_10_y_ls_100, col sep=comma] {\loadFyL};
         \addlegendentry{};
         
        \addplot[color=red,dashdotted,very thick, mark = none] table[x=sb_multi_3_numel_10_x_ls_100, y=sb_multi_3_numel_10_y_ls_100, col sep=comma] {\loadFyL};
        \addlegendentry{ISB3};

          \addplot[color=green,very thick,  mark=triangle*, mark repeat=10] table[x=sb_single_4_numel_10_x_ls_100, y=sb_single_4_numel_10_y_ls_100, col sep=comma] {\loadFyL};
          \addlegendentry{};
          
          \addplot[color=green,dotted,very thick, mark = none] table[x=sb_multi_4_numel_10_x_ls_100, y=sb_multi_4_numel_10_y_ls_100, col sep=comma] {\loadFyL};
          \addlegendentry{ISB4};

        \addlegendentry{};
        \addplot[color=brown,very thick, , mark=diamond*, mark repeat=7] table[x=simo_X, y=simo_Y, col sep=comma] {\loadFyL};
        
        \addlegendimage{empty legend};
        \addlegendentry{Simo \& Vu-Quoc  \cite{Simo.1986}};

        \addlegendentry{};
        \addplot[color=violet,very thick,  mark=asterisk, mark repeat=2] table[x=Smolenski_X, y=Smolenski_Y, col sep=comma] {\loadFyL};
        \addlegendimage{empty legend}
        \addlegendentry{Smole\'nski \cite{Smoleski1999}}

      \node [above right] at (rel axis cs:0.5,0.52) {\Large $P_{cr}=1.075$};
        \addplot [style=very thick,dashed,color=magenta] table[x=ref_x
        ,y=ref_y, col sep=comma]{\loadFyL};

    \end{axis}
\end{tikzpicture}}
    }
     \caption{Lateral buckling of right-angle frame}
    \label{fig:Numerical example_4.1}
\end{figure}

In this benchmark example, which was investigated by Simo \& Vu-Quoc \cite{Simo.1986}, Crisfield \cite{M.A.Crisfield.1990}, and Smole\'nski \cite{Smoleski1999}  using beam elements, the stability analysis of a right-angle frame is investigated. The geometry and parameters of the frame can be seen in \cref{fig:Numerical example_4.1a}. For lateral buckling analysis, it is clamped at one end and subject to an in-plane load $P = 1.485$ and an out-of-plane perturbation load $P_{perb}=0.0002$ at the other end. The linear isotropic St.~Venant-Kirchhoff material model is considered with $E = 71240$, $\nu = 0.31$. 
The right angle frame is discretized by considering two cases of using only one single patch along the right angle frame, which has $C^0$-continuity at the kink, and using two patches connected at the kink with enforced $C^0$-continuity. In both cases, the right angle frame is discretized into a total of only 10 elements. Different degrees of NURBS functions are used for the discretization along the beam direction, while it is linear along the cross-section. 
\Cref{fig:Numerical example_4.1b} shows the undeformed and deformed configurations of the buckled right angle frame, while the load-displacement response curves of the frame at the center of the free end can be seen in \cref{fig:Numerical example_4.1c}. 
The response curves show that the results for the single- and multi-patch IGA solid-beam are almost identical, which serves as a verification of the multi-patch coupling. Furthermore, the convergence of the ISB$p$ results can be observed, with the load-displacement curves and critical load for ISB4 being very close to Simo \& Vu-Quoc \cite{Simo.1986} and Smole\'nski \cite{Smoleski1999}, who employed beam elements. 

\subsection{Homogenization of simple cubic lattice structure under shear}
\label{sec:num:5}

\begin{figure}[ht!]
    \centering
    \subfigure[ Simple cubic lattice structure \label{fig:geo_5}]{\resizebox{0.5\linewidth}{!}{
        \input{lattice_structure_model}
    }}
    \subfigure[ $L = 20$ ($a=1/41$) \label{fig:def_L_20}]{\resizebox{0.49\linewidth}{!}{
        \includegraphics[]{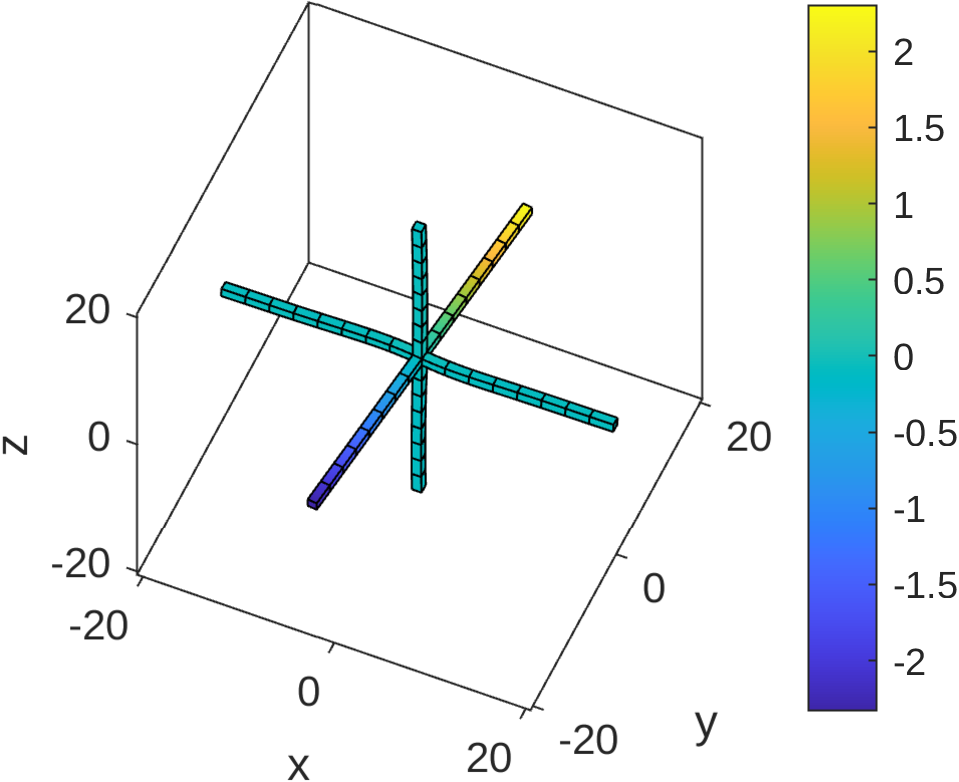}
    }}
    \subfigure[ $L = 10$ ($a=1/21$) \label{fig:def_L_10}]{\resizebox{0.49\linewidth}{!}{
        \includegraphics[]{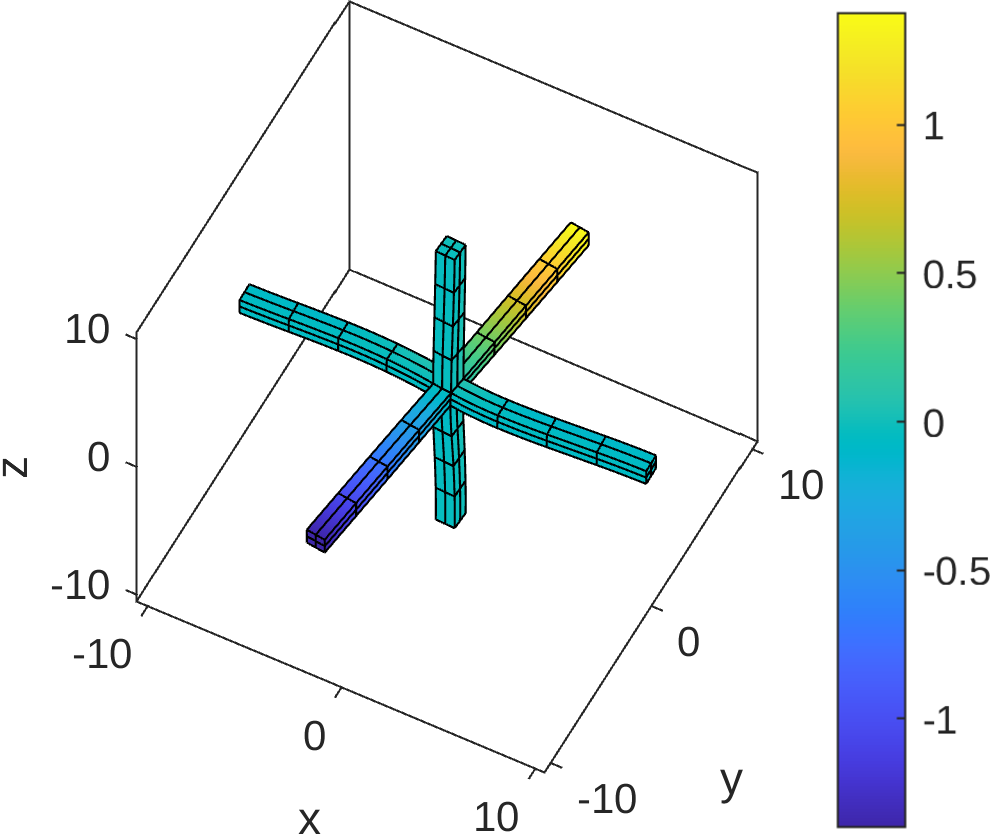}
    }}
    \subfigure[ $L = 5$ ($a=1/11$) \label{fig:def_L_5}]{\resizebox{0.49\linewidth}{!}{
        \includegraphics[]{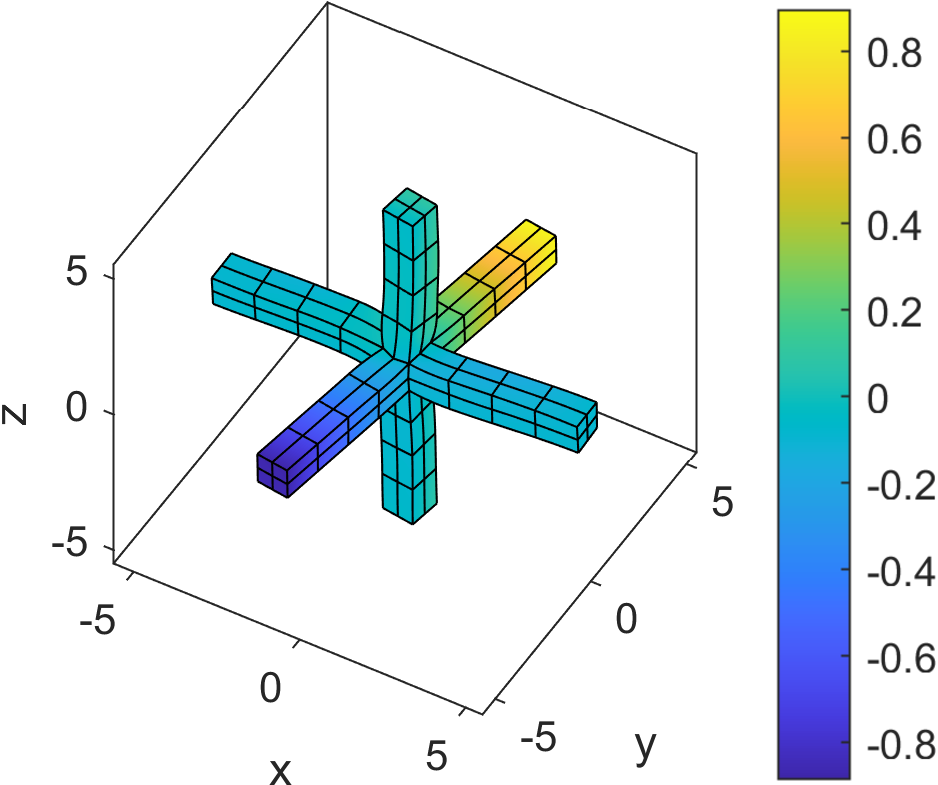}
    }}
    \caption{(a) Geometry and boundary conditions applied on SC unit cell. (b)-(d) Deformed configurations of SC lattice cells subject to simple shear in the ($x,y$) plane}
    \label{fig:geom_5.1}
\end{figure}

In the last example, the homogenization of simple cubic (SC) lattice structures is presented. In the SC unit cell, 6 patches of beams are connected together by a cubic solid element at the center of the lattice, as shown in \cref{fig:geo_5}.
A unit cross-section $H,~W = 1$ is assumed for each beam, but different lengths $L = 5,~10,~ \text{and} ~20$ are considered. Thus, the ``aspect ratios'' of beam thickness over total unit cell size are $a=1/(2L+1)=1/11,~1/21,~\text{and}~1/41$. Material properties of amorphous silicon are considered by using the isotropic Saint Venant-Kirchhoff material model with $E = 80~\text{GPa},~ \nu = 0.22$, as mentioned in \cite{DiLeo.2015}. 

This example aims to investigate the effective mechanical response of the periodic lattice structure under shear loading and compare the results with the ones obtained by using geometrically exact Timoshenko beam theory \cite{WEEGER2017100}.
For the numerical homogenization, the effective deformation gradient for the simple shear case is considered as
\begin{equation}
    \label{N.E.5.DG}
    \mathbf{F}^*=\begin{bmatrix}
1 & \lambda & 0 \\
0 & 1 & 0 \\
0 & 0 & 1  \\
\end{bmatrix}, \qquad 0\leq\lambda\leq 0.1.
\end{equation}
To avoid rigid body motions, the effective deformation is applied using a combination of Dirichlet (DBC) and periodic boundary conditions (PBC):
\begin{equation}
\label{N.E.5.BCs}
    \mathbf{u} = (\mathbf{F}^*-\mathbf{I})\cdot \mathbf{X},\qquad \quad
    \mathbf{u}^{+}-\mathbf{u}^{-} = (\mathbf{F}^*-\mathbf{I})\cdot (\mathbf{X}^{+}-\mathbf{X}^{-}),
\end{equation}
Here, $\mathbf{u}^{-}$ and $\mathbf{u}^{+}$ indicate the coupled displacements in the respective directions of the corresponding boundary nodes of the conformal mesh, see \cref{fig:geo_5}.
The resulting effective first Piola-Kirchhoff stress is defined by the mean value of the stresses with respect to the representative volume of the unit cell ${\Omega}_m$ of the microstructure:
\begin{equation}
    \label{N.E.5.PK}
     \mathbf{P}^* = \left< \mathbf{P} \right>:=\frac{1}{\left| {\Omega}_{m}\right|}\int_{{\Omega}_{m}}\mathbf{P} ~dV.
\end{equation}

In the thinnest case for $L=20$, each beam is discretized with $8\times 1\times 1$ ISB2 elements of degree $p=2$.
For $L=5 ~\text{and} ~ 10$, each beam is meshed by $4\times 2\times 2$ ISB2 elements with degree $p=2$ along the beam direction. Here, a finer discretization of the cross-section directions is applied due to the need to capture the substantial deformation of the cubic patch at the center.  

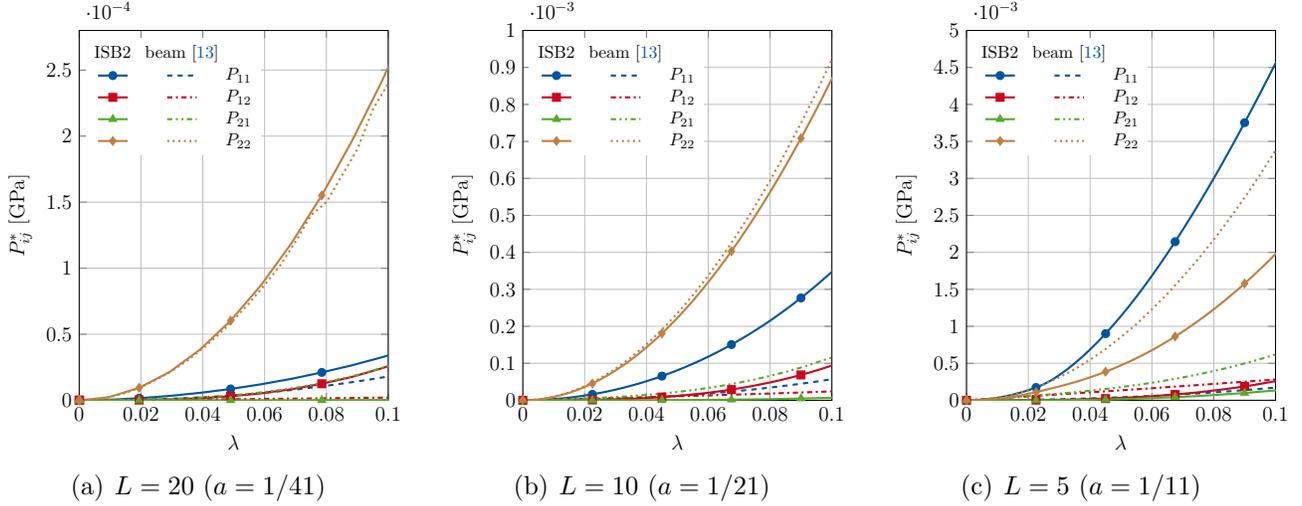
\begin{figure}[t!]
	\centering
        \subfigure[$L = 20$ ($a=1/41$)\label{fig:Numerical example_5.3c}]{\resizebox{0.32\textwidth}{!}{%
        \def\mystrut{\vphantom{hg}}

\pgfplotsset{
    legend image with text/.style={
        legend image code/.code={%
            \node[anchor=center] at (0.3cm,0cm) {#1};
        }
    },
    every axis plot/.append style={very thick},every mark/.append style={mark size=2pt}

}
\pgfplotstableread[col sep=comma]{comparison_L_5_10_20.csv}{\loadL}

\begin{tikzpicture}[tight background]
    \begin{axis}[legend columns=2,
            legend style={
                legend cell align=left,
            },
        width=19em,
        height=22em,
        grid,
        xmin=0,
        ymin=0,
        ymax = 0.00028,
        xmax = 0.1,
        xtick = {0,0.02,0.04,0.06,0.08,0.1},
        legend style={at={(0.02,0.98)},anchor=north west,{draw=none},font=\footnotesize},
        xlabel=\text{$\lambda$},
        ylabel=\text{$P^*_{ij} ~[\text{GPa}]$},
        xticklabel style={
        /pgf/number format/fixed,
        /pgf/number format/precision=2,
        }
                ]
                
        \addlegendimage{legend image with text=ISB2}
        \addlegendentry{}
        
        \addlegendimage{legend image with text=beam \cite{WEEGER2017100}}
        \addlegendentry{}

      \addplot[ mark=*, mark repeat=3, color=blue] table[x=val_sb_L_20_811_2, y=P11_sb_L_20_811_2, col sep=comma] {\loadL};
        \addlegendentry{}
        
        \addplot[mark=none, color=blue,dashed] table[x=val_beam_L_20, y=P11_L_20, col sep=comma] {\loadL};
        \addlegendentry{$P_{11}$}
        
        \addplot[mark=square*, mark repeat=3, color=red] table[x=val_sb_L_20_811_2, y=P12_sb_L_20_811_2, col sep=comma] {\loadL};
        \addlegendentry{}
        
        \addplot[mark=none, color=red,dashdotted] table[x=val_beam_L_20, y=P12_L_20, col sep=comma] {\loadL};
        \addlegendentry{$P_{12}$}
        
       \addplot[mark=triangle*, mark repeat=3, color=green] table[x=val_sb_L_20_811_2, y=P21_sb_L_20_811_2, col sep=comma] {\loadL};
        \addlegendentry{}
        
        \addplot[mark=none, color=green,dashdotdotted] table[x=val_beam_L_20, y=P21_L_20, col sep=comma] {\loadL};
        \addlegendentry{$P_{21}$}

       \addplot[mark=diamond*, mark repeat=3, color=brown] table[x=val_sb_L_20_811_2, y=P22_sb_L_20_811_2, col sep=comma] {\loadL};
        \addlegendentry{}
        
        \addplot[mark=none, color=brown,dotted] table[x=val_beam_L_20, y=P22_L_20, col sep=comma] {\loadL};
        \addlegendentry{$P_{22}$}
        
    \end{axis}
\end{tikzpicture}}}
        \subfigure[$L = 10$ ($a=1/21$)\label{fig:Numerical example_5.3b}]{\resizebox{0.32\textwidth}{!}{%
        \def\mystrut{\vphantom{hg}}

\pgfplotsset{
    legend image with text/.style={
        legend image code/.code={%
            \node[anchor=center] at (0.3cm,0cm) {#1};
        }
    },
    every axis plot/.append style={very thick},every mark/.append style={mark size=2.5pt}

}
\pgfplotstableread[col sep=comma]{comparison_L_5_10_20.csv}{\loadL}

\begin{tikzpicture}[tight background]
    \begin{axis}[legend columns=2,
            legend style={
                legend cell align=left,
            },
        width=19em,
        height=22em,
        grid,
        xmin=0,
        ymin=0,
        ymax = 0.001,
        xmax = 0.1,
        xtick = {0,0.02,0.04,0.06,0.08,0.1},
        ytick = {0,0.0001,...,0.001},
        legend style={at={(0.02,0.98)},anchor=north west,,{draw=none},font=\footnotesize},
        xlabel=\text{$\lambda$},
        ylabel=\text{$P^*_{ij} ~[\text{GPa}]$},
        xticklabel style={
        /pgf/number format/fixed,
        /pgf/number format/precision=2,
        }
                ]
                
        \addlegendimage{legend image with text=ISB2}
        \addlegendentry{}
        
        \addlegendimage{legend image with text=beam \cite{WEEGER2017100}}
        \addlegendentry{}

      \addplot[mark=*, mark repeat=9, color=blue] table[x=val_sb_L_10_422_2, y=P11_sb_L_10_422_2, col sep=comma] {\loadL};
        \addlegendentry{}
        
        \addplot[mark=none, color=blue,dashed] table[x=val_beam_L_10, y=P11_L_10, col sep=comma] {\loadL};
        \addlegendentry{$P_{11}$}
        
        \addplot[mark=square*, mark repeat=9, color=red] table[x=val_sb_L_10_422_2, y=P12_sb_L_10_422_2, col sep=comma] {\loadL};
        \addlegendentry{}
        
        \addplot[mark=none, color=red,dashdotted] table[x=val_beam_L_10, y=P12_L_10, col sep=comma] {\loadL};
        \addlegendentry{$P_{12}$}
        
       \addplot[mark=triangle*, mark repeat=9, color=green] table[x=val_sb_L_10_422_2, y=P21_sb_L_10_422_2, col sep=comma] {\loadL};
        \addlegendentry{}
        
        \addplot[mark=none, color=green,dashdotdotted] table[x=val_beam_L_10, y=P21_L_10, col sep=comma] {\loadL};
        \addlegendentry{$P_{21}$}

       \addplot[mark=diamond*, mark repeat=9, color=brown] table[x=val_sb_L_10_422_2, y=P22_sb_L_10_422_2, col sep=comma] {\loadL};
        \addlegendentry{}
        
        \addplot[mark=none, color=brown,dotted] table[x=val_beam_L_10, y=P22_L_10, col sep=comma] {\loadL};
        \addlegendentry{$P_{22}$}
        
    \end{axis}
\end{tikzpicture}}}
        \subfigure[ $L = 5$ ($a=1/11$) \label{fig:Numerical example_5.3a}]{\resizebox{0.32\textwidth}{!}{%
        \def\mystrut{\vphantom{hg}}

\pgfplotsset{
    legend image with text/.style={
        legend image code/.code={%
            \node[anchor=center] at (0.3cm,0cm) {#1};
        }
    },
    every axis plot/.append style={very thick},every mark/.append style={mark size=2.5pt}

}
\pgfplotstableread[col sep=comma]{comparison_L_5_10_20.csv}{\loadL}

\begin{tikzpicture}[tight background]
    \begin{axis}[legend columns=2,
            legend style={
                legend cell align=left,
            },
        width=19em,
        height=22em,
        grid,
        xmin=0,
        ymin=0,
        ymax = 0.005,
        ytick = {0,0.0005,...,0.005},
        xmax = 0.1,
        xtick = {0,0.02,0.04,0.06,0.08,0.1},
        legend style={at={(0.02,0.98)},anchor=north west,{draw=none},font=\footnotesize},
        xlabel=\text{$\lambda$},
        ylabel=\text{$P^*_{ij} ~[\text{GPa}]$},
        xticklabel style={
        /pgf/number format/fixed,
        /pgf/number format/precision=2,
        }
                ]
                
        \addlegendimage{legend image with text=ISB2}
        \addlegendentry{}
        
        \addlegendimage{legend image with text=beam \cite{WEEGER2017100}}
        \addlegendentry{}

      \addplot[mark=*, mark repeat=9, color=blue] table[x=val_sb_L_5_422_2, y=P11_sb_L_5_422_2, col sep=comma] {\loadL};
        \addlegendentry{}
        
        \addplot[mark=none, color=blue,dashed] table[x=val_beam_L_5, y=P11_L_5, col sep=comma] {\loadL};
        \addlegendentry{$P_{11}$}
        
        \addplot[mark=square*, mark repeat=9, color=red] table[x=val_sb_L_5_422_2, y=P12_sb_L_5_422_2, col sep=comma] {\loadL};
        \addlegendentry{}
        
        \addplot[mark=none, color=red,dashdotted] table[x=val_beam_L_5, y=P12_L_5, col sep=comma] {\loadL};
        \addlegendentry{$P_{12}$}
        
       \addplot[mark=triangle*, mark repeat=9, color=green] table[x=val_sb_L_5_422_2, y=P21_sb_L_5_422_2, col sep=comma] {\loadL};
        \addlegendentry{}
        
        \addplot[mark=none, color=green,dashdotdotted] table[x=val_beam_L_5, y=P21_L_5, col sep=comma] {\loadL};
        \addlegendentry{$P_{21}$}

       \addplot[mark=diamond*, mark repeat=9, color=brown] table[x=val_sb_L_5_422_2, y=P22_sb_L_5_422_2, col sep=comma] {\loadL};
        \addlegendentry{}
        
        \addplot[mark=none, color=brown,dotted] table[x=val_beam_L_5, y=P22_L_5, col sep=comma] {\loadL};
        \addlegendentry{$P_{22}$}
        
    \end{axis}
\end{tikzpicture}}}
	\caption{Effective, homogenized stresses $P_{ij}^*$ for SC unit cells with different aspect ratios}
	\label{fig:Numerical example_5.3}
\end{figure}

Deformed configurations for all 3 aspect ratios can be seen in \crefrange{fig:def_L_20}{fig:def_L_5} and the effective stresses are plotted over the load parameter $\lambda$ in \cref{fig:Numerical example_5.3}.
In the case of $L=20$, the results of the solid-beam elements are in close correspondence with the geometrically exact beam elements and the $P_{22}$-component is clearly the most dominant. 
However, for the thicker structures with $L=10$ and $L=5$, in the solid-beam results $P_{11}$ becomes larger. This feature is well recaptured by the solid-beam elements but cannot be observed in the beam results. This can be attributed to the cross-section deformations and the influence of the center element, both of which go beyond the capability of the beam theory.

This demonstrates the advantages of employing the isogeometric solid-beam formulation for lattice structures, where the beams are often relatively thick and nodal effects cannot be neglected.

\section{Conclusion}
\label{Sec:Conclusion}

In this contribution, we developed the robust isogeometric solid-beam element formulation which requires only displacement degrees of freedom for finite strain and hyperelastic material models. Locking effects are cured by the combination of the assumed natural strain method for membrane and transversal shear locking and the enhanced assumed strain method for Poisson thickness locking. To further enhance the robustness of the formulation, which is in particular necessary when using the EAS method, the mixed integration point method is embedded with the isogeometric solid-beam element, which uses B-splines and NURBS for the discretization of the displacement variables. Further advantages of the solid-beam element are that coupling with solid elements and the implementation of arbitrary three-dimensional material models are straightforward. 
The developed formulation is tested and verified in application to several non-linear benchmark problems and the results are quite satisfying and comparative with the isoparametric solid-beam element \cite{Frischkorn.2013} and beam elements based on the geometrically exact theory. As locking is alleviated, in particular for slender beams reliable and robust results are obtained compared to solid elements. Furthermore, in most cases, linear interpolation with one element along the cross-section direction is enough to give significant results. The formulation is also tested on a stability analysis problem using a single- and a two-patch geometry to get the buckling behaviour of a right-angle frame. In the end, homogenization of a simple cubic lattice structure is performed to obtain the effective stress behaviour of the meta-material.      

In the near future work, the isogeometric solid-beam element could be extended to multi-physics problems, e.g., thermo- or chemo-mechanical diffusion behaviour  to analyze  lattice-based lithium-ion electrodes. Furthermore, the isogeometric analysis approach enables a straightforward embedding in a shape optimization framework.


\section{Acknowledgements}
\label{Sec:ACKNOWLEDGMENTS}

This research was supported by the Deutsche Forschungsgemeinschaft (DFG -- German Research Foundation) -- Grant No.~460684687.
The authors want to thank Dr. Xianglong Peng and Dr. Shahed Rezeai from the Technical University of Darmstadt and Access e.V., respectively, for the fruitful discussions.


\appendix

\section{}

\subsection{Deformation dependent/ Non-Conservative loading}
\label{subsec: Deformation dependent}
The deformation-dependent load also referred to as follower load, changes with the configuration of the body. For instance, a pressure load needs to be normal to the applied surface, so as the surface moves the nodal external load also changes accordingly but the pressure remains constant. As external loads are dependent on the current configuration, the variational form, compare \cref{eqn:Governing equations.5}, can be written as
\begin{equation}
\label{eqn:Deformation dependent loading.1}
g_u^{ext}(\mathbf{u},\delta \mathbf{u}) =\int_{\Gamma_{t}}\delta \mathbf{u\cdot t}~dA = \int_{\Gamma _{t}}\delta \mathbf{u}\cdot(p~\mathbf{n})~dA~,
\end{equation}
where the normal vector $\mathbf{n}$ can be defined by the covariant base vectors of current configuration $\mathbf{g}_{a}=\frac{\partial \mathbf{x}}{\partial \xi_{a}}, ~~ \xi_{a}=[\xi, \eta, \zeta]$ as: 
\begin{equation}
\label{eqn:Deformation dependent loading.2}
\textbf{n}=\frac{\mathbf{g}_{\eta}\times \mathbf{g}_{\zeta}}{\left\| \mathbf{g}_{\eta}\times \mathbf{g}_{\zeta}\right\|},
\end{equation}

and the surface on the parametric domain $dA$ are given as:

\begin{equation}
\label{eqn:surface area}
dA = \left\|\mathbf{g}_{\eta}\times \mathbf{g}_{\zeta} \right\|d\eta d\zeta~.
\end{equation}

The discretized equation of the weak form is then expressed as:
\begin{equation}
\label{eqn:Deformation dependent loading.3}
\begin{aligned}
\int_{\Gamma _{t}}\delta \mathbf{u}~p~\mathbf{n}~dA &=\delta \mathbf{u}_{I}\int_{\eta_{1}^{e}}^{\eta_{2}^{e}}\int_{\zeta_{1}^{e}}^{\zeta_{2}^{e}}N_{I}~p(\eta ,\zeta )\left [ \mathbf{g}_{\eta} \times \mathbf{g}_{\zeta}\right ]d\eta d\zeta~ ,\\[4pt] 
\mathbf{f}_{I} &=\int_{\eta_{1}^{e}}^{\eta_{2}^{e}}\int_{\zeta_{1}^{e}}^{\zeta_{2}^{e}}N_{I}~p(\eta ,\zeta )\left [ \mathbf{g}_{\eta} \times \mathbf{g}_{\zeta}\right ]d\eta d\zeta~.
\end{aligned}
\end{equation}
Furthermore, the directional derivative of the weak form and its discretized form are
\begin{equation}
\label{eqn:Deformation dependent loading.4}
\begin{aligned}
Dg_u^{ext}(\mathbf{u},\delta \mathbf{u})\cdot \Delta \mathbf{u} &= \int_{\Gamma_t } \delta \mathbf{u}\cdot p(\eta,\zeta )(\Delta \mathbf{u}_{,\eta }\times \mathbf{g}_{\zeta}+\mathbf{g}_{\eta}\times \Delta \mathbf{u}_{,\zeta }) ~d\eta d\zeta ~,\\[4pt]
\mathbf{K}_{L}&= -\frac{\partial \mathbf{f}}{\partial \mathbf{u}}=\int_{\eta_{1}^{e}}^{\eta_{2}^{e}}\int_{\zeta_{1}^{e}}^{\zeta_{2}^{e}}N_{I}~p(\eta ,\zeta )\left [ N_{J,\eta }\mathbf{\hat{x}}_{,\zeta}-N_{J,\zeta } \mathbf{\hat{x}}_{,\eta }\right ] ~d\eta d\zeta ~,\\[4pt]
\text{with}\quad\mathbf{\hat{x}}&=\begin{bmatrix}
0 & -x_{3} & x_{2} \\
x_{3} & 0 & -x_{1} \\
-x_{2} & x_{1} & 0 \\
\end{bmatrix}, \quad \mathbf{x}=\mathbf{X}+\mathbf{u}.
\end{aligned}
\end{equation}

\subsection{Lagrange interpolation functions}
\label{subsec: Lagrange interpolation functions}

The one-dimensional Lagrangian interpolation functions $L_{i,q}({\xi})$ of degree $q \geq 1$ associated with the nodes $\{ {\xi}_1,\ldots,{\xi}_{q+1} \}$ for a given parametric coordinate ${\xi}$ are defined as: 
\begin{equation}
\label{eqn:Lagrange interpolation functions.1}
L_{i,q}({\xi}) = \prod_{\substack{j=1\\ j \neq i}}^{q+1} \frac{{\xi}-{\xi}_j}{{\xi}_i-{\xi}_j}, \quad 1 \leq i \leq q+1.
\end{equation}
Typically, for finite element shape functions, the nodes $\{ {\xi}_1,\ldots,{\xi}_{q+1} \}$ are equi-distributed in the reference interval $[-1,1]$. 
Importantly, the Lagrange polynomials satisfy the interpolation property 
\begin{equation}
\label{eqn:Lagrange interpolation functions.2}
L_{i,q}({\xi}_j)= \delta _{ij}, \quad 1\leq i,j\leq q+1.
\end{equation}

\bibliographystyle{elsarticle-num}
\addcontentsline{toc}{section}{References}

\begin{thebibliography}{10}
\expandafter\ifx\csname url\endcsname\relax
  \def\url#1{\texttt{#1}}\fi
\expandafter\ifx\csname urlprefix\endcsname\relax\def\urlprefix{URL }\fi
\expandafter\ifx\csname href\endcsname\relax
  \def\href#1#2{#2} \def\path#1{#1}\fi

\bibitem{Greer.2019}
J.~R. Greer, V.~S. Deshpande, Three-dimensional architected materials and
  structures: design, fabrication, and mechanical behavior, MRS Bulletin
  44~(10) (2019) 750--757.
\newblock \href {https://doi.org/10.1557/mrs.2019.232}
  {\path{doi:10.1557/mrs.2019.232}}.

\bibitem{Pang.2020}
Y.~Pang, Y.~Cao, Y.~Chu, M.~Liu, K.~Snyder, D.~MacKenzie, C.~Cao, Additive
  manufacturing of batteries, Advanced Functional Materials 30~(1) (2020)
  1906244.
\newblock \href {https://doi.org/10.1002/adfm.201906244}
  {\path{doi:10.1002/adfm.201906244}}.

\bibitem{Narita.2022}
K.~Narita, M.~A. Saccone, Y.~Sun, J.~R. Greer, Additive manufacturing of 3d
  batteries: a perspective, Journal of Materials Research 37~(9) (2022)
  1535--1546.
\newblock \href {https://doi.org/10.1557/s43578-022-00562-w}
  {\path{doi:10.1557/s43578-022-00562-w}}.

\bibitem{jamshidian2020}
M.~Jamshidian, N.~Boddeti, D.~W. Rosen, O.~Weeger, Multiscale modelling of soft
  lattice metamaterials: {{Micromechanical}} nonlinear buckling analysis,
  experimental verification, and macroscale constitutive behaviour,
  International Journal of Mechanical Sciences 188 (2020) 105956.
\newblock \href {https://doi.org/10.1016/j.ijmecsci.2020.105956}
  {\path{doi:10.1016/j.ijmecsci.2020.105956}}.

\bibitem{zhang2022a}
Y.~Zhang, K.~Yu, K.~H. Lee, K.~Li, H.~Du, Q.~Wang, Mechanics of stretchy
  elastomer lattices, Journal of the Mechanics and Physics of Solids 159 (2022)
  104782.
\newblock \href {https://doi.org/10.1016/j.jmps.2022.104782}
  {\path{doi:10.1016/j.jmps.2022.104782}}.

\bibitem{weeger2022}
O.~Weeger, I.~Valizadeh, Y.~Mistry, D.~Bhate, Inelastic finite deformation beam
  modeling, simulation, and validation of additively manufactured lattice
  structures, Additive Manufacturing Letters 4 (2022) 100111.
\newblock \href {https://doi.org/10.1016/j.addlet.2022.100111}
  {\path{doi:10.1016/j.addlet.2022.100111}}.

\bibitem{antman2005}
S.~Antman, Nonlinear {{Problems}} of {{Elasticity}}, Vol. 107 of Applied
  {{Mathematical Sciences}}, {Springer New York}, 2005.
\newblock \href {https://doi.org/10.1007/0-387-27649-1}
  {\path{doi:10.1007/0-387-27649-1}}.

\bibitem{eugster2015}
S.~Eugster, Geometric {{Continuum Mechanics}} and {{Induced Beam Theories}},
  Vol.~75 of Lecture {{Notes}} in {{Applied}} and {{Computational Mechanics}},
  {Springer International Publishing}, 2015.
\newblock \href {https://doi.org/10.1007/978-3-319-16495-3}
  {\path{doi:10.1007/978-3-319-16495-3}}.

\bibitem{Bathe.1979}
K.-J. Bathe, S.~Bolourchi, Large displacement analysis of three-dimensional
  beam structures, International Journal for Numerical Methods in Engineering
  14~(7) (1979) 961--986.
\newblock \href {https://doi.org/10.1002/nme.1620140703}
  {\path{doi:10.1002/nme.1620140703}}.

\bibitem{Simo.1986}
J.~C. Simo, L.~Vu-Quoc, A three-dimensional finite-strain rod model. part ii:
  Computational aspects, Computer Methods in Applied Mechanics and Engineering
  58~(1) (1986) 79--116.
\newblock \href {https://doi.org/10.1016/0045-7825(86)90079-4}
  {\path{doi:10.1016/0045-7825(86)90079-4}}.

\bibitem{jung2011}
P.~Jung, S.~Leyendecker, J.~Linn, M.~Ortiz, A discrete mechanics approach to
  the {{Cosserat}} rod theory—{{Part}} 1: Static equilibria, International
  Journal for Numerical Methods in Engineering 85~(1) (2011) 31--60.
\newblock \href {https://doi.org/10.1002/nme.2950}
  {\path{doi:10.1002/nme.2950}}.

\bibitem{meier2014}
C.~Meier, A.~Popp, W.~A. Wall, An objective {{3D}} large deformation finite
  element formulation for geometrically exact curved {{Kirchhoff}} rods,
  Computer Methods in Applied Mechanics and Engineering 278 (2014) 445--478.
\newblock \href {https://doi.org/10.1016/j.cma.2014.05.017}
  {\path{doi:10.1016/j.cma.2014.05.017}}.

\bibitem{WEEGER2017100}
O.~Weeger, S.-K. Yeung, M.~L. Dunn, Isogeometric collocation methods for
  cosserat rods and rod structures, Computer Methods in Applied Mechanics and
  Engineering 316 (2017) 100--122.
\newblock \href {https://doi.org/10.1016/j.cma.2016.05.009}
  {\path{doi:10.1016/j.cma.2016.05.009}}.

\bibitem{simo1984}
J.~C. Simo, K.~D. Hjelmstad, R.~L. Taylor, Numerical formulations of
  elasto-viscoplastic response of beams accounting for the effect of shear,
  Computer Methods in Applied Mechanics and Engineering 42~(3) (1984) 301--330.
\newblock \href {https://doi.org/10.1016/0045-7825(84)90011-2}
  {\path{doi:10.1016/0045-7825(84)90011-2}}.

\bibitem{gruttmann2000b}
F.~Gruttmann, R.~Sauer, W.~Wagner, Theory and numerics of three-dimensional
  beams with elastoplastic material behaviour, International Journal for
  Numerical Methods in Engineering 48~(12) (2000) 1675--1702.
\newblock \href
  {https://doi.org/10.1002/1097-0207(20000830)48:12<1675::AID-NME957>3.0.CO;2-6}
  {\path{doi:10.1002/1097-0207(20000830)48:12<1675::AID-NME957>3.0.CO;2-6}}.

\bibitem{lestringant2020}
C.~Lestringant, B.~Audoly, D.~M. Kochmann, A discrete, geometrically exact
  method for simulating nonlinear, elastic and inelastic beams, Computer
  Methods in Applied Mechanics and Engineering 361 (2020) 112741.
\newblock \href {https://doi.org/10.1016/j.cma.2019.112741}
  {\path{doi:10.1016/j.cma.2019.112741}}.

\bibitem{herrnbock2022}
L.~Herrnböck, A.~Kumar, P.~Steinmann, Two-scale off-and online approaches to
  geometrically exact elastoplastic rods, Computational Mechanics 71 (2022)
  1--24.
\newblock \href {https://doi.org/10.1007/s00466-022-02204-8}
  {\path{doi:10.1007/s00466-022-02204-8}}.

\bibitem{Weeger.2022b}
O.~Weeger, D.~Schillinger, R.~M{\"u}ller, Mixed isogeometric collocation for
  geometrically exact 3d beams with elasto-visco-plastic material behavior and
  softening effects, Computer Methods in Applied Mechanics and Engineering 399
  (2022) 115456.
\newblock \href {https://doi.org/10.1016/j.cma.2022.115456}
  {\path{doi:10.1016/j.cma.2022.115456}}.

\bibitem{favata2016}
A.~Favata, A beam theory consistent with three-dimensional thermo-elasticity,
  Mathematics and Mechanics of Solids 21~(4) (2016) 426--443.
\newblock \href {https://doi.org/10.1177/1081286514524974}
  {\path{doi:10.1177/1081286514524974}}.

\bibitem{ebrahimi2016}
F.~Ebrahimi, M.~R. Barati, Dynamic modeling of a thermo–piezo-electrically
  actuated nanosize beam subjected to a magnetic field, Applied Physics A
  122~(4) (2016) 451.
\newblock \href {https://doi.org/10.1007/s00339-016-0001-3}
  {\path{doi:10.1007/s00339-016-0001-3}}.

\bibitem{smriti2019}
Smriti, A.~Kumar, A.~Großmann, P.~Steinmann, A thermoelastoplastic theory for
  special {{Cosserat}} rods, Mathematics and Mechanics of Solids 24~(3) (2019)
  686--700.
\newblock \href {https://doi.org/10.1177/1081286517754132}
  {\path{doi:10.1177/1081286517754132}}.

\bibitem{Xia.2019}
X.~Xia, A.~Afshar, H.~Yang, C.~M. Portela, D.~M. Kochmann, C.~V. {Di Leo},
  J.~R. Greer, Electrochemically reconfigurable architected materials, Nature
  573~(7773) (2019) 205--213.
\newblock \href {https://doi.org/10.1038/s41586-019-1538-z}
  {\path{doi:10.1038/s41586-019-1538-z}}.

\bibitem{Stein.2014}
P.~Stein, B.~Xu, {3D} isogeometric analysis of intercalation-induced stresses
  in {Li}-ion battery electrode particles, Computer Methods in Applied
  Mechanics and Engineering 268 (2014) 225--244.
\newblock \href {https://doi.org/10.1016/j.cma.2013.09.011}
  {\path{doi:10.1016/j.cma.2013.09.011}}.

\bibitem{Zhao.2015}
Y.~Zhao, P.~Stein, B.-X. Xu, Isogeometric analysis of mechanically coupled
  cahn--hilliard phase segregation in hyperelastic electrodes of li-ion
  batteries, Computer Methods in Applied Mechanics and Engineering 297 (2015)
  325--347.
\newblock \href {https://doi.org/10.1016/j.cma.2015.09.008}
  {\path{doi:10.1016/j.cma.2015.09.008}}.

\bibitem{SUSSMAN1987357}
T.~Sussman, K.-J. Bathe, A finite element formulation for nonlinear
  incompressible elastic and inelastic analysis, Computers \& Structures 26~(1)
  (1987) 357--409.
\newblock \href {https://doi.org/10.1016/0045-7949(87)90265-3}
  {\path{doi:10.1016/0045-7949(87)90265-3}}.

\bibitem{Zienkiewicz1971ReducedIT}
O.~C. Zienkiewicz, R.~Taylor, J.~J.~M. Too, Reduced integration technique in
  general analysis of plates and shells, International Journal for Numerical
  Methods in Engineering 3 (1971) 275--290.
\newblock \href {https://doi.org/10.1002/NME.1620030211}
  {\path{doi:10.1002/NME.1620030211}}.

\bibitem{MALKUS197863}
D.~S. Malkus, T.~J. Hughes, Mixed finite element methods — reduced and
  selective integration techniques: A unification of concepts, Computer Methods
  in Applied Mechanics and Engineering 15~(1) (1978) 63--81.
\newblock \href {https://doi.org/10.1016/0045-7825(78)90005-1}
  {\path{doi:10.1016/0045-7825(78)90005-1}}.

\bibitem{BbarMethod}
T.~J.~R. Hughes, Generalization of selective integration procedures to
  anisotropic and nonlinear media, International Journal for Numerical Methods
  in Engineering 15~(9) (1980) 1413--1418.
\newblock \href {https://doi.org/10.1002/nme.1620150914}
  {\path{doi:10.1002/nme.1620150914}}.

\bibitem{Bathe1985}
K.-J. Bathe, E.~N. Dvorkin, A four-node plate bending element based on
  {M}indlin/{R}eissner plate theory and a mixed interpolation, International
  Journal for Numerical Methods in Engineering 21~(2) (1985) 367--383.
\newblock \href {https://doi.org/10.1002/nme.1620210213}
  {\path{doi:10.1002/nme.1620210213}}.

\bibitem{Simo1986}
J.~C. Simo, T.~J.~R. Hughes, On the variational foundations of assumed strain
  methods, Journal of Applied Mechanics 53~(1) (1986) 51--54.
\newblock \href {https://doi.org/10.1115/1.3171737}
  {\path{doi:10.1115/1.3171737}}.

\bibitem{Simo.1990}
J.~C. Simo, M.~S. Rifai, A class of mixed assumed strain methods and the method
  of incompatible modes, International Journal for Numerical Methods in
  Engineering 29~(8) (1990) 1595--1638.
\newblock \href {https://doi.org/10.1002/nme.1620290802}
  {\path{doi:10.1002/nme.1620290802}}.

\bibitem{Simo1992}
J.~C. Simo, F.~Armero, Geometrically non-linear enhanced strain mixed methods
  and the method of incompatible modes, International Journal for Numerical
  Methods in Engineering 33~(7) (1992) 1413--1449.
\newblock \href {https://doi.org/10.1002/nme.1620330705}
  {\path{doi:10.1002/nme.1620330705}}.

\bibitem{Hughes.2005}
T.~J.~R. Hughes, J.~A. Cottrell, Y.~Bazilevs, Isogeometric analysis: {CAD},
  finite elements, {NURBS}, exact geometry and mesh refinement, Computer
  Methods in Applied Mechanics and Engineering 194~(39-41) (2005) 4135--4195.
\newblock \href {https://doi.org/10.1016/j.cma.2004.10.008}
  {\path{doi:10.1016/j.cma.2004.10.008}}.

\bibitem{Cottrell.2007}
J.~A. Cottrell, T.~Hughes, A.~Reali, Studies of refinement and continuity in
  isogeometric structural analysis, Computer Methods in Applied Mechanics and
  Engineering 196~(41-44) (2007) 4160--4183.
\newblock \href {https://doi.org/10.1016/j.cma.2007.04.007}
  {\path{doi:10.1016/j.cma.2007.04.007}}.

\bibitem{Wall.2008}
W.~A. Wall, M.~A. Frenzel, C.~Cyron, Isogeometric structural shape
  optimization, Computer Methods in Applied Mechanics and Engineering
  197~(33-40) (2008) 2976--2988.
\newblock \href {https://doi.org/10.1016/j.cma.2008.01.025}
  {\path{doi:10.1016/j.cma.2008.01.025}}.

\bibitem{Nagy.2010}
A.~P. Nagy, M.~M. Abdalla, Z.~G{\"u}rdal, Isogeometric sizing and shape
  optimisation of beam structures, Computer Methods in Applied Mechanics and
  Engineering 199~(17-20) (2010) 1216--1230.
\newblock \href {https://doi.org/10.1016/j.cma.2009.12.010}
  {\path{doi:10.1016/j.cma.2009.12.010}}.

\bibitem{Echter.2010}
R.~Echter, M.~Bischoff, Numerical efficiency, locking and unlocking of {NURBS}
  finite elements, Computer Methods in Applied Mechanics and Engineering
  199~(5-8) (2010) 374--382.
\newblock \href {https://doi.org/10.1016/j.cma.2009.02.035}
  {\path{doi:10.1016/j.cma.2009.02.035}}.

\bibitem{DaBeiraoVeiga.2012}
L.~{Da Beir{\~a}o Veiga}, C.~Lovadina, A.~Reali, Avoiding shear locking for the
  {T}imoshenko beam problem via isogeometric collocation methods, Computer
  Methods in Applied Mechanics and Engineering 241-244 (2012) 38--51.
\newblock \href {https://doi.org/10.1016/j.cma.2012.05.020}
  {\path{doi:10.1016/j.cma.2012.05.020}}.

\bibitem{Adam.2014}
C.~Adam, S.~Bouabdallah, M.~Zarroug, H.~Maitournam, Improved numerical
  integration for locking treatment in isogeometric structural elements, {P}art
  {I}: Beams, Computer Methods in Applied Mechanics and Engineering 279 (2014)
  1--28.
\newblock \href {https://doi.org/10.1016/j.cma.2014.06.023}
  {\path{doi:10.1016/j.cma.2014.06.023}}.

\bibitem{Kiendl.2015}
J.~Kiendl, F.~Auricchio, T.~Hughes, A.~Reali, Single-variable formulations and
  isogeometric discretizations for shear deformable beams, Computer Methods in
  Applied Mechanics and Engineering 284 (2015) 988--1004.
\newblock \href {https://doi.org/10.1016/j.cma.2014.11.011}
  {\path{doi:10.1016/j.cma.2014.11.011}}.

\bibitem{Bouclier.2012}
R.~Bouclier, T.~Elguedj, A.~Combescure, Locking free isogeometric formulations
  of curved thick beams, Computer Methods in Applied Mechanics and Engineering
  245-246 (2012) 144--162.
\newblock \href {https://doi.org/10.1016/j.cma.2012.06.008}
  {\path{doi:10.1016/j.cma.2012.06.008}}.

\bibitem{Vo.2020}
D.~Vo, P.~Nanakorn, A total lagrangian {Timoshenko} beam formulation for
  geometrically nonlinear isogeometric analysis of planar curved beams, Acta
  Mechanica 231~(7) (2020) 2827--2847.
\newblock \href {https://doi.org/10.1007/s00707-020-02675-x}
  {\path{doi:10.1007/s00707-020-02675-x}}.

\bibitem{Auricchio.2013}
F.~Auricchio, L.~{Da Beir{\~a}o Veiga}, J.~Kiendl, C.~Lovadina, A.~Reali,
  Locking-free isogeometric collocation methods for spatial {Timoshenko} rods,
  Computer Methods in Applied Mechanics and Engineering 263 (2013) 113--126.
\newblock \href {https://doi.org/10.1016/j.cma.2013.03.009}
  {\path{doi:10.1016/j.cma.2013.03.009}}.

\bibitem{Elguedj.2008}
T.~Elguedj, Y.~Bazilevs, V.~M. Calo, T.~Hughes, \={B} and \={F} projection
  methods for nearly incompressible linear and non-linear elasticity and
  plasticity using higher-order {NURBS} elements, Computer Methods in Applied
  Mechanics and Engineering 197~(33-40) (2008) 2732--2762.
\newblock \href {https://doi.org/10.1016/j.cma.2008.01.012}
  {\path{doi:10.1016/j.cma.2008.01.012}}.

\bibitem{Echter.2013b}
R.~Echter, B.~Oesterle, M.~Bischoff, A hierarchic family of isogeometric shell
  finite elements, Computer Methods in Applied Mechanics and Engineering 254
  (2013) 170--180.
\newblock \href {https://doi.org/10.1016/j.cma.2012.10.018}
  {\path{doi:10.1016/j.cma.2012.10.018}}.

\bibitem{MI2021113693}
Y.~Mi, X.~Yu, Isogeometric {MITC} shell, Computer Methods in Applied Mechanics
  and Engineering 377 (2021) 113693.
\newblock \href {https://doi.org/10.1016/j.cma.2021.113693}
  {\path{doi:10.1016/j.cma.2021.113693}}.

\bibitem{Bouclier.2013b}
R.~Bouclier, T.~Elguedj, A.~Combescure, On the development of {NURBS}-based
  isogeometric solid shell elements: 2{D} problems and preliminary extension to
  3{D}, Computational Mechanics 52~(5) (2013) 1085--1112.
\newblock \href {https://doi.org/10.1007/s00466-013-0865-4}
  {\path{doi:10.1007/s00466-013-0865-4}}.

\bibitem{Bouclier.2013}
R.~Bouclier, T.~Elguedj, A.~Combescure, Efficient isogeometric {NURBS}-based
  solid-shell elements: Mixed formulation and \={B}-method, Computer Methods in
  Applied Mechanics and Engineering 267 (2013) 86--110.
\newblock \href {https://doi.org/10.1016/j.cma.2013.08.002}
  {\path{doi:10.1016/j.cma.2013.08.002}}.

\bibitem{Hosseini.2013}
S.~Hosseini, J.~J.~C. Remmers, C.~V. Verhoosel, R.~de~Borst, An isogeometric
  solid-like shell element for nonlinear analysis, International Journal for
  Numerical Methods in Engineering 95~(3) (2013) 238--256.
\newblock \href {https://doi.org/10.1002/nme.4505}
  {\path{doi:10.1002/nme.4505}}.

\bibitem{Hosseini.2014}
S.~Hosseini, J.~J. Remmers, C.~V. Verhoosel, R.~de~Borst, An isogeometric
  continuum shell element for non-linear analysis, Computer Methods in Applied
  Mechanics and Engineering 271 (2014) 1--22.
\newblock \href {https://doi.org/10.1016/j.cma.2013.11.023}
  {\path{doi:10.1016/j.cma.2013.11.023}}.

\bibitem{caseiro2014assumed}
J.~F. Caseiro, R.~A.~F. Valente, A.~Reali, J.~Kiendl, F.~Auricchio, R.~J.
  {Alves de Sousa}, On the assumed natural strain method to alleviate locking
  in solid-shell nurbs-based finite elements, Computational Mechanics 53~(6)
  (2014) 1341--1353.
\newblock \href {https://doi.org/10.1007/s00466-014-0978-4}
  {\path{doi:10.1007/s00466-014-0978-4}}.

\bibitem{caseiro2015assumed}
J.~F. Caseiro, R.~A.~F. Valente, A.~Reali, J.~Kiendl, F.~Auricchio, R.~J.
  Alves~de Sousa, Assumed natural strain {NURBS}-based solid-shell element for
  the analysis of large deformation elasto-plastic thin-shell structures,
  Computer Methods in Applied Mechanics and Engineering 284 (2015) 861--880.
\newblock \href {https://doi.org/10.1016/j.cma.2014.10.037}
  {\path{doi:10.1016/j.cma.2014.10.037}}.

\bibitem{Bouclier.2015}
R.~Bouclier, T.~Elguedj, A.~Combescure, An isogeometric locking-free
  {NURBS}-based solid-shell element for geometrically nonlinear analysis,
  International Journal for Numerical Methods in Engineering 101~(10) (2015)
  774--808.
\newblock \href {https://doi.org/10.1002/nme.4834}
  {\path{doi:10.1002/nme.4834}}.

\bibitem{LEONETTI2018159}
L.~Leonetti, F.~Liguori, D.~Magisano, G.~Garcea, An efficient isogeometric
  solid-shell formulation for geometrically nonlinear analysis of elastic
  shells, Computer Methods in Applied Mechanics and Engineering 331 (2018)
  159--183.
\newblock \href {https://doi.org/10.1016/j.cma.2017.11.025}
  {\path{doi:10.1016/j.cma.2017.11.025}}.

\bibitem{Antolin.2020}
P.~Antolin, J.~Kiendl, M.~Pingaro, A.~Reali, A simple and effective method
  based on strain projections to alleviate locking in isogeometric solid
  shells, Computational Mechanics 65~(6) (2020) 1621--1631.
\newblock \href {https://doi.org/10.1007/s00466-020-01837-x}
  {\path{doi:10.1007/s00466-020-01837-x}}.

\bibitem{Frischkorn.2013}
J.~Frischkorn, S.~Reese, A solid-beam finite element and non-linear
  constitutive modelling, Computer Methods in Applied Mechanics and Engineering
  265 (2013) 195--212.
\newblock \href {https://doi.org/10.1016/j.cma.2013.06.009}
  {\path{doi:10.1016/j.cma.2013.06.009}}.

\bibitem{Schwarze.2011b}
M.~Schwarze, S.~Reese, A reduced integration solid-shell finite element based
  on the {EAS} and the {ANS} concept-large deformation problems, International
  Journal for Numerical Methods in Engineering 85~(3) (2011) 289--329.
\newblock \href {https://doi.org/10.1002/nme.2966}
  {\path{doi:10.1002/nme.2966}}.

\bibitem{Choi.2021c}
M.-J. Choi, R.~A. Sauer, S.~Klinkel, An isogeometric finite element formulation
  for geometrically exact {Timoshenko} beams with extensible directors,
  Computer Methods in Applied Mechanics and Engineering 385 (2021) 113993.
\newblock \href {https://doi.org/10.1016/j.cma.2021.113993}
  {\path{doi:10.1016/j.cma.2021.113993}}.

\bibitem{Choi.2023}
M.-J. Choi, R.~A. Sauer, S.~Klinkel, A selectively reduced degree basis for
  efficient mixed nonlinear isogeometric beam formulations with extensible
  directors, arXiv preprint (2023).
\newblock \href {https://doi.org/10.48550/arXiv.2306.13354}
  {\path{doi:10.48550/arXiv.2306.13354}}.

\bibitem{Magisano.2017}
D.~Magisano, L.~Leonetti, G.~Garcea, How to improve efficiency and robustness
  of the {Newton} method in geometrically non-linear structural problem
  discretized via displacement-based finite elements, Computer Methods in
  Applied Mechanics and Engineering 313 (2017) 986--1005.
\newblock \href {https://doi.org/10.1016/j.cma.2016.10.023}
  {\path{doi:10.1016/j.cma.2016.10.023}}.

\bibitem{Pfefferkorn.2021}
R.~Pfefferkorn, S.~Bieber, B.~Oesterle, M.~Bischoff, P.~Betsch, Improving
  efficiency and robustness of enhanced assumed strain elements for nonlinear
  problems, International Journal for Numerical Methods in Engineering 122~(8)
  (2021) 1911--1939.
\newblock \href {https://doi.org/10.1002/nme.6605}
  {\path{doi:10.1002/nme.6605}}.

\bibitem{Leonetti.2023}
L.~Leonetti, J.~Kiendl, A mixed integration point (mip) formulation for
  hyperelastic kirchhoff--love shells for nonlinear static and dynamic
  analysis, Computer Methods in Applied Mechanics and Engineering 416 (2023)
  116325.
\newblock \href {https://doi.org/10.1016/j.cma.2023.116325}
  {\path{doi:10.1016/j.cma.2023.116325}}.

\bibitem{PiegTill96}
L.~Piegl, W.~Tiller, The {NURBS} Book, 2nd Edition, Springer-Verlag, New York,
  NY, USA, 1996.
\newblock \href {https://doi.org/10.1007/978-3-642-97385-7}
  {\path{doi:10.1007/978-3-642-97385-7}}.

\bibitem{Cottrell.2009}
J.~A. Cottrell, T.~J.~R. Hughes, J.~J. Bazilevs, Isogeometric analysis: Toward
  integration of {CAD} and {FEA}, Wiley, Chichester, 2009.
\newblock \href {https://doi.org/10.1002/9780470749081}
  {\path{doi:10.1002/9780470749081}}.

\bibitem{Du.2020}
X.~Du, G.~Zhao, W.~Wang, M.~Guo, R.~Zhang, J.~Yang, {NLIGA}: A {MATLAB}
  framework for nonlinear isogeometric analysis, Computer Aided Geometric
  Design 80 (2020) 101869.
\newblock \href {https://doi.org/10.1016/j.cagd.2020.101869}
  {\path{doi:10.1016/j.cagd.2020.101869}}.

\bibitem{E.Ramm1976}
E.~Ramm, A plate/shell element for large deflections and rotations, MIT Press,
  Cambridge, MA, Symp. on formulations and computational algorithms in finite
  element analysis (1976).

\bibitem{Rhim.1998}
J.~Rhim, S.~Lee, A vectorial approach to computational modelling of beams
  undergoing finite rotations, International Journal for Numerical Methods in
  Engineering 41~(3) (1998) 527--540.
\newblock \href
  {https://doi.org/10.1002/(SICI)1097-0207(19980215)41:3<527::AID-NME297>3.0.CO;2-7}
  {\path{doi:10.1002/(SICI)1097-0207(19980215)41:3<527::AID-NME297>3.0.CO;2-7}}.

\bibitem{M.A.Crisfield.1990}
M.~A. Crisfield, A consistent co-rotational formulation for non-linear,
  three-dimensional, beam-elements, Computer Methods in Applied Mechanics and
  Engineering 81~(2) (1990) 131--150.
\newblock \href {https://doi.org/10.1016/0045-7825(90)90106-V}
  {\path{doi:10.1016/0045-7825(90)90106-V}}.

\bibitem{Cardona.1988}
A.~Cardona, M.~Geradin, A beam finite element non-linear theory with finite
  rotations, International Journal for Numerical Methods in Engineering 26~(11)
  (1988) 2403--2438.
\newblock \href {https://doi.org/10.1002/nme.1620261105}
  {\path{doi:10.1002/nme.1620261105}}.

\bibitem{Dvorkin.1988}
E.~N. Dvorkin, E.~Onte, J.~Oliver, On a non-linear formulation for curved
  {Timoshenko} beam elements considering large displacement/rotation
  increments, International Journal for Numerical Methods in Engineering 26~(7)
  (1988) 1597--1613.
\newblock \href {https://doi.org/10.1002/nme.1620260710}
  {\path{doi:10.1002/nme.1620260710}}.

\bibitem{Smoleski1999}
W.~Smoleński, Statically and kinematically exact nonlinear theory of rods and
  its numerical verification, Computer Methods in Applied Mechanics and
  Engineering 178 (1999) 89--113.
\newblock \href {https://doi.org/10.1016/S0045-7825(99)00006-7}
  {\path{doi:10.1016/S0045-7825(99)00006-7}}.

\bibitem{DiLeo.2015}
C.~V. {Di Leo}, E.~Rejovitzky, L.~Anand, Diffusion--deformation theory for
  amorphous silicon anodes: The role of plastic deformation on electrochemical
  performance, International Journal of Solids and Structures 67-68 (2015)
  283--296.
\newblock \href {https://doi.org/10.1016/j.ijsolstr.2015.04.028}
  {\path{doi:10.1016/j.ijsolstr.2015.04.028}}.

\end{thebibliography}

\newpage

\end{document}